\documentclass[11pt]{amsart}
\usepackage[lmargin=1in,rmargin=1in,tmargin=1in,bmargin=1in]{geometry}
\usepackage[unicode,bookmarks=false]{hyperref}
\usepackage{amsfonts,amssymb,amsmath,amsthm,mathrsfs}
\usepackage{array,esint}

\numberwithin{equation}{section}
\allowdisplaybreaks[4]
\newtheorem{lemma}{Lemma}[section]
\newtheorem{theorem}{Theorem}[section]
\newtheorem{corollary}[lemma]{Corollary}

\newcommand{\Z}{\mathbb{Z}}
\newcommand{\R}{\mathbb{R}}
\newcommand{\C}{\mathbb{C}}
\newcommand{\N}{\mathbb{N}}
\newcommand{\e}{\operatorname{e}}
\newcommand{\sumtwo}{\operatorname*{\sum\sum}}
\newcommand{\quadand}{\quad \textnormal{and} \quad}

\newcommand{\E}[1]{\mathbb{E}[#1]}
\newcommand{\Pb}[1]{\mathbb{P}[#1]}
\newcommand{\eps}{\varepsilon}
\renewcommand{\Re}{\operatorname{Re}}
\newcommand{\sumfour}{\operatorname*{\sum\sum\sum\sum}}
\newcommand{\nn}[1]{\lVert #1 \rVert_2}
\newcommand{\Bcal}{\mathcal{B}}
\newcommand{\Cw}{\mathcal{C}_w}
\newcommand{\Dw}{\mathcal{D}_w}
\newcommand{\Hcal}{\mathcal{H}}
\newcommand{\Jcal}{\mathcal{J}}
\newcommand{\Qcal}{\mathcal{Q}}
\newcommand{\Acal}{\mathcal{A}}
\newcommand{\Pcal}{\mathcal{P}}
\newcommand{\ind}{\mathbf{1}}
\begin{document}

\title[Bilinear forms with Kloosterman fractions]{Bilinear forms with Kloosterman fractions and applications}

\author[A.~Dong]{Anji Dong}
\address{Department of Mathematics, University of Illinois, 1409 West Green Street, Urbana, IL 61801, USA}
\email{anjid2@illinois.edu}

\author[N.~Robles]{Nicolas Robles}
\address{RAND Corporation, Engineering and Applied Sciences, 1200 S Hayes St, Arlington, VA, 22202, USA}
\email{robles.nicolas.m@gmail.com}

\author[D.~Zeindler]{Dirk Zeindler}
\address{Department of Mathematics and Statistics, Lancaster University, Fylde College, Bailrigg, Lancaster LA1 4YF, United Kingdom}
\email{d.zeindler@lancaster.ac.uk}

\begin{abstract}
We establish bounds for bilinear forms $\sumtwo_{m,n}\alpha_m\beta_n\,\e(a\overline{m}/n)$ with Kloosterman fractions, where $M<m\le 2M$, $N<n\le 2N$ and $(m,n)=1$, for arbitrary sequences of complex numbers $(\alpha_m)$ and $(\beta_n)$ and for positive integer numerators $a\ll X^{1/10}$, where $X=\max(M,N)$, and for the corresponding trilinear forms with an average over $a$. In the balanced case, where $M$ and $N$ are of the same order of magnitude, and for bounded $a$, the saving over the trivial bound is $X^{1/16-\eps}$, which improves upon the savings $X^{1/48-\eps}$ of Duke, Friedlander and Iwaniec and $X^{1/20-\eps}$ of Bettin and Chandee. The proof follows the amplification method of these authors; where they apply Weil's bound, we use for most of the terms the estimates of Deshouillers and Iwaniec for sums of Kloosterman sums, in the form given by Pascadi. As an application, we prove an asymptotic formula for the twisted second moment of the Riemann zeta-function with Dirichlet polynomials of length $T^{\theta}$ for $\theta<14/27$, extending the range $\theta<17/33$ of Bettin, Chandee and Radziwi{\l}{\l}. We also give applications to a determinant equation and to convolutions in arithmetic progressions with a fixed residue, on average over moduli as large as $x^{14/27-\eps}$.
\end{abstract}

\maketitle

\tableofcontents

\section{Introduction and results}
\label{sec:intro}

\subsection{Previous results in the literature}
\label{sec:previous_results}
Let $\alpha=(\alpha_m)$, $\beta=(\beta_n)$ and $\nu=(\nu_a)$ be sequences of complex numbers supported on $m\sim M$, $n\sim N$ and $a\sim A$, respectively, and set
\begin{align}
	\Bcal(A,M,N)
	:=
	\sum_{a\sim A}
	\sumtwo_{\substack{m\sim M,\, n\sim N\\(m,n)=1}}
	\nu_a\alpha_m\beta_n\,\e\Big(\frac{a\overline{m}}{n}\Big),
	\label{eq:def_B_AMN}
\end{align}
where $\overline{m}$ denotes the multiplicative inverse of $m$ modulo $n$.
The numbers $a\overline{m}/n$ are called Kloosterman fractions. The form \eqref{eq:def_B_AMN} is trilinear in $\nu$, $\alpha$ and $\beta$, and it is a bilinear form in $\alpha$ and $\beta$ when $\nu$ is supported at a single integer $a$. The notation is explained in Section~\ref{sec:notation}. We write $\nn{\alpha}=(\sum_{m\sim M}|\alpha_m|^2)^{1/2}$ for the $\ell^2$-norm of $\alpha$, and similarly for $\beta$ and $\nu$. For $A\ge1/2$ and $M,N\ge1$ the Cauchy-Schwarz inequality gives the trivial bound
\[
	|\Bcal(A,M,N)|\ll\nn{\alpha}\nn{\beta}\nn{\nu}\,(AMN)^{1/2},
\]
and one asks for a saving of a power of the lengths when the coefficients are arbitrary. Such estimates are a substitute for bounds for individual incomplete Kloosterman sums over ranges that are too short for the known methods, see the introduction of \cite{DFI}.

The first estimates of this kind are due to Duke, Friedlander and Iwaniec \cite{DFI}. For a positive integer $a$ they proved in \cite[Theorem~2]{DFI} that
\[
	\sumtwo_{\substack{m\sim M,\, n\sim N\\(m,n)=1}}
	\alpha_m\beta_n\,\e\Big(\frac{a\overline{m}}{n}\Big)
	\ll_\eps
	\nn{\alpha}\nn{\beta}\,(a+MN)^{3/8}(M+N)^{11/48+\eps}.
\]
For $a\ll MN$ and $M$, $N$ of the same order of magnitude $X$, this saves a factor $X^{1/48-\eps}$ over the trivial bound. A second estimate, \cite[Theorem~1]{DFI}, is nontrivial when $M$ and $N$ are neither nearly equal nor of very different size.

Bettin and Chandee \cite{BC} refined the argument of \cite{DFI} and introduced an average over the numerator, which is present in many applications. With a slightly different normalization of the dyadic ranges, \cite[Theorem~1]{BC} states that for every nonzero integer $\vartheta$,
\begin{align*}
	&\sum_{a\sim A}
	\sumtwo_{\substack{m\sim M,\, n\sim N\\(m,n)=1}}
	\nu_a\alpha_m\beta_n\,\e\Big(\frac{\vartheta a\overline{m}}{n}\Big)\\
	&\qquad\ll_\eps
	\nn{\alpha}\nn{\beta}\nn{\nu}\Big(1+\frac{|\vartheta|A}{MN}\Big)^{1/2}
	\big((AMN)^{7/20+\eps}(M+N)^{1/4}+(AMN)^{3/8+\eps}(AM+AN)^{1/8}\big).
\end{align*}
For a single numerator $a\ll MN$ and $M$, $N$ of the same order of magnitude $X$, the saving over the trivial bound is $X^{1/20-\eps}$.

These estimates have been applied to several problems, among them those taken up in Section~\ref{sec:new_results}, where the results that we use are recalled. In \cite{DFI2} the estimates of \cite{DFI} are applied to the number of solutions of determinant equations and to mean values of Dirichlet $L$-functions and Dirichlet polynomials, and the result on determinant equations was improved in \cite[Corollary~1]{BC}. The estimate of \cite{BC} was used in \cite{BCR} for the mean square of the Riemann zeta-function multiplied by a Dirichlet polynomial, in \cite{BPRZ} for the corresponding problem for Dirichlet $L$-functions to a prime modulus, and in \cite{FR,FR2}, combined with the dispersion method, for the distribution of convolutions in arithmetic progressions. Wright \cite{Wr1,Wr2} enlarged the ranges of \cite{FR,FR2} by means of variants of the estimate of \cite{BC} for moduli with a fixed factor and for coefficients supported on short intervals.

In \cite{DFI} and \cite{BC} the sum is set up by the Cauchy-Schwarz inequality and an amplification argument, and the incomplete Kloosterman sums that result are estimated by Weil's bound \cite{We}. Duke, Friedlander and Iwaniec remark \cite[p.~35]{DFI} that one could attempt to use at this point the spectral theory of automorphic forms, in the form of the estimates of Deshouillers and Iwaniec \cite{DI}, and that the large level of the group involved makes the quality of the resulting bound unclear. For sums of Kloosterman fractions with smooth weights, estimates of this type are given in \cite[Theorem~12]{DI}. They have since been refined in the exceptional spectrum, see \cite{T} and \cite{P}.

\subsection{New results}
\label{sec:new_results}
Our main result is the following trilinear estimate.
\begin{theorem}
\label{thm:main}
Let $M,N\ge 2$, set $X=\max(M,N)$, and let $1/2\le A\le X^{1/10}$. Let $\alpha$, $\beta$ and $\nu$ be sequences of complex numbers supported on $m\sim M$, $n\sim N$ and $a\sim A$, respectively. Then for every $\eps>0$,
\begin{align}
	|\Bcal(A,M,N)|
	\ll_\eps
	\nn{\alpha}\nn{\beta}\nn{\nu}\,
	X^{15/16+\eps}A^{3/8}.
	\label{eq:main_bound}
\end{align}
\end{theorem}
Since $a\sim A$ means $A<a\le 2A$, the choice $A=1/2$ restricts the sum over $a$ to the single numerator $a=1$. A sequence $\nu$ supported at a single integer $a$ gives back the bilinear forms of \cite{DFI} and \cite{BC}.
\begin{corollary}
\label{cor:bilinear}
Let $M,N\ge 2$, $X=\max(M,N)$, and let $a$ be a positive integer with $a\le X^{1/10}$. Let $\alpha$ and $\beta$ be sequences of complex numbers supported on $m\sim M$ and $n\sim N$, respectively. Then for every $\eps>0$,
\begin{align}
	\Bcal_a(M,N)
	:=
	\sumtwo_{\substack{m\sim M,\, n\sim N\\(m,n)=1}}
	\alpha_m\beta_n\,\e\Big(\frac{a\overline{m}}{n}\Big)
	\ll_\eps
	\nn{\alpha}\nn{\beta}\,
	X^{15/16+\eps}a^{3/8}.
	\label{eq:bilinear_bound}
\end{align}
\end{corollary}
For a single numerator $a\ll MN$, Theorem~1 of \cite{BC} implies the uniform bound $X^{19/20+\eps}$. The powers of $X$ in \eqref{eq:bilinear_bound} and in this bound agree at $a=X^{1/30}$. For $a\le X^{1/30-\delta}$, with fixed $\delta>0$, \eqref{eq:bilinear_bound} gives a power saving over it. A bound for larger numerators, $X^{1/10}\le A\le X^{1/4-\eta}$, is given in Section~\ref{sec:larger_numerators}.

We give four applications, in which Theorem~\ref{thm:main} and Corollary~\ref{cor:bilinear} are used in place of the estimates of \cite{BC}. They are proved in Section~\ref{sec:moment}.

The first application concerns the twisted second moment of the Riemann zeta-function. In its statement, and in Section~\ref{sec:app_zeta}, $\alpha$ and $\beta$ denote complex numbers and not sequences. Let $T\ge2$, let $\Phi$ be a smooth function supported in $[1,2]$ with $\Phi^{(j)}(x)\ll_j(\log T)^j$ for all $j\ge0$, and let
\[
	P(s)=\sum_{n\le Y}\frac{c_n}{n^s},
	\qquad c_n\ll_\eps n^\eps,
	\qquad Y=T^\theta,
\]
with $0<\theta<1$. For $\alpha,\beta\ll(\log T)^{-1}$ we consider
\begin{align}
	I(\alpha,\beta)
	:=
	\int_\R\zeta(\tfrac12+\alpha+it)\,\zeta(\tfrac12+\beta-it)\,|P(\tfrac12+it)|^2\,\Phi\Big(\frac tT\Big)\,dt.
	\label{eq:zeta_I}
\end{align}
The implied constants in Theorem~\ref{thm:zeta} and in its proof may depend on those in the bounds for $c_n$, $\Phi^{(j)}$ and $\alpha$, $\beta$, and $T$ is assumed to be sufficiently large in terms of the last of these. For $\alpha=\beta=0$ an asymptotic formula for $I(\alpha,\beta)$ was proved in \cite{BCH} for $\theta<1/2$ and in \cite{BCR} for $\theta<17/33$, the latter by means of the estimate of \cite{BC}. The shifts $\alpha$, $\beta$ were included in \cite{PR}. Using \cite{BCR} and \cite{PR} with Theorem~\ref{thm:main} in place of the estimate of \cite{BC}, we obtain the following.

\begin{theorem}
\label{thm:zeta}
Let $0<\theta\le14/27$ and $\alpha,\beta\ll(\log T)^{-1}$. Then for every $\eps>0$,
\begin{align}
	I(\alpha,\beta)
	={}&
	\sum_{m,n\le Y}\frac{c_m\overline{c_n}}{[m,n]}\,\frac{(m,n)^{\alpha+\beta}}{m^\beta n^\alpha}
	\int_\R\bigg(\zeta(1+\alpha+\beta)+\zeta(1-\alpha-\beta)\Big(\frac{2\pi mn}{t(m,n)^2}\Big)^{\alpha+\beta}\bigg)\Phi\Big(\frac tT\Big)\,dt
	\notag\\
	&+O_\eps\big(T^{1/8+\eps}\,Y^{27/16}+T^{1/2+\eps}\,Y^{1/2}\big),
	\label{eq:zeta_result}
\end{align}
where for $\alpha+\beta=0$ the expression in the large brackets is defined by continuity.
\end{theorem}

The error term in \eqref{eq:zeta_result} is $O(T^{1-\delta})$ with some $\delta=\delta(\theta)>0$ for every $\theta<\frac{14}{27}=\frac{1}{2}+\frac{1}{54}$, whereas \cite{BCR} and \cite{PR} have $\theta<\frac{17}{33}=\frac{1}{2}+\frac{1}{66}$. In \cite[Theorem~1.1]{PR} the factor $m^\beta n^\alpha$ is printed as $m^\alpha n^\beta$. The placement in \eqref{eq:zeta_result} is the one that corresponds to \eqref{eq:zeta_I}, and the resulting double sums agree when the coefficients $c_n$ are real.

The second application is to the determinant equation $m_1n_2-m_2n_1=\Delta$ in the setting of \cite[Corollary~1]{BC}: the variables $m_1,m_2$ carry smooth weights and the variables $n_1,n_2$ carry arbitrary coefficients. Here $\Delta$ denotes a nonzero integer, and $f,g:\R\to\C$ are fixed smooth functions with compact support in $(1,2)$. For $N\ge2$ and sequences $\alpha$, $\beta$ of complex numbers supported on $n\sim N$, set
\begin{align}
	\mathcal S_\Delta(N)
	&:=
	\sumfour_{\substack{m_1,m_2\in\Z,\ n_1,n_2\sim N\\ m_1n_2-m_2n_1=\Delta}}
	f\Big(\frac{m_1}{N}\Big)g\Big(\frac{m_2}{N}\Big)\alpha_{n_1}\beta_{n_2}
	\label{eq:det_S}
	\\
	\mathcal R_\Delta(N)
	&:=
	\sumtwo_{\substack{n_1,n_2\sim N\\(n_1,n_2)\mid\Delta}}
	\frac{(n_1,n_2)}{n_1n_2}\,\alpha_{n_1}\beta_{n_2}
	\int_\R f\Big(\frac{x+\Delta}{Nn_2}\Big)g\Big(\frac{x}{Nn_1}\Big)\,dx.
	\label{eq:det_R}
\end{align}

\begin{theorem}
\label{thm:det}
Let $N\ge2$ and let $\Delta$ be a nonzero integer. Then for every $\eps>0$,
\begin{align}
	\mathcal S_\Delta(N)
	=
	\mathcal R_\Delta(N)
	+O_{f,g,\eps}\big(\nn{\alpha}\nn{\beta}\,N^{\eps}\min\{|\Delta|^{3/8}N^{15/16},\,N^{19/20}\}\big).
	\label{eq:det_result}
\end{align}
\end{theorem}

When all four variables have size $N$ and the smooth weights are fixed, \cite[Corollary~1]{BC} has the same main term and the error term $O(\nn{\alpha}\nn{\beta}N^{19/20+\eps})$, uniformly in $\Delta$. This is the second term of the minimum in \eqref{eq:det_result}, and it is the smaller one for $|\Delta|\ge N^{1/30}$. The first term is new, and for fixed $\Delta$ it gives the error term $O(\nn{\alpha}\nn{\beta}N^{15/16+\eps})$.

The last two applications concern convolutions in arithmetic progressions. In their statements, and in Section~\ref{sec:app_conv}, $a$ denotes a nonzero integer, and $\tau_k(n)$ is the number of ways of writing $n$ as a product of $k$ positive integers. Following \cite{FR}, we say that a sequence $\beta=(\beta_n)$ of complex numbers satisfies the Siegel-Walfisz condition if there is an integer $k\ge1$ such that for every $B>0$, uniformly in $y\ge2$, $q>|c|\ge1$ with $(c,q)=1$, and $r\ge1$,
\[
	\sum_{\substack{y<n\le2y\\ n\equiv c\bmod q\\ (n,r)=1}}\beta_n
	=\frac{1}{\varphi(q)}\sum_{\substack{y<n\le2y\\ (n,qr)=1}}\beta_n
	+O_B\big(\tau_k(r)\,y(\log y)^{-B}\big).
\]
When $\beta$ varies in a family of sequences, these estimates are assumed to hold uniformly in the family. The implied constants in Theorems~\ref{thm:conv} and~\ref{thm:nb}, in Corollary~\ref{cor:conv} and in their proofs may depend on $k$, on $\eps$ and on the implied constants in the Siegel-Walfisz condition.

For sequences $\alpha$, $\beta$ supported on $m\sim M$ and $n\sim N$, for $x>0$ and $q\ge1$, set
\[
	E(q):=\sumtwo_{\substack{x<mn\le2x\\ mn\equiv a\bmod q}}\alpha_m\beta_n
	-\frac{1}{\varphi(q)}\sumtwo_{\substack{x<mn\le2x\\ (mn,q)=1}}\alpha_m\beta_n.
\]
Suppose that $|\alpha_m|\le\tau_k(m)$ and $|\beta_n|\le\tau_k(n)$ for all $m,n$, and that $\beta$ satisfies the Siegel-Walfisz condition. Using the estimate of \cite{BC}, Fouvry and Radziwi{\l}{\l} \cite[Corollary~1.1(i)]{FR} proved that for every $B>0$,
\begin{align}
	\sum_{\substack{q\sim Q\\(q,a)=1}}|E(q)|\ll_B x(\log x)^{-B}
	\label{eq:conv_goal}
\end{align}
uniformly for $M,N\ge2$ with $MN/2\le x\le4MN$, provided that
\[
	\exp((\log x)^\eps)\le N\le Q^{-11/12}x^{17/36-\eps}
	\qquad\text{and}\qquad
	1\le|a|\le x/12.
\]
Wright \cite[Corollary~2.2(i)]{Wr1} extended the range of $N$ to $N\le Q^{-33/28}x^{17/28-\eps}$. In both results $Q$ can be as large as $x^{17/33-\eps}$ when $N$ is smaller than every power of $x$.

\begin{theorem}
\label{thm:conv}
Let $k\ge1$ and $\eps>0$, and let $\alpha$, $\beta$ be sequences of complex numbers supported on $m\sim M$ and $n\sim N$ with $|\alpha_m|\le\tau_k(m)$ and $|\beta_n|\le\tau_k(n)$. Suppose that $\beta$ satisfies the Siegel-Walfisz condition. Then for every $B>0$ the bound \eqref{eq:conv_goal} holds uniformly for $M,N\ge2$ with $MN/2\le x\le4MN$, provided that
\begin{align}
	\exp((\log x)^\eps)\le N,
	\qquad
	|a|^6\,Q^{27}N^{35}\le x^{14-\eps},
	\qquad
	1\le|a|\le x/12.
	\label{eq:conv_range}
\end{align}
\end{theorem}

For fixed $a$ and $N\le x^{\eps/2}$, Theorem~\ref{thm:conv} allows $Q\le x^{14/27-\eps}$. Write $Q=x^{\lambda}$ and $N=x^{\varrho}$ with $a$ fixed. Up to the terms in $\eps$, the range \eqref{eq:conv_range} is $27\lambda+35\varrho<14$ and that of \cite[Corollary~2.2(i)]{Wr1} is $33\lambda+28\varrho<17$. The first contains points outside the second if and only if $\varrho<1/133$. For $\varrho\ge1/133$ the range of \cite{Wr1} is the wider one. For $|a|>x^{1/40}$ the range \eqref{eq:conv_range} is contained in that of \cite[Corollary~1.1(i)]{FR}, as the proof shows, so that Theorem~\ref{thm:conv} is new only for smaller residues.

\begin{corollary}
\label{cor:conv}
Let $k\ge1$ and $\eps>0$. Corollaries~1.2, 1.3 and 1.4 of \cite{FR}, under their respective hypotheses (in particular $k\ge2$ in Corollary~1.4), hold with the range $Q\le x^{14/27-\eps}$ in place of $Q\le x^{17/33-\eps}$, for $1\le|a|\le x^{\eps}$. In particular, uniformly for $x\ge2$, $Q\le x^{14/27-\eps}$ and $1\le|a|\le x^\eps$,
\[
	\sum_{\substack{q\sim Q\\(q,a)=1}}\bigg|\sum_{\substack{x<n\le2x\\ n\equiv a\bmod q}}\tau_k(n)
	-\frac{1}{\varphi(q)}\sum_{\substack{x<n\le2x\\(n,q)=1}}\tau_k(n)\bigg|
	\ll_{k,\eps}\frac{x}{(\log x)^{1-\eps}}.
\]
\end{corollary}

We turn to convolutions in which $N$ is slightly larger than $X^{1/2}$, where $X=MN$. For sequences $\alpha$, $\beta$ supported on $m\sim M$ and $n\sim N$, and for $q\ge1$, set
\[
	E_0(q):=\sumtwo_{\substack{m\sim M,\,n\sim N\\ mn\equiv a\bmod q}}\alpha_m\beta_n
	-\frac{1}{\varphi(q)}\sumtwo_{\substack{m\sim M,\,n\sim N\\ (mn,q)=1}}\alpha_m\beta_n.
\]
Suppose that $|\alpha_m|\le\tau_k(m)$ and $|\beta_n|\le\tau_k(n)$ for all $m,n$, and that $\beta$ satisfies the Siegel-Walfisz condition. Using the estimate of \cite{BC}, Fouvry and Radziwi{\l}{\l} \cite[Theorem~1.1]{FR2} proved that for every $\eps>0$ and $B>0$,
\begin{align}
	\sum_{\substack{q\sim Q\\(q,a)=1}}|E_0(q)|\ll_B X(\log X)^{-B}
	\label{eq:nb_goal}
\end{align}
uniformly for $N^{56/23}X^{-17/23+\eps}\le Q\le NX^{-\eps}$ and $1\le|a|\le X$. Wright \cite[Theorem~2.2]{Wr2} extended the range of $Q$ to $N^{34}X^{-17+\eps}\le Q\le NX^{-\eps}$.

\begin{theorem}
\label{thm:nb}
Let $k\ge1$ and $\sigma>0$, and let $\alpha$, $\beta$ be as in Theorem~\ref{thm:conv}. Then for every $B>0$ the bound \eqref{eq:nb_goal} holds uniformly for $M,N\ge2$ and for nonzero integers $a$, provided that
\begin{align}
	X^{1/2}\le Q\le NX^{-\sigma},
	\qquad
	|a|^6\,N^{35}\le X^{14-\sigma}Q^{8}.
	\label{eq:nb_range}
\end{align}
The implied constant may also depend on $\sigma$.
\end{theorem}

Theorem~\ref{thm:nb} concerns the moduli in a dyadic range: for given $M$ and $N$ the second condition in \eqref{eq:nb_range} bounds $Q$ from below. Write $Q=X^{\lambda}$ and $N=X^{\varrho}$ with $a$ fixed, so that $\lambda<\varrho$. Up to the terms in $\eps$ and $\sigma$, the ranges of \cite{FR2}, of \cite{Wr2} and of Theorem~\ref{thm:nb} are
\[
	56\varrho-23\lambda<17,\qquad 34\varrho-\lambda<17,\qquad 35\varrho-8\lambda<14,
\]
respectively. The first two require $\lambda<17/33$, the third allows $1/2\le\lambda<14/27$. For $\lambda=1/2$ they give $\varrho<\tfrac12+\tfrac1{112}$, $\varrho<\tfrac12+\tfrac1{68}$ and $\varrho<\tfrac12+\tfrac1{70}$. Moreover, the estimate for the sums over short intervals displayed in \cite[Section~4]{Wr2} gives, when the factor $(Q/N)^{2/5}$ is retained, the wider range $28\varrho+5\lambda<17$ for $\lambda\ge1/2$, which contains the first two. The range of Theorem~\ref{thm:nb} contains points outside this wider range if and only if $\lambda>29/57$.

\subsection{Ideas of the proofs}
\label{sec:ideas}
We describe the proof of Theorem~\ref{thm:main}. The applications are proved in Section~\ref{sec:moment} by following the arguments of the papers quoted in Section~\ref{sec:new_results}, with Theorem~\ref{thm:main} in place of \cite[Theorem~1]{BC}. The point to be verified in each case is the condition $A\le X^{1/10}$ on the numerators.

The first steps are those of \cite{DFI} and \cite{BC}. By the reciprocity law we may assume that $N\le M$ (Section~\ref{sec:reduction_N_le_M}). The Cauchy-Schwarz inequality in the variable $m$ removes $\alpha$ and leads to the second moment $\Cw(\beta)$ defined in \eqref{eq:def_Cw}, in which $m$ carries a smooth weight. An amplifier built from the primes in $(L,2L]$ and the Dirichlet characters modulo $m$ leads to the amplified second moment $\Dw(\beta)$ defined in \eqref{eq:def_Dw}, and the orthogonality of the characters produces the congruences $p_1n_1\equiv p_2n_2\bmod m$ for primes $p_1,p_2\in(L,2L]$. The off-diagonal terms, in which $p_1\ne p_2$, lead to the correlation sum $\Hcal_g$ defined in \eqref{eq:def_Hg}, where $g$ is the greatest common divisor of the two indices. In $\Hcal_g$ the variable $m$ is replaced by the complementary divisor $d=(p_1n_1-p_2n_2)/m$, the exponential is rewritten by means of the reciprocity law, and a second application of the Cauchy-Schwarz inequality removes $\beta$ (Sections~\ref{sec:complementary_divisor}-\ref{sec:cauchy_schwarz}). After the square is expanded (Section~\ref{sec:expanding_square}), the variables $n_1$ and $n_2$ carry smooth weights, and $n_2$ runs over an arithmetic progression $n_2\equiv\kappa n_1\bmod k$.

At this point \cite{DFI} and \cite{BC} estimate the sum over one of the two variables by Weil's bound. We complete the sum over $n_2$ by Poisson summation (Lemma~\ref{lem:completion_identity}), which gives Kloosterman sums whose modulus is a fixed multiple of $n_1$, and we estimate the sum of these Kloosterman sums over $n_1$ by the spectral theory of automorphic forms. This is the route mentioned in \cite[p.~35]{DFI}. We use the estimates of Deshouillers and Iwaniec in the form given by Pascadi \cite[Corollary~5.5]{P}, which includes the treatment of a fixed index due to Topacogullari \cite{T}. They are stated in Lemma~\ref{lem:spectral} and applied to sums over arithmetic progressions in Lemma~\ref{lem:progression}. The level of the group is divisible by the four amplifier primes, by $g$ and by the complementary divisors, and the exceptional eigenvalues enter through the bound $\theta\le7/64$ of Kim and Sarnak \cite{K}. Their contribution is limited by the fact that the fixed index of the Kloosterman sums is determined by an integer $\Delta$ that is coprime to the amplifier primes (Lemma~\ref{lem:Delta_generic}): these primes enlarge the level but not its greatest common divisor with $\Delta$. The spectral estimate is applied to the terms in which the four amplifier primes are distinct and coprime to the complementary divisors (Section~\ref{sec:generic}). The other terms are fewer, and they are estimated by Weil's bound (Lemma~\ref{lem:completion}) if $\Delta\ne0$ and trivially if $\Delta=0$ (Sections~\ref{sec:delta_zero} and~\ref{sec:nongeneric}).

Two further points differ from \cite{DFI} and \cite{BC}. In these papers the estimates are proved first for sequences $\beta$ supported on squarefree integers, and the general case is recovered by the Cauchy-Schwarz inequality. The reciprocity step only requires that $(n_1,n_2)$ be coprime to $n_1/(n_1,n_2)$ and to $n_2/(n_1,n_2)$ whenever $\beta_{n_1}\beta_{n_2}\ne0$. We say that such sequences have pairwise unitary support, and we write an arbitrary sequence as an average of sequences with this property, by choosing at random, for each prime $p$, the exponent with which $p$ is allowed to occur (Section~\ref{sec:unitary_support}). Secondly, the amplified second moment contains diagonal terms, off-diagonal terms, and terms in which an amplifier prime divides $n_1n_2$. We separate them by an inequality between operators (Lemma~\ref{lem:absorb}), which bounds $\Dw(\beta)$ by its off-diagonal part with indices coprime to the amplifier primes and by a multiple of the second moment without amplifier.

The outcome, for $N\le M$ and with $L$ still free, is the bound \eqref{eq:balancepair},
\[
	\Cw(\beta)\ll_\eps\nn{\beta}^2\nn{\nu}^2\Big(\frac{A^{1/2}M^2}{L}+ALM^{7/4}\Big)M^\eps.
\]
The first term contains the diagonal terms and the other terms controlled by the second moment without amplifier, together with the terms with $\Delta=0$. The second term bounds the remaining contributions. When $N$ is of the same order of magnitude as $M$, its size is that of the terms treated by Weil's bound in Section~\ref{sec:nongeneric}. For the amplifier length chosen below, the contributions of the spectral estimate in Section~\ref{sec:generic}, including that of the exceptional spectrum, are smaller. We take $L=3M^{1/8}A^{-1/4}$, see \eqref{eq:amplifierchoice}, which balances the two terms up to a constant and gives $\Cw(\beta)\ll_\eps\nn{\beta}^2\nn{\nu}^2M^{15/8+\eps}A^{3/4}$, and hence \eqref{eq:main_bound}. The condition $A\le M^{1/10}$ gives $L\ge3A$, so that every amplifier prime is larger than every numerator.

\subsection{Organization of the paper}
\label{sec:organization}
Section~\ref{sec:preliminary_lemmas} contains the preliminary lemmas: the reciprocity law, the completion of incomplete Kloosterman sums and Weil's bound, the spectral estimate for sums of Kloosterman sums over a smooth range of moduli, and its consequence for Kloosterman fractions over an arithmetic progression. In Section~\ref{sec:argument} we state three auxiliary results, Theorems~\ref{thm:Cw}, \ref{thm:Dw} and~\ref{thm:H}, on the second moment, the amplified second moment and the correlation sum, and we deduce Theorem~\ref{thm:main} from them. This section also contains the reduction to sequences with pairwise unitary support and the operator inequality. Theorem~\ref{thm:H} is proved in Section~\ref{sec:correlation}. The applications, Theorems~\ref{thm:zeta}-\ref{thm:nb} and Corollary~\ref{cor:conv}, are proved in Section~\ref{sec:moment}. Section~\ref{sec:conclusion} discusses the critical range of the lengths, the exceptional spectrum, and larger numerators.

\subsection{Notation}
\label{sec:notation}
Throughout, $\e(t)=\exp(2\pi i t)$, and $m\sim M$ stands for $M<m\le 2M$. For integers $u,v$, not both zero, $(u,v)$ denotes their greatest common divisor and $[u,v]$ their least common multiple. We write $\tau(n)$ for the number of positive divisors of $n$, $\varphi(n)$ for Euler's function, $\mu(n)$ for the M\"obius function, and $v_p(n)$ for the exponent of the prime $p$ in $n$, so that $p^v\parallel n$ means $v_p(n)=v$. We shall use the divisor bound $\tau(n)\ll_\eps n^\eps$ repeatedly, see \cite[(1.81)]{IK}.

Multiplicative inverses occur with respect to several moduli at once. When the modulus has to be made explicit we write $\overline{m}^{(n)}$ for the inverse of $m$ modulo $n$, and we take the inverse modulo $1$ to be $0$. As a rule, however, the modulus is the denominator of the fraction in which the inverse appears: in an expression such as
\[
	\frac{\overline{m}}{n_1}+\frac{\overline{m}}{n_2},
\]
the first $\overline{m}$ is the inverse of $m$ modulo $n_1$ and the second is the inverse of $m$ modulo $n_2$. This is the convention of \cite{DFI} and \cite{BC}. A bar over a coefficient, as in $\overline{\beta_n}$, denotes complex conjugation.

Sequences of complex numbers are indexed by positive integers and are supported on dyadic intervals, that is, they vanish outside of the indicated interval.

We write \(\sum_{\pm}\) to mean the sum of the \(+\) and \(-\) cases. In each case, all occurrences of \(\pm\) take the chosen sign. Thus the condition \(\pm d>0\) means \(d>0\) in the \(+\) case and \(d<0\) in the \(-\) case. A sum $\sum_{1\le R<D}$ over a letter $R$ runs over the powers of two in the indicated range. Two summation signs, as in $\sumtwo_{a_1,a_1'\in\Acal}$, denote a sum over two variables subject to a common set of conditions. We write $\log_2$ for the logarithm to the base $2$. 

The letter $\eps$ denotes an arbitrarily small positive number, not necessarily the same at each occurrence, and all implied constants are allowed to depend on $\eps$. Further notation, in particular the smooth weight $w$, the set $\Pcal$ of amplifier primes and the sums $\Cw$, $\Dw$ and $\Hcal_g$, is introduced at the beginning of Section~\ref{sec:argument}.

\section{Preliminary lemmas}
\label{sec:preliminary_lemmas}

\subsection{Reciprocity and the Chinese remainder theorem}
\label{sec:reciprocity}
The following identity is the elementary reciprocity law for Kloosterman fractions, see \cite[(2.15)]{DFI}; it is used throughout \cite{DFI} and \cite{BC}.
\begin{lemma}
\label{lem:reciprocity}
Let $x$ and $y$ be coprime positive integers. Then
\begin{align}
	\frac{\overline{x}}{y}+\frac{\overline{y}}{x}
	\equiv
	\frac{1}{xy} \bmod 1.
	\label{eq:reciprocity}
\end{align}
\end{lemma}
\begin{proof}
We have $x\overline{x}\equiv 1 \bmod y$ and $y\overline{y}\equiv 1 \bmod x$. Hence $x\overline{x}+y\overline{y}\equiv 1$ both modulo $y$ and modulo $x$, and therefore modulo $xy$, since $(x,y)=1$. Dividing by $xy$ gives \eqref{eq:reciprocity}.
\end{proof}
We also need to split an inverse modulo a product of two coprime moduli into inverses modulo the factors.
\begin{lemma}
\label{lem:crt_inverse}
Let $y$ and $z$ be coprime positive integers, and let $x$ be an integer coprime to $yz$. Then
\begin{align}
	\frac{\overline{x}}{yz}
	\equiv
	\frac{\overline{xz}}{y}+\frac{\overline{xy}}{z} \bmod 1,
	\label{eq:crt_inverse}
\end{align}
where, by our convention, the three inverses are taken modulo $yz$, $y$ and $z$, respectively.
\end{lemma}
\begin{proof}
Set $u=\overline{xz}^{(y)}z+\overline{xy}^{(z)}y$. Then $xu\equiv x\overline{xz}^{(y)}z\equiv1 \bmod y$ and $xu\equiv x\overline{xy}^{(z)}y\equiv 1\bmod z$. Since $(y,z)=1$, we have $xu\equiv1\bmod yz$. Hence $u\equiv\overline{x}\bmod yz$, and dividing by $yz$ gives \eqref{eq:crt_inverse}.
\end{proof}

\subsection{Kloosterman sums and Weil's bound}
\label{sec:weil}
For $c\in\N$ and $u,v\in\Z$, the Kloosterman sum is
\begin{align}
	S(u,v;c):=\sum_{\substack{z\bmod c\\(z,c)=1}}\e\Big(\frac{u\overline{z}+vz}{c}\Big).
	\label{eq:def_kloosterman}
\end{align}
By Weil's bound,
\begin{align}
	|S(u,v;c)|\le\tau(c)\,(u,v,c)^{1/2}\,c^{1/2},
	\label{eq:weil}
\end{align}
which is valid for all integers $u,v$ and all $c\in\N$, see \cite{We} and \cite[Corollary 11.12]{IK}. For $v=0$ the sum \eqref{eq:def_kloosterman} is a Ramanujan sum, and a better bound is available. Detecting the condition $(z,c)=1$ by M\"obius inversion, we get
\[
	S(u,0;c)
	=\sum_{z\bmod c}\e\Big(\frac{uz}{c}\Big)\sum_{r\mid(z,c)}\mu(r)
	=\sum_{r\mid c}\mu(r)\sum_{z'\bmod c/r}\e\Big(\frac{uz'}{c/r}\Big)
	=\sum_{\substack{r\mid c\\(c/r)\mid u}}\mu(r)\,\frac{c}{r},
\]
and therefore
\begin{align}
	|S(u,0;c)|\le\sum_{s\mid(u,c)}s\le\tau(c)\,(u,c).
	\label{eq:ramanujan}
\end{align}

To estimate smoothly weighted sums of Kloosterman fractions over an arithmetic progression, we first express them in terms of Kloosterman sums, for which bounds are available. The following lemma gives this exact transformation. Its proof uses a Fourier expansion modulo $ck$, followed by Poisson summation.

\begin{lemma}
\label{lem:completion_identity}
Fix $c,k\in\N$ with $(c,k)=1$, and $r,u\in\Z$. Let $\Omega:\R\to\C$ be a smooth function with compact support, with Fourier transform $\widehat{\Omega}(\xi)=\int_\R\Omega(x)\e(-\xi x)\,dx$. Then
\begin{align}
	\sum_{\substack{x\equiv r\bmod k\\(x,c)=1}}\Omega(x)\,\e\Big(\frac{u\overline{x}}{c}\Big)
	=
	\frac{1}{ck}\sum_{h\in\Z}\e\Big(\frac{hr\overline{c}}{k}\Big)\,S(u\overline{k},h;c)\,\widehat{\Omega}\Big(\frac{h}{ck}\Big),
	\label{eq:completion_identity}
\end{align}
where $\overline{x}$ and $\overline{k}$ are inverses modulo $c$ and $\overline{c}$ is an inverse modulo $k$.
\end{lemma}
\begin{proof}
Let 
\[
f(x)=
\begin{cases}
\e\!\Big(\dfrac{u\overline{x}}{c}\Big),
& \text{if } x\equiv r\bmod{k}\ \text{and } (x,c)=1,\\[6pt]
0, & \text{otherwise}.
\end{cases}
\]
The function $f$ depends only on $x$ modulo $k$ and on $x$ modulo $c$, hence it has period $ck$, and
\[
	f(x)=\frac{1}{ck}\sum_{h\bmod ck}\widehat{f}(h)\,\e\Big(\frac{hx}{ck}\Big),
	\qquad
	\widehat{f}(h):=\sum_{z\bmod ck}f(z)\,\e\Big(-\frac{hz}{ck}\Big).
\]
We compute $\widehat{f}(h)$. Since $(c,k)=1$, the residue classes $z\bmod ck$ with $f(z)\ne0$ are those of the form $z\equiv rc\,\overline{c}+wk\,\overline{k}\bmod ck$, where $w$ runs over the reduced residues modulo $c$. For such $z$ we have $z/(ck)\equiv r\overline{c}/k+w\overline{k}/c\bmod 1$ and $z\equiv w\bmod c$, so that $\overline{z}\equiv\overline{w}\bmod c$. Hence
\[
	\widehat{f}(h)
	=\e\Big(-\frac{hr\overline{c}}{k}\Big)\sum_{\substack{w\bmod c\\(w,c)=1}}\e\Big(\frac{u\overline{w}-h\overline{k}\,w}{c}\Big)
	=\e\Big(-\frac{hr\overline{c}}{k}\Big)\,S(u\overline{k},-h;c),
\]
where the last step follows from replacing the summation variable $w$ by $kw$. By Poisson summation,
\begin{align*}
	\sum_{x\in\Z}\Omega(x)f(x)
	=\frac{1}{ck}\sum_{h\bmod ck}\widehat{f}(h)\sum_{x\in\Z}\Omega(x)\e\Big(\frac{hx}{ck}\Big)
	&=\frac{1}{ck}\sum_{h\bmod ck}\widehat{f}(h)\sum_{m\in\Z}\widehat{\Omega}\Big(m-\frac{h}{ck}\Big)
	\\
	&=\frac{1}{ck}\sum_{h\in\Z}\widehat{f}(h)\,\widehat{\Omega}\Big(-\frac{h}{ck}\Big),
\end{align*}
where in the last step we combined $h\bmod ck$ and $m\in\Z$ into the single integer $h-mck$, using the periodicity of $\widehat{f}$. Inserting the value of $\widehat{f}(h)$ and replacing $h$ by $-h$ gives \eqref{eq:completion_identity}.
\end{proof}
We estimate the term \(h=0\) in \eqref{eq:completion_identity} using the Ramanujan sum bound in \eqref{eq:ramanujan}, and the remaining terms \(h\ne 0\) using Weil's bound in \eqref{eq:weil}. This gives the following estimate, which will be used in Section~\ref{sec:nongeneric}. It holds for every \(\ell>0\), uniformly in the location of the interval supporting the weight.

\begin{lemma}
\label{lem:completion}
Fix $c,k\in\N$ with $(c,k)=1$, $r,u\in\Z$,
and $\ell>0$. Suppose that $F:\R\to\C$ is smooth,
supported in an interval of length at most $\ell$, and satisfies
\[
|F^{(j)}(x)|\le C_j\ell^{-j}
\qquad (j=0,1,2)
\]
for all $x\in\R$.
Then
\begin{align}
	\Big|\sum_{\substack{x\equiv r\bmod k\\(x,c)=1}}F(x)\,\e\Big(\frac{u\overline{x}}{c}\Big)\Big|
	\le
	C_0\,\frac{\ell}{ck}\,\tau(c)\,(u,c)+8\max\Big(C_0,\frac{C_2}{4\pi^2}\Big)\tau(c)^2\sqrt{c}.
	\label{eq:completion_bound}
\end{align}
\end{lemma}
\begin{proof}
We apply Lemma~\ref{lem:completion_identity} with $\Omega=F$ 
to the left side of \eqref{eq:completion_bound}, and obtain
\[
\sum_{\substack{x\equiv r\bmod k\\(x,c)=1}}F(x)\,\e\Big(\frac{u\overline{x}}{c}\Big)=\frac{1}{ck}\sum_{h\in\Z}
\e\Big(\frac{hr\bar c}{k}\Big)
S(u\bar k,h;c)\,
\widehat F\Big(\frac{h}{ck}\Big).
\]
Taking absolute values therefore gives
\begin{align}
\Big|\sum_{\substack{x\equiv r\bmod k\\(x,c)=1}}F(x)\,\e\Big(\frac{u\overline{x}}{c}\Big)\Big|\le \frac{1}{ck}\sum_{h\in\Z}
|S(u\bar k,h;c)|\,
\Big|\widehat F\Big(\frac{h}{ck}\Big)\Big|.\label{eq:bound for ks in Lem 2.4}
\end{align}
We now need bounds for the two factors in each summand in \eqref{eq:bound for ks in Lem 2.4}.
Since $\overline{k}$ is a unit modulo $c$, $(u\overline{k},h,c)=(u,h,c)$, so that by \eqref{eq:ramanujan} and \eqref{eq:weil}
\begin{align}
	|S(u\overline{k},0;c)|\le\tau(c)\,(u,c),
	\qquad
	|S(u\overline{k},h;c)|\le\tau(c)\,(u,h,c)^{1/2}\sqrt{c}\qquad(h\ne0).
	\label{eq:fhat_bounds}
\end{align}
We now bound $\widehat{F}$. By definition,
\[
\widehat F(\xi)=\int_{\R}F(x)\e(-\xi x)\,dx.
\]
Since $F$ is supported in an interval of length at most $\ell$ and $|F|\le C_0$, we have 
\[
|\widehat{F}(\xi)|\le \int_{\R}|F(x)|\,dx
\le C_0\ell.
\]
For $\xi\ne0$, integrating by parts twice gives
\[
\widehat F(\xi)
=\frac{1}{(2\pi i\xi)^2}
\int_{\R}F''(x)\e(-\xi x)\,dx.
\]
Since $|F''(x)|\le C_2\ell^{-2}$,
\[
\int_{\R}|F''(x)|\,dx\le C_2\ell^{-1},
\]
and hence
\[
|\widehat F(\xi)|
\le \frac{C_2}{4\pi^2\ell\xi^2}.
\]

 Writing $C'=\max(C_0,C_2/(4\pi^2))$ and combining the two bounds, we get 
\[
 |\widehat{F}(\xi)|\le C'\ell\min(1,(\ell|\xi|)^{-2})\le4C'\ell\,(1+\ell|\xi|)^{-2}.
\] 

Combining with the Ramanujan sum bound in \eqref{eq:fhat_bounds}, we see that the term $h=0$ of \eqref{eq:completion_identity} contributes at most 
\[
\frac{1}{ck}|S(u\bar k,0;c)|\,|\widehat F(0)|
\le \frac{C_0\ell}{ck}\tau(c)(u,c),
\]
 which is the first term in \eqref{eq:completion_bound}. For the terms $h\ne0$, Weil's bound in \eqref{eq:fhat_bounds} gives
\[
	\frac{1}{ck}\sum_{h\ne0}
|S(u\bar k,h;c)|\,
\Big|\widehat F\Big(\frac{h}{ck}\Big)\Big|
\le
\frac{4C'\ell\tau(c)\sqrt c}{ck}
\sum_{h\ne0}(u,h,c)^{1/2}
\Big(1+\frac{\ell|h|}{ck}\Big)^{-2}.
\]
Since $(u,h,c)$ is a positive divisor of $(h,c)$, we have
\[
(u,h,c)^{1/2}\le \sum_{j\mid(h,c)}j^{1/2}.
\]
Therefore,
\begin{align}
\sum_{h\ne0}(u,h,c)^{1/2}
\Big(1+\frac{\ell|h|}{ck}\Big)^{-2}
\le
\sum_{j\mid c}j^{1/2}
\sum_{\substack{h\ne0\\j\mid h}}
\Big(1+\frac{\ell|h|}{ck}\Big)^{-2}=
2\sum_{j\mid c}j^{1/2}
\sum_{m\ge1}
\Big(1+\frac{\ell jm}{ck}\Big)^{-2},\label{eq: bound for h neq 0 in Lem 2.4}
\end{align}
where the last step follows from writing $h = \pm jm$ with $m\ge 1$.

Since
\[
\sum_{m\ge1}(1+am)^{-2}\le\int_0^\infty(1+at)^{-2}\,dt=1/a
\]
for $a>0$, the inner sum in \eqref{eq: bound for h neq 0 in Lem 2.4} is at most $ck/(\ell j)$, and
\[
	\sum_{h\ne0}|S(u\overline{k},h;c)|\,\Big|\widehat{F}\Big(\frac{h}{ck}\Big)\Big|
	\le
	8C'ck\,\tau(c)\sqrt{c}\sum_{j\mid c}j^{-1/2}
	\le
	8C'ck\,\tau(c)^2\sqrt{c}.
\]
Dividing by $ck$ gives the second term in \eqref{eq:completion_bound}.
\end{proof}

\subsection{Sums of Kloosterman sums over a smooth range of moduli}
\label{sec:spectral_input}
The second ingredient is an estimate that exploits cancellation among Kloosterman sums as the modulus varies. Rather than bounding each Kloosterman sum separately as in Lemma~\ref{lem:completion}, this estimate treats the sums collectively and allows a sum over the second index, with arbitrary coefficients. Such bounds come from the Kuznetsov formula for the congruence group $\Gamma_0(q)$ and the spectral large sieve of Deshouillers and Iwaniec \cite{DI}. Their shape depends on the exceptional spectrum: for $q\in\N$ let $\lambda_1(q)$ be the smallest positive eigenvalue of the hyperbolic Laplacian on $L^2(\Gamma_0(q)\backslash\mathbb H)$, and let
\begin{align}
	\theta(q):=\sqrt{\max\big(0,\tfrac14-\lambda_1(q)\big)}.
	\label{eq:def_theta}
\end{align}
Selberg's conjecture asserts that $\theta(q)=0$ for every $q$; the best known bound is that of Kim and Sarnak \cite[Appendix 2]{K},
\begin{align}
	\theta(q)\le\frac{7}{64}\qquad(q\in\N).
	\label{eq:kim_sarnak}
\end{align}
We use Pascadi's formulation \cite{P} of the Deshouillers-Iwaniec estimates, which allows arbitrary coefficients and includes Topacogullari's refinement \cite{T} for the fixed index. The following lemma is a special case of \cite[Corollary 5.5]{P}, written in terms of classical Kloosterman sums.
\begin{lemma}
\label{lem:spectral}
Fix $r,s\in\N$ with $(r,s)=1$, and set $q=rs$.
 Suppose that $m_0\in\N$, $H\ge1/2$, and $\mathscr C$ is larger than a suitable absolute constant. Let $(\gamma_h)$ be a sequence of complex numbers supported on $h\sim H$, and $\Phi:(0,\infty)^2\to\C$ be a smooth function supported in a fixed compact subset of $(0,\infty)^2$ with $\partial_x^j\partial_y^k\Phi\ll_{j,k}1$ for all $j,k\ge0$. Then, for either choice of the sign $\pm$ and every $\eps>0$,
\begin{align}
	&\Big|\sum_{h\sim H}\gamma_h\sum_{\substack{v\ge1\\(v,r)=1}}\Phi\Big(\frac{h}{H},\frac{sv\sqrt{r}}{\mathscr C}\Big)S(m_0\overline{r},\pm h;sv)\Big|
	\notag\\
	&\quad\ll_\eps
	(qm_0H\mathscr C)^{\eps}\,(1+T_0)^{2\theta(q)}\,\frac{\mathscr C^2\,\nn{\gamma}}{\mathscr C+\sqrt{m_0H}}
	\Big(1+\frac{m_0H}{\mathscr C^2}+\frac{\sqrt{(q,m_0)\,m_0}}{q}\Big)^{1/2}\Big(1+\frac{m_0H}{\mathscr C^2}+\frac{H}{q}\Big)^{1/2},
	\label{eq:spectral_bound}
\end{align}
where $\nn{\gamma}=(\sum_h|\gamma_h|^2)^{1/2}$ and
\begin{align}
	T_0:=\frac{\mathscr C}{\max\big(m_0,\,q^2/(q,m_0)\big)^{1/2}\max(H,q)^{1/2}}
	\le\frac{\mathscr C\,(q,m_0)^{1/2}}{q^{3/2}}.
	\label{eq:def_T0}
\end{align}
The implied constant depends on $\eps$, on the compact set containing the support of $\Phi$ and on the constants in the bounds for the derivatives of $\Phi$.
\end{lemma}
\begin{proof}
We apply \cite[Corollary~5.5]{P} with the cusps
\[
\mathfrak a=\infty,\qquad \mathfrak b=1/s
\]
of \(\Gamma_0(q)\), using the scaling matrices specified in \cite[(3.9)]{P}. Since \(q=rs\) and \((r,s)=1\), \cite[(3.6)]{P} gives
\[
\mu(\mathfrak a)=\mu(\mathfrak b)=q^{-1},
\]
which is required by  \cite[Corollary~5.5]{P}. 

We first express our sum in Pascadi's notation. 
For the cusps $\mathfrak a=\infty$ and $\mathfrak b=1/s$, Lemma~3.2 of \cite{P}, which is \cite[(3.5)]{Wa}, shows that
\begin{align}
	S_{\infty,1/s}(m_0,\pm h;sv\sqrt{r})=\ind_{v\in\N,\ (v,r)=1}\,S(m_0\overline{r},\pm h;sv).
	\label{eq:lemma_B}
\end{align}
Thus, our sum has the form required by \cite[Corollary~5.5]{P} with  
\[
m=m_0, \qquad n=h,\qquad N=H,\qquad Z=1,\qquad
a_h=\gamma_h,\qquad A_N=\|\gamma\|_2,\qquad Y_N=1.
\]
It remains to verify \cite[Assumption 5.4]{P} on these coefficients and the weight hypotheses in \cite[Corollary 5.5]{P}. For every $\xi\in\R$, 
\[
\sum_{h\sim H}
\Big|\e\Big(\frac{h\xi}{H}\Big)\gamma_h\Big|^2
=
\|\gamma\|_2^2.
\]

Applying \cite[Theorem 1.2]{P} to the sequence
$(\e(h\xi/H)\gamma_h)$, we obtain
\begin{align}
	\sum_{\lambda_j<1/4}X^{2\theta_j}\Big|\sum_{h\sim H}\e\Big(\frac{h\xi}{H}\Big)\gamma_h\,\rho_{j\mathfrak a}(h)\Big|^2
	\ll_\eps(qH)^{\eps}\Big(1+\frac{H}{q}\Big)\nn{\gamma}^2,
	\qquad\theta_j=\sqrt{\tfrac14-\lambda_j},
	\label{eq:assumption_14}
\end{align}
when $X\ll\max(1,q/H)$. This holds for every cusp $\mathfrak a$ with
$\mu(\mathfrak a)=q^{-1}$, with the spectral notation and
normalization of \cite[(3.15)]{P}. In particular, it holds in the range
\[
X\ll\frac{\max(1,q/H)}{1+|\xi|^2}
\]
required by \cite[Assumption 5.4]{P} with $Y_N=1$. Assumption 5.4 is therefore satisfied. 
The weight hypotheses in \cite[Corollary 5.5]{P} are also satisfied:
the fixed compact support of $\Phi$ gives $x,y\asymp1$, and with
$Z=1$ the required derivative bounds are precisely
\[
\partial_x^j\partial_y^k\Phi\ll_{j,k}1
\]
for $j,k\ge 0$. We may now apply \cite[Corollary 5.5]{P}. Since $Y_N=1$, we have
\[
T=T_0
=
\frac{\mathscr C}{
\max\big(m_0,q^2/(q,m_0)\big)^{1/2}
\max(H,q)^{1/2}}.
\]
Consequently, with a sufficiently small choice of its
$\eps$, \cite[(5.20)]{P} gives \eqref{eq:spectral_bound}. Finally,
\[
T_0
\le
\frac{\mathscr C}{
\big(q^2/(q,m_0)\big)^{1/2}q^{1/2}}
=
\frac{\mathscr C\,(q,m_0)^{1/2}}{q^{3/2}},
\]
which proves the inequality in \eqref{eq:def_T0}.
\end{proof}
We record the special case in which Lemma~\ref{lem:spectral} will be applied. Suppose that, for some $\Lambda\ge1$,
\begin{align}
	m_0\le\Lambda q,\qquad H\le\Lambda q,\qquad m_0H\le\Lambda \mathscr C^2.
	\label{eq:restricted_hypotheses}
\end{align}
Then, since $(q,m_0)\le m_0$, 
\[
\frac{\sqrt{(q,m_0)m_0}}{q}\le \frac{m_0}{q}\le\Lambda.
\]
Thus, each of the two square root factors in \eqref{eq:spectral_bound} is at most $(1+2\Lambda)^{1/2}$. Their product is therefore at most $(1+2\Lambda)\le 3\Lambda$. Moreover,
\[
\frac{\mathscr C^2}{\mathscr C+\sqrt{m_0H}}\le \mathscr C.
\]
Together with \eqref{eq:def_T0} and \eqref{eq:kim_sarnak}, this gives
\begin{align}
	\Big|\sum_{h\sim H}\gamma_h\sum_{\substack{v\ge1\\(v,r)=1}}\Phi\Big(\frac{h}{H},\frac{sv\sqrt{r}}{\mathscr C}\Big)S(m_0\overline{r},\pm h;sv)\Big|
	\ll_\eps
	(qm_0H\mathscr C)^{\eps}\,\Lambda\,\mathscr C\,\nn{\gamma}\,\Big(1+\frac{\mathscr C\,(q,m_0)^{1/2}}{q^{3/2}}\Big)^{7/32}
	\label{eq:spectral_restricted}
\end{align}
under the assumptions in \eqref{eq:restricted_hypotheses}. In our application, the amplifier primes divide \(q\) but are coprime to \(m_0\). They therefore do not contribute to \((q,m_0)^{1/2}\) in \eqref{eq:spectral_restricted}, which helps control the factor arising from the exceptional spectrum.

\subsection{A bound for Kloosterman fractions over an arithmetic progression}
\label{sec:progression_lemma}
We now combine the completion identity in Lemma~\ref{lem:completion_identity} with the spectral estimate in Lemma~\ref{lem:spectral} to bound a sum over two variables $n_1,n_2$, with the restriction that $n_2\equiv\kappa n_1\bmod k$, where $\kappa$ is a unit modulo $k$. Completing the sum over $n_2$ produces Kloosterman sums $S(\cdot,h;Pn_1)$.
\begin{lemma}
\label{lem:progression}
Fix $Z\ge1$ and $C\ge1$. Suppose that
\[
P,Q,k\in\N,\qquad
\kappa,\Delta\in\Z,\qquad
\Delta\ne0,\qquad 0<\delta\le1,
\]
with
\[
(P,Qk)=(\kappa,k)=(\Delta,P)=1,
\]
and that
\begin{align}
	P\le Z,\qquad
	|\Delta|\le C\,PQk,\qquad
	\delta Q\ge1,\qquad
	k\le C\,\delta Z^2,\qquad
	\max(P,Q,k,|\Delta|)\le Z^{C}.
	\label{eq:progression_hypotheses}
\end{align}
Consider a smooth function  $W:\R^2\to\C$ supported in $[1/2,3]^2$, and set $F(t,y):=W(t,ty)$. Suppose that $F(t,y)=0$ unless $y$ belongs to a fixed interval $I$
of length at most $C\delta$, and that
\begin{align}
	|\partial_t^i\partial_y^jF(t,y)|\le C_{ij}\,\delta^{-j}\qquad(i,j\ge0,\ t,y\in\R)
	\label{eq:F_derivatives}
\end{align}
for some constants $C_{ij}$. Then for every $\eps>0$, the sum
\begin{align}
	\Jcal:=\sum_{\substack{n_1\ge1\\(n_1,PQk)=1}}\ \sum_{\substack{n_2\ge1,\ n_2\equiv\kappa n_1\bmod k\\(n_2,PQn_1)=1}}W\Big(\frac{n_1}{Z},\frac{n_2}{Z}\Big)\e\Big(\frac{\Delta\,\overline{Qn_2}}{Pn_1}\Big)
	\label{eq:def_J_lemma}
\end{align}
satisfies 
\begin{align}
	|\Jcal|\ll_\eps Z^{\eps}\Big\{\frac{\delta Z}{Pk}+Z\sqrt{PQ\delta}\,\Big[1+\Big(\frac{Z\sqrt{(Qk,\Delta)}}{\sqrt{P}\,Qk}\Big)^{7/32}\Big]\Big\}.
	\label{eq:progression_bound}
\end{align}
The implied constant depends on $\eps$, $C$ and the constants $C_{ij}$.
\end{lemma}
The summation conditions give $(Qn_2,Pn_1)=1$, so the inverse in
\eqref{eq:def_J_lemma} is well-defined. They also imply $(n_2,k)=1$,
since $n_2\equiv\kappa n_1\bmod k$ and $(\kappa n_1,k)=1$. The hypotheses allow $Q$ and $k$ to have common factors, and
$\Delta$ need not be coprime to either of them. We require only
$(P,Qk)=(P,\Delta)=1$. 

The number of pairs $n_1,n_2$ in the sum is $\ll Z+\delta Z^2/k$. In \eqref{eq:progression_bound}, the first term comes from the term $h=0$ after completion, and the second term comes from the terms with $h\ne0$. The factor in square brackets accounts for the contribution of the exceptional spectrum, which is at most $2$ when 
\[
Z\sqrt{(Qk,\Delta)}\le\sqrt{P}\,Qk.
\]
 The derivative condition \eqref{eq:F_derivatives} allows $W$ to vary at scale $\delta$ in $y=n_2/n_1$, while requiring uniform
smoothness in $t=n_1/Z$. This smoothness is needed to apply
Lemma~\ref{lem:spectral} to the sum over moduli.

We begin the proof with an exact completion identity. For $f\mid Q$ with $(f,k)=1$, let
\begin{align}
	H_0:=fPk
	\quadand
	\alpha_0:=\frac{\kappa\,\overline{fP}}{k},
	\label{eq:def_H0_alpha0}
\end{align}
where both $H_0$ and $\alpha_0$ depend on $f$. 
\begin{lemma}
\label{lem:J_completion}
Fix $P,Q,k,\kappa,\Delta,Z$ as in Lemma~\ref{lem:progression}. Let $W:\R^2\to\C$ be a smooth function with compact support in $(0,\infty)^2$, and set $F(t,y):=W(t,ty)$. Then
\begin{align}
	\Jcal
	=
	\sum_{e\mid P}\ \sum_{\substack{f\mid Q\\(f,k)=1}}\mu(e)\mu(f)
	\sum_{\substack{v\ge1\\(v,Qk)=1}}
	\frac{1}{H_0}\sum_{h\in\Z}S\big(\Delta\,\overline{Qfk},h;Pev\big)
	\int_\R F\Big(\frac{ev}{Z},y\Big)\e\Big(\Big(\alpha_0-\frac{y}{H_0}\Big)h\Big)\,dy,
	\label{eq:J_completed}
\end{align}
where $\Jcal$ is defined in \eqref{eq:def_J_lemma} and $H_0,\alpha_0$ are
defined in \eqref{eq:def_H0_alpha0}.
\end{lemma}
\begin{proof}
The original coprimality conditions in \eqref{eq:def_J_lemma} split as
\begin{align*}
(n_1,PQk)=1
&\iff (n_1,P)=1\ \text{and}\ (n_1,Qk)=1,\\
(n_2,PQn_1)=1
&\iff (n_2,Q)=1\ \text{and}\ (n_2,Pn_1)=1.
\end{align*}
We detect the conditions $(n_1,P)=1$ and $(n_2,Q)=1$ by M\"obius inversion:
\[
	\ind_{(n_1,P)=1}=\sum_{e\mid(n_1,P)}\mu(e),
	\qquad
	\ind_{(n_2,Q)=1}=\sum_{f\mid(n_2,Q)}\mu(f),
\]
and write $n_1=ev$, $n_2=fj$. We now handle the remaining conditions. Since $e\mid P$ and
$(P,Qk)=1$, we have $(e,Qk)=1$, so
\[
(n_1,Qk)=1
\iff (ev,Qk)=1
\iff (v,Qk)=1.
\]
Next, the original progression condition states that
\[
n_2\equiv\kappa n_1\bmod k.
\]
As discussed earlier, since both $\kappa$ and $n_1$ are coprime to $k$, we must have $(n_2,k)=1$, and thus $(f,k)=1$. Therefore, the congruence $n_2\equiv\kappa n_1\bmod k$ is equivalent to 
\[
j\equiv\kappa\overline{f}n_1\bmod k.
\]
Finally, we have $f\mid Q$, $(P,Q)=1$, and $(n_1,Q)=1$. Hence $(f,Pn_1)=1$. Therefore,
\[
(n_2,Pn_1)=1
\iff(fj,Pn_1)=1
\iff(j,Pn_1)=1.
\]

Putting these together, we obtain
\begin{align}
	\Jcal
	=
	\sum_{e\mid P}\ \sum_{\substack{f\mid Q\\(f,k)=1}}\mu(e)\mu(f)
	\sum_{\substack{v\ge1\\(v,Qk)=1}}
	\sum_{\substack{j\equiv\kappa\overline{f}n_1\bmod k\\(j,Pn_1)=1}}
	W\Big(\frac{n_1}{Z},\frac{fj}{Z}\Big)\e\Big(\frac{\Delta\,\overline{Qfj}}{Pn_1}\Big),
	\qquad n_1=ev.
	\label{eq:J_after_mobius}
\end{align}
Fix $e$, $f$, and $v$, and set $t=n_1/Z$. We apply Lemma~\ref{lem:completion_identity} to the inner sum in
\eqref{eq:J_after_mobius}, with
\[
c=Pn_1,\qquad
r=\kappa\bar f n_1,\qquad
u=\Delta\,\overline{Qf},
\qquad
\Omega(x)=W\Big(t,\frac{fx}{Z}\Big).
\]
These choices give
\[
u\bar j
\equiv
\Delta\,\overline{Qf}\,\bar j
\equiv
\Delta\,\overline{Qfj}
\bmod{Pn_1}.
\]
The hypothesis $(c,k)=1$ holds because $(P,k)=1$ and
$(n_1,k)=1$. Thus, Lemma~\ref{lem:completion_identity} gives
\begin{align}
\sum_{\substack{j\equiv\kappa\overline{f}n_1\bmod k\\(j,Pn_1)=1}}
	W\Big(\frac{n_1}{Z},\frac{fj}{Z}\Big)\e\Big(\frac{\Delta\,\overline{Qfj}}{Pn_1}\Big)=\frac{1}{kPn_1}
\sum_{h\in\Z}
\e\Big(\frac{hr\bar c}{k}\Big)
S(u\bar k,h;Pn_1)\,
\widehat\Omega\Big(\frac{h}{kPn_1}\Big).\label{eq:inner sum in Lem 2.7}
\end{align}

We now simplify each factor in this expression. First, 
we have 
\[
u\overline{k}\equiv\Delta\,\overline{Qfk}\bmod Pn_1,
\]
so the Kloosterman sum is $S(\Delta\,\overline{Qfk},h;Pn_1)$. For the exponential factor, since $n_1\overline{Pn_1}\equiv\overline{P}\bmod k$,
\[
	\frac{r\,\overline{c}}{k}=\frac{\kappa\overline{f}\,n_1\overline{Pn_1}}{k}\equiv\frac{\kappa\,\overline{fP}}{k}=\alpha_0\bmod 1.
\]
Observe that the factor $\e(\alpha_0h)$ does not depend on $n_1=ev$. Next, by the definition of the Fourier transform,
\[
\widehat\Omega\Big(\frac{h}{kPn_1}\Big)
=
\int_\R
W\Big(t,\frac{fx}{Z}\Big)
\e\Big(-\frac{hx}{kPn_1}\Big)\,dx.
\]
The substitution $x=Zty/f$ gives
\begin{align*}
	\widehat{\Omega}\Big(\frac{h}{kPn_1}\Big)
	=\int_\R W\Big(t,\frac{fx}{Z}\Big)\e\Big(-\frac{hx}{kPn_1}\Big)\,dx
	&=\frac{Zt}{f}\int_\R W(t,ty)\,\e\Big(-\frac{hZty}{fkPn_1}\Big)\,dy
	\\
	&=\frac{n_1}{f}\int_\R F(t,y)\,\e\Big(-\frac{hy}{H_0}\Big)\,dy,
\end{align*}
where the last step follows from $Zt=n_1$ and $H_0$ is defined in \eqref{eq:def_H0_alpha0}. Substituting back into \eqref{eq:inner sum in Lem 2.7}, 
\begin{align*}
	\sum_{\substack{j\equiv\kappa\overline{f}n_1\bmod k\\(j,Pn_1)=1}}
	W\Big(\frac{n_1}{Z},\frac{fj}{Z}\Big)&\e\Big(\frac{\Delta\,\overline{Qfj}}{Pn_1}\Big)
	\\
	&\quad=
	\frac{1}{H_0}\sum_{h\in\Z}\e(\alpha_0h)\,S\big(\Delta\,\overline{Qfk},h;Pn_1\big)\int_\R F(t,y)\,\e\Big(-\frac{hy}{H_0}\Big)\,dy.
\end{align*}
Substituting this expression into \eqref{eq:J_after_mobius} with $n_1=ev$ and $t=ev/Z$ proves \eqref{eq:J_completed}.
\end{proof}

\begin{proof}[Proof of Lemma~\ref{lem:progression}]
We start from the completed expression in Lemma~\ref{lem:J_completion}. We write
\begin{align}
	\Jcal=\Jcal_{h=0}+\Jcal_{h\ne0}
	\label{eq:J_split}
\end{align}
for the contributions to \eqref{eq:J_completed} from $h=0$
and $h\ne0$, respectively. The divisor
bound and the hypotheses $P,Q\le Z^C$ in \eqref{eq:progression_hypotheses} give
\[
\tau(P)\tau(Q)\ll_\eps Z^\eps.
\]
We therefore fix $e\mid P$ and $f\mid Q$ in \eqref{eq:J_completed}, with $(f,k)=1$, and
sum over these divisors at the end. Throughout, $\eps>0$ is fixed, and the implied constants may depend on $\eps$, $C$ and the constants $C_{ij}$. 

Since
\[
F\Big(\frac{ev}{Z},y\Big)
=
W\Big(\frac{ev}{Z},\frac{ev}{Z}y\Big)
\]
and $W$ is supported in $[1/2,3]^2$, the support of $F$ restricts $v$ to $v\le3Z/e$.

We begin with $h=0$. The support and size assumptions on $F$ give
\begin{align}
\Big|\int_{\R}F(t,y)\,dy\Big|
\le CC_{00}\delta.\label{eq:F_integral_trivial}
\end{align}
The Kloosterman sum in \eqref{eq:J_completed} is now a Ramanujan sum. Moreover, since $e\mid P$ and $(\Delta,P)=1$, we have $(\Delta,Pe)=1$. By \eqref{eq:ramanujan},
\[
	|S(\Delta\,\overline{Qfk},0;Pev)|\le\tau(Pev)\,(\Delta\,\overline{Qfk},Pev)=\tau(Pev)\,(\Delta,Pev)=\tau(Pev)\,(\Delta,v).
\]
Using the identity $(\Delta,v)=\sum_{s\mid(\Delta,v)}\varphi(s)$, we get
\[
	\sum_{v\le3Z/e}(\Delta,v)
	=\sum_{s\mid|\Delta|}\varphi(s)\,\#\{v\le3Z/e : s\mid v\}
	\le\sum_{s\mid|\Delta|}\varphi(s)\,\frac{3Z}{es}
	\le\frac{3Z}{e}\,\tau(|\Delta|).
\]
Since $\tau(Pev)\tau(|\Delta|)\ll Z^{\eps}$ by \eqref{eq:progression_hypotheses}, this and \eqref{eq:F_integral_trivial} give, for each $e$ and $f$,
\[
	\frac{1}{H_0}\sum_{\substack{v\le3Z/e\\(v,Qk)=1}}|S(\Delta\,\overline{Qfk},0;Pev)|\,\Big|\int_\R F\Big(\frac{ev}{Z},y\Big)dy\Big|
	\ll Z^{\eps}\,\frac{\delta}{fPk}\cdot\frac{Z}{e}
	\le Z^{\eps}\,\frac{\delta Z}{Pk},
\]
and summing over $e$ and $f$,
\begin{align}
	|\Jcal_{h=0}|\ll Z^{2\eps}\,\frac{\delta Z}{Pk}.
	\label{eq:J_zero_bound}
\end{align}
This is the first term of \eqref{eq:progression_bound}, with $2\eps$ in place of $\eps$.

We now consider $h\ne 0$. Recall that $H_0=fPk$, and set
\[
H_1:=\frac{H_0}{\delta}.
\]

Using \eqref{eq:F_derivatives} and the length of the support, for any integer $J\ge2$, integration by parts $J$ times gives
\begin{align*}
\Big|
\int_{\R}F(t,y)\e\Big(-\frac{hy}{H_0}\Big)\,dy
\Big|\le
\Big(\frac{H_0}{2\pi|h|}\Big)^J
\int_{\R}|\partial_y^JF(t,y)|\,dy&\le
CC_{0J}\delta
\Big(\frac{H_0}{2\pi|h|\delta}\Big)^J\\
&\ll_J
\delta\Big(\frac{H_1}{|h|}\Big)^J.
\end{align*}
There are no boundary terms, since $F$ is smooth and compactly supported. 

We use this decay to estimate the part with
$|h|>H_1Z^\eps$. The trivial Kloosterman sum bound gives
\[
|S(\Delta\,\overline{Qfk},h;Pev)|\le Pev\le3PZ,
\]
and we also have
\[
\sum_{|h|>H_1Z^{\eps}}\Big(\frac{H_1}{|h|}\Big)^{J}
	\le4H_1^{J}\,(H_1Z^{\eps})^{1-J}
	\ll H_1\,Z^{-\eps(J-1)},
\]
which follows from $H_1Z^{\eps}\ge1$ and $\sum_{h>X}h^{-J}\le2X^{1-J}$ for $X\ge1$. Since there are $O(Z)$ possible values of \(v\), the contribution for fixed \(e,f\) is
\begin{align*}
	\frac{1}{H_0}\sum_{v\le3Z/e}\sum_{|h|>H_1Z^{\eps}}|S(\Delta\,\overline{Qfk},h;Pev)|\,&\Big|\int_\R F\Big(\frac{ev}{Z},y\Big)\e\Big(-\frac{hy}{H_0}\Big)dy\Big|
	\\
	&\quad\ll\frac{1}{H_0}\cdot Z\cdot PZ\cdot\delta H_1Z^{-\eps(J-1)}
	=PZ^{2-\eps(J-1)},
\end{align*}
where the last equality follows from $\delta H_1=H_0$. Since $P\le Z$, we choose $J\ge1+3/\eps$, which gives $PZ^{2-\eps(J-1)}\le1$. Therefore, the part $\Jcal_{\rm tail}$ of $\Jcal_{h\ne0}$ with $|h|>H_1Z^\eps$, summed over $e$ and $f$, satisfies
\begin{align}
	|\Jcal_{\rm tail}|\ll Z^{\eps}.
	\label{eq:J_tail_bound}
\end{align}
This is absorbed by the second term in \eqref{eq:progression_bound}, since $Z\sqrt{PQ\delta}\ge Z\sqrt{\delta Q}\ge Z\ge1$ by \eqref{eq:progression_hypotheses}.

It remains to bound the part of $\Jcal_{h\ne0}$ with $0<|h|\le H_1Z^{\eps}$. We separate the two signs and divide the absolute values into dyadic intervals, that is, 
\[
h\sim H,  \qquad h\le H_1Z^{\eps}
\]
and 
\[
-h\sim H, \qquad -h\le H_1Z^{\eps}
\]
with $H=\frac{1}{2}, 1, 2,4,\cdots,$ and $H\le H_1Z^{\eps}$. 

Fix $e,f$, one sign, and $H$, and set 
\[
\gamma_h=
\mathbf1_{h\le H_1Z^\eps}
\e\bigg(\pm\Big(\alpha_0-\frac y{H_0}\Big)h\bigg)\qquad(h\sim H).
\]
Interchanging the sum and the integral, the contribution of this block is 
\[
	\Jcal(e,f,\pm,H):=\frac{1}{H_0}\sum_{\substack{v\ge1\\(v,Qk)=1}}\sum_{\substack{h\sim H\\h\le H_1Z^{\eps}}}S\big(\Delta\,\overline{Qfk},\pm h;Pev\big)\int_\R F\Big(\frac{ev}{Z},y\Big)\e\Big(\pm\Big(\alpha_0-\frac{y}{H_0}\Big)h\Big)\,dy,
\]
so that $\Jcal_{h\ne0}-\Jcal_{\rm tail}$ is the sum of $\mu(e)\mu(f)\Jcal(e,f,\pm,H)$ over $e$, $f$, the sign and $H$. For each fixed $y$, apply Lemma~\ref{lem:spectral} with
\begin{align}
	s=Pe,\qquad r=Qfk,\qquad q=sr=PeQfk,\qquad m_0=|\Delta|,\qquad \mathscr C=PZ\sqrt{Qfk}.
	\label{eq:spectral_parameters}
\end{align}
We first verify its hypotheses and match its sum with the expression above. The required coprimality follows from
\[
(r,s)=(Qfk,Pe)=1,
\]
since $(P,Qk)=1$, $e\mid P$, and $f\mid Q$. The restriction $(v,r)=1$ in Lemma~\ref{lem:spectral} becomes
$(v,Qfk)=1$. Since $f\mid Q$,
\[
(v,Qfk)=1\iff(v,Qk)=1,
\]
which is precisely the restriction on $v$ in
\eqref{eq:J_completed}.
 Moreover,
\[
sv=Pev,\qquad \bar r=\overline{Qfk}\bmod{Pev},
\]
so the modulus and inverse are those required in $S(m_0\bar r,\pm h;sv)$. Since $\Delta\ne0$, we have $m_0=|\Delta|\in\N$. If $\Delta<0$, we use
\[
S(-u,h;c)=S(u,-h;c),
\]
which follows by replacing $z$ by $-z$ in the defining sum. This makes the first parameter positive, with a possible change in the sign of the second index.

We next check the coefficient and weight hypotheses. The coefficients $\gamma_h$ are supported on $h\sim H$, are independent of $v$, and satisfy $|\gamma_h|\le1$. Since $H\ge1/2$,
\[
\|\gamma\|_2\le\sqrt{H+1}\le2\sqrt H.
\]
For the weight, choose
\[
\Phi(x,x')=\phi(x)F(x',y),
\]
where $\phi$ is a fixed smooth function supported in $[1/2,3]$ and equal to $1$ on $[1,2]$. Since
\[
\frac{sv\sqrt r}{\mathscr C}
=
\frac{Pev\sqrt{Qfk}}{PZ\sqrt{Qfk}}
=
\frac{ev}{Z},
\]
we have, for $h\sim H$,
\[
\Phi\Big(\frac hH,\frac{sv\sqrt r}{\mathscr C}\Big)
=
F\Big(\frac{ev}{Z},y\Big),
\]
which is the weight in
\eqref{eq:J_completed}.
Its support lies in $[1/2,3]^2$, and its derivatives involve only derivatives of $\phi$ and $t$-derivatives of $F$. The assumed bounds with $j=0$ therefore give
\[
\partial_x^i\partial_{x'}^{i'}\Phi\ll_{i,i'}1
\]
uniformly in $y$. Finally, $\mathscr C\ge Z$, so the required lower bound on $\mathscr C$ holds for sufficiently large $Z$. For bounded $Z$, the conclusion follows directly from the support and size bounds on $W$. Therefore, all hypotheses of Lemma~\ref{lem:spectral} are verified.

For each fixed $y$, we apply \eqref{eq:spectral_restricted},
the consequence of Lemma~\ref{lem:spectral} obtained under
\eqref{eq:restricted_hypotheses}, with
\[
\Lambda=C^2Z^\eps.
\]
The extra inequality hypotheses in \eqref{eq:restricted_hypotheses} hold because 
\[
m_0=|\Delta|\le CPQk\le Cq,
\qquad
\frac Hq
\le
\frac{fPkZ^\eps}{\delta PeQfk}
=
\frac{Z^\eps}{\delta eQ}
\le Z^\eps,
\]
and
\[
\frac{m_0H}{\mathscr C^2}
\le
\frac{CPQk\cdot fPkZ^\eps/\delta}
{P^2Z^2Qfk}
=
\frac{CkZ^\eps}{\delta Z^2}
\le C^2Z^\eps,
\]
using $\delta Q\ge1$ and $k\le C\delta Z^2$. Applying \eqref{eq:spectral_restricted} gives, for each fixed $y$,
\begin{align}
&\Big|
\sum_{h\sim H}\gamma_h
\sum_{\substack{v\ge1\\(v,Qk)=1}}
F\Big(\frac{ev}{Z},y\Big)
S(\Delta\,\overline{Qfk},\pm h;Pev)
\Big|\ll_\eps
(qm_0H\mathscr C)^\eps
Z^\eps\mathscr C\sqrt H
\Big(
1+\frac{\mathscr C\sqrt{(q,m_0)}}{q^{3/2}}
\Big)^{7/32}.\label{eq:summation bound in lem 2.7}
\end{align}
The bounds in \eqref{eq:progression_hypotheses}, together with
$\delta^{-1}\le Q$ and $H\le fPkZ^\eps/\delta$, imply
\[
(qm_0H\mathscr C)^\eps\ll Z^{O(\eps)}.
\]
Since $(\Delta,Pe)=1$, we also have
\[
(q,m_0)=(Qfk,\Delta)\le f(Qk,\Delta).
\]
Consequently,
\[
\frac{\mathscr C\sqrt{(q,m_0)}}{q^{3/2}}
\le
\frac{PZ\sqrt{Qfk}\sqrt{f(Qk,\Delta)}}{(PeQfk)^{3/2}}=
\frac{Z\sqrt{(Qk,\Delta)}}{\sqrt P\,e^{3/2}\sqrt f\,Qk}\le
\frac{Z\sqrt{(Qk,\Delta)}}{\sqrt P\,Qk}.
\]
Substituting all back into \eqref{eq:summation bound in lem 2.7} and integrating the result over $I$, which has length at most $C\delta$,
 we obtain
\[
|\Jcal(e,f,\pm,H)|
\ll
Z^{O(\eps)}
\frac{\delta\mathscr C\sqrt H}{H_0}
\Big(
1+\frac{Z\sqrt{(Qk,\Delta)}}{\sqrt P\,Qk}
\Big)^{7/32}.
\]
The dyadic sum satisfies
\[
\sum_{\substack{H\text{ dyadic}\\H\le H_1Z^\eps}}\sqrt H
\ll
\sqrt{H_1}Z^{\eps/2},
\]
and
\[
\frac{\delta\mathscr C\sqrt{H_1}}{H_0}
=
\frac{\delta PZ\sqrt{Qfk}\sqrt{fPk/\delta}}{fPk}
=
Z\sqrt{PQ\delta}.
\]
Summing over the dyadic intervals, the two signs, and the
divisors $e,f$ therefore gives
\begin{align}
|\Jcal_{h\ne0}-\Jcal_{\rm tail}|
&\ll
Z^{O(\eps)}Z\sqrt{PQ\delta}
\Big(
1+\frac{Z\sqrt{(Qk,\Delta)}}{\sqrt P\,Qk}
\Big)^{7/32}.
\label{eq:J_nonzero_bound}
\end{align}

Combining \eqref{eq:J_zero_bound}, \eqref{eq:J_tail_bound},
and \eqref{eq:J_nonzero_bound}, and using
\[
(1+x)^{7/32}\le1+x^{7/32}
\]
for $x\ge 0$, we obtain
\[
|\Jcal|
\ll
Z^{O(\eps)}
\Big\{
\frac{\delta Z}{Pk}
+
Z\sqrt{PQ\delta}
\Big[
1+
\Big(
\frac{Z\sqrt{(Qk,\Delta)}}{\sqrt P\,Qk}
\Big)^{7/32}
\Big]
\Big\}.
\]
Applying the argument with a sufficiently small multiple of $\eps$ in place of
$\eps$ gives \eqref{eq:progression_bound}.
\end{proof}

\section{The argument for the trilinear estimate}
\label{sec:argument}
In this section we state three auxiliary results, Theorems~\ref{thm:Cw}, \ref{thm:Dw} and \ref{thm:H}, and show how they imply Theorem~\ref{thm:main}. Theorem~\ref{thm:H} is proved in Section~\ref{sec:correlation}.

We first introduce the objects involved. Throughout this section and the next, $\alpha$, $\beta$ and $\nu$ are as in Theorem~\ref{thm:main}, and we set
\begin{align}
	\Acal:=\{a\in\N : A<a\le 2A\}.
	\label{eq:def_Acal}
\end{align}
Fix once and for all a real function $w\in C_c^\infty((1/2,3))$ with $0\le w\le 1$ and $w(t)=1$ for $t\in[1,2]$. The smooth weight $w(m/M)$ replaces the sharp cutoff $m\sim M$; it is harmless in an upper bound, and it is needed later, when the summation over $m$ is converted into a summation over $n_1$ and $n_2$ that is completed by Poisson summation. All implied constants may depend on $w$.

As in \cite{DFI} and \cite{BC}, we amplify with the primes in a dyadic interval. Let $L\ge 2$ be a parameter, to be chosen in Section~\ref{sec:Dw_implies_Cw}, and set
\begin{align}
	\Pcal:=\{p \text{ prime} : L<p\le 2L\}.
	\label{eq:def_Pcal}
\end{align}
Unless stated otherwise, the letters $p_1$, $p_2$, $p_1'$ and $p_2'$ denote elements of $\Pcal$.

For a sequence $\beta$ supported on $n\sim N$ and for $(m,n)=1$, set
\begin{align}
	V_n(m)
	:=
	\sum_{a\in\Acal}\nu_a\,\e\Big(\frac{a\overline{m}}{n}\Big),
	\label{eq:def_Vn}
\end{align}
and define the second moment
\begin{align}
	\Cw(\beta)
	:=
	\sum_{m\in\N} w\Big(\frac{m}{M}\Big)\bigg|\sum_{\substack{n\sim N\\(m,n)=1}}\beta_nV_n(m)\bigg|^2
	\label{eq:def_Cw}
\end{align}
and the amplified second moment
\begin{align}
	\Dw(\beta)
	:=
	\sum_{m\in\N}\frac{w(m/M)}{\varphi(m)}\sum_{\chi \bmod m}
	\bigg|\sum_{p\in\Pcal}\chi(p)\bigg|^2
	\bigg|\sum_{\substack{n\sim N\\(m,n)=1}}\beta_n\chi(n)V_n(m)\bigg|^2,
	\label{eq:def_Dw}
\end{align}
where $\chi$ runs over the Dirichlet characters modulo $m$. Both quantities depend on $A$, $M$, $N$ and $\nu$, and $\Dw$ also on $L$; we suppress this dependence from the notation.

Finally, for $g\in\N$, real $N>0$ and a sequence $\beta$ supported on $n\sim N$, set
\begin{align}
	\Hcal_g(\beta,N)
	:=&
	\sum_{\substack{p_1,p_2\in\Pcal\\p_1\ne p_2\\(g,p_1p_2)=1}}
	\sum_{\substack{m\in\N\\(m,gp_1p_2)=1}} w\Big(\frac{m}{M}\Big)
	\sum_{\substack{n_1,n_2\sim N\\(n_1,n_2)=(mgp_1p_2,n_1n_2)=1\\p_1n_1\equiv p_2n_2 \bmod m}}
	\beta_{n_1}\overline{\beta_{n_2}}
	\notag\\
	&\quad\times
	\sumtwo_{a_1,a_2\in\Acal}\nu_{a_1}\overline{\nu_{a_2}}\,
	\e\Big(\frac{a_1\overline{m}}{gn_1}-\frac{a_2\overline{m}}{gn_2}\Big).
	\label{eq:def_Hg}
\end{align}
Here, by the convention of Section~\ref{sec:notation}, the first $\overline{m}$ is the inverse of $m$ modulo $gn_1$ and the second is the inverse of $m$ modulo $gn_2$; an empty sum is zero. The sum $\Hcal_g$ will arise from the off-diagonal terms of $\Dw$ after the common divisor $g=(n_1,n_2)$ of the two indices has been taken out; it will therefore be applied with $N/g$ in place of $N$.

An important step in \cite{DFI} and \cite{BC} is to establish the estimates first for sequences supported on squarefree integers and to recover the general case afterwards by the Cauchy-Schwarz inequality. We use a different reduction, which rests on the following property. A divisor $g$ of a positive integer $n$ is called unitary if $(g,n/g)=1$. We say that a sequence $\beta=(\beta_n)$ of complex numbers has pairwise unitary support if, whenever $\beta_n\ne 0$ and $\beta_u\ne 0$, the greatest common divisor $(n,u)$ is a unitary divisor of both $n$ and $u$, that is,
\begin{align}
	\Big((n,u),\frac{n}{(n,u)}\Big)=1 \quadand \Big((n,u),\frac{u}{(n,u)}\Big)=1.
	\label{eq:def_pairwise_unitary}
\end{align}
In other words, every prime dividing both $n$ and $u$ occurs with the same exponent in $n$ and in $u$. Every sequence supported on squarefree integers has pairwise unitary support. The point of the definition is that $(n,u)$ is then invertible modulo $n/(n,u)$ and modulo $u/(n,u)$, which is what the reciprocity step in Section~\ref{sec:correlation} requires.

We can now state the three auxiliary results. The first one is the bound for the second moment $\Cw$ from which Theorem~\ref{thm:main} follows.
\begin{theorem}
\label{thm:Cw}
Let $M\ge N\ge 2$ and $1/2\le A\le M^{1/10}$. Let $\beta$ be a sequence supported on $n\sim N$ with pairwise unitary support and let $\nu$ be a sequence supported on $a\sim A$. Then for every $\eps>0$,
\begin{align}
	\Cw(\beta)
	\ll_\eps
	\nn{\beta}^2\nn{\nu}^2\,M^{15/8+\eps}A^{3/4}.
	\label{eq:bound_Cw}
\end{align}
\end{theorem}
The second one is the bound for the amplified second moment $\Dw$, in which the amplifier length $L$ is still free. It holds under the hypotheses
\begin{align}
	M\ge N\ge 2,\qquad
	L\ge 2\log(2M),\qquad
	\frac{1}{2}\le A<\frac{\min(N,L)}{2},\qquad
	M\ge 2^{26}L^4.
	\label{eq:blockrange}
\end{align}
The last condition keeps the amplifier short compared with the modulus. It is used in Section~\ref{sec:correlation}, where the long variables have to be longer than the products of two amplifier primes, and in Section~\ref{sec:Dw_implies_Cw} in the form $L\le M^{1/4}$.
\begin{theorem}
\label{thm:Dw}
Assume that $M$, $N$, $A$ and $L$ satisfy \eqref{eq:blockrange}. Let $\beta$ and $\nu$ be sequences supported on $n\sim N$
and $a\sim A$, respectively, with $\beta$ having pairwise
unitary support. Then for every $\eps>0$,
\begin{align}
	\Dw(\beta)
	\ll_\eps
	\frac{L^2}{M}\nn{\beta}^2\nn{\nu}^2
	\Big(\frac{M^2}{L}+\frac{A^{1/2}MN}{L}+ALN^{7/4}+AL^2N^{3/2}+AL^{25/16}M^{7/64}N^{3/2}\Big)M^\eps.
	\label{eq:bound_Dw}
\end{align}
\end{theorem}
The third one is the bound for the correlation sum $\Hcal_g$, from which Theorem~\ref{thm:Dw} is deduced in Section~\ref{sec:H_implies_Dw}.
\begin{theorem}
\label{thm:H}
Assume that $M,N,A,L$ satisfy \eqref{eq:blockrange}, and fix
$g\in\N$. Let $\beta$ and $\nu$ be sequences supported on
$n\sim N/g$ and $a\sim A$, respectively. Suppose that $\beta$
has pairwise unitary support and that $\beta_n=0$ whenever
$(n,g)>1$.  Then for every $\eps>0$,
\begin{align}
	|\Hcal_g(\beta,N/g)|
	\ll_\eps
	\nn{\beta}^2\nn{\nu}^2
	\Big(\frac{A^{1/2}LN}{g^{1/2}}
	+\frac{AL^3N^{7/4}}{Mg^{3/4}}
	+\frac{AL^4N^{3/2}}{Mg^{1/4}}
	+\frac{AL^{57/16}N^{3/2}}{M^{57/64}g^{23/64}}\Big)M^\eps.
	\label{eq:bound_H}
\end{align}
\end{theorem}

\subsection{Reduction to \texorpdfstring{$N\le M$}{N <= M}}
\label{sec:reduction_N_le_M}
We show that it suffices to prove Theorem~\ref{thm:main} when $N\le M$. Assume that the theorem holds in this case, and let $N>M$. Since only coprime pairs $m,n$ occur in $\Bcal(A,M,N)$, Lemma~\ref{lem:reciprocity} gives
\[
	\e\Big(\frac{a\overline{m}}{n}\Big)
	=
	\e\Big(\frac{a}{mn}\Big)\e\Big(-\frac{a\overline{n}}{m}\Big),
\]
where $\overline{n}$ is the inverse of $n$ modulo $m$. Inserting this into \eqref{eq:def_B_AMN} and taking complex conjugates, we get
\[
	|\Bcal(A,M,N)|
	=
	\bigg|\sum_{a\sim A}\sumtwo_{\substack{m\sim M,\, n\sim N\\(m,n)=1}}
	\overline{\nu_a}\,\overline{\alpha_m}\,\overline{\beta_n}\,
	\e\Big(\frac{a\overline{n}}{m}\Big)\e\Big(-\frac{a}{mn}\Big)\bigg|.
\]
We expand the last factor into its power series,
\[
	\e\Big(-\frac{a}{mn}\Big)
	=
	\sum_{j=0}^{\infty}\frac{(-2\pi i)^j}{j!}\,\frac{a^j}{m^jn^j},
\]
and obtain
\begin{align}
	|\Bcal(A,M,N)|
	\le
	\sum_{j=0}^{\infty}\frac{(2\pi)^j}{j!}
	\bigg|\sum_{a\sim A}\sumtwo_{\substack{m\sim M,\, n\sim N\\(m,n)=1}}
	\big(a^j\overline{\nu_a}\big)\big(n^{-j}\overline{\beta_n}\big)\big(m^{-j}\overline{\alpha_m}\big)
	\e\Big(\frac{a\overline{n}}{m}\Big)\bigg|.
	\label{eq:reduction_series}
\end{align}
The $j$-th sum on the right is the trilinear form \eqref{eq:def_B_AMN} with the roles of $m$ and $n$ interchanged: the variable in the modulus is now $m\sim M$ and the inverted variable is $n\sim N$. Since $N>M$, the modulus variable is the shorter one, which is the case already established. Theorem~\ref{thm:main} therefore applies to it, with the sequences $a^j\overline{\nu_a}$, $n^{-j}\overline{\beta_n}$ and $m^{-j}\overline{\alpha_m}$ in place of $\nu$, $\alpha$ and $\beta$, and with the same $X=\max(M,N)$. Since
\[
	\nn{a^j\overline{\nu_a}}\le(2A)^j\nn{\nu},
	\qquad
	\nn{n^{-j}\overline{\beta_n}}\le N^{-j}\nn{\beta},
	\qquad
	\nn{m^{-j}\overline{\alpha_m}}\le M^{-j}\nn{\alpha},
\]
and since the implied constant in \eqref{eq:main_bound} depends only on $\eps$, we get
\[
	|\Bcal(A,M,N)|
	\ll_\eps
	\nn{\alpha}\nn{\beta}\nn{\nu}\,X^{15/16+\eps}A^{3/8}
	\sum_{j=0}^{\infty}\frac{1}{j!}\Big(\frac{4\pi A}{MN}\Big)^j
	\ll
	\nn{\alpha}\nn{\beta}\nn{\nu}\,X^{15/16+\eps}A^{3/8},
\]
since $MN\ge 2X$ because $M,N\ge 2$, and $A\le X^{1/10}$, so that $4\pi A/(MN)\le 2\pi$. This proves Theorem~\ref{thm:main} for $N>M$.

\subsection{Two elementary estimates}
\label{sec:elementary_estimates}
We relate $\Bcal(A,M,N)$, $\Cw(\beta)$ and $\Dw(\beta)$ to one another. Neither of the two estimates below requires the sequence $\beta$ to have pairwise unitary support. The first one is
\begin{align}
	|\Bcal(A,M,N)|^2
	\le
	\nn{\alpha}^2\,\Cw(\beta).
	\label{eq:firstCS}
\end{align}
Indeed, writing the sum over $a$ inside the sum over $n$ by means of \eqref{eq:def_Vn} and applying the Cauchy-Schwarz inequality to the sum over $m$, we get
\[
	|\Bcal(A,M,N)|^2
	=
	\bigg|\sum_{m\sim M}\alpha_m\sum_{\substack{n\sim N\\(m,n)=1}}\beta_nV_n(m)\bigg|^2
	\le
	\nn{\alpha}^2\sum_{m\sim M}\bigg|\sum_{\substack{n\sim N\\(m,n)=1}}\beta_nV_n(m)\bigg|^2
	\le
	\nn{\alpha}^2\,\Cw(\beta),
\]
since $w\ge0$ and $w(m/M)=1$ for $m\sim M$.

The second estimate is the amplification inequality. By the prime number theorem there is an absolute constant $L_0$ such that
\begin{align}
	\pi(2L)-\pi(L)\ge\frac{3L}{4\log L}
	\qquad (L\ge L_0),
	\label{eq:pnt}
\end{align}
where $\pi$ denotes the prime counting function. We claim that if $L\ge\max(L_0,8)$ and $L\ge 2\log(2M)$, then
\begin{align}
	\Cw(\beta)
	\le
	\frac{192\,M(\log L)^2}{L^2}\,\Dw(\beta).
	\label{eq:amplification}
\end{align}
Let $\chi_0$ denote the principal character modulo $m$, and let $m<3M$, which is the case on the support of $w(m/M)$. The number of primes $p>L$ dividing $m$ is at most $\log m/\log L<\log(3M)/\log L$, and $\log(3M)\le\log(2M)+1\le L/2+1$. Hence, by \eqref{eq:pnt},
\[
	\sum_{p\in\Pcal}\chi_0(p)
	=
	\#\{p\in\Pcal : p\nmid m\}
	\ge
	\pi(2L)-\pi(L)-\frac{\log(3M)}{\log L}
	\ge
	\frac{3L}{4\log L}-\frac{L/2+1}{\log L}
	\ge
	\frac{L}{8\log L},
\]
where the last inequality uses $L\ge 8$. Keeping only the principal character in \eqref{eq:def_Dw} and using $\varphi(m)\le m<3M$, we get
\begin{align*}
	\Dw(\beta)
	&\ge
	\sum_{m\in\N}\frac{w(m/M)}{\varphi(m)}\bigg|\sum_{p\in\Pcal}\chi_0(p)\bigg|^2
	\bigg|\sum_{\substack{n\sim N\\(m,n)=1}}\beta_nV_n(m)\bigg|^2
	\\
	&\ge
	\frac{L^2}{64(\log L)^2}\cdot\frac{1}{3M}\sum_{m\in\N}w\Big(\frac{m}{M}\Big)\bigg|\sum_{\substack{n\sim N\\(m,n)=1}}\beta_nV_n(m)\bigg|^2,
\end{align*}
which is \eqref{eq:amplification}. The factor $M/L^2$ in \eqref{eq:amplification} is the price of the amplifier. It is compensated by the saving in the off-diagonal terms of $\Dw(\beta)$, where the congruence $p_1n_1\equiv p_2n_2 \bmod m$ produced by the characters can be exploited.

\subsection{Sequences with pairwise unitary support}
\label{sec:unitary_support}
We show that a bound for $\Cw$ on sequences with pairwise unitary support implies a bound for $\Cw$ on all sequences. The idea is to write an arbitrary sequence as an average of sequences with pairwise unitary support, which is done by choosing at random, for each prime $p$, the only exponent with which $p$ is allowed to occur.
\begin{lemma}
\label{lem:random_weights}
Let $N\ge1$. There exist a probability space $(\Omega,\mathbb{P})$ and random variables $P_n\ge0$, $n\sim N$, on $\Omega$ such that
\begin{align}
	\E{P_n}=1
	\quadand
	\E{P_n^2}\le\tau(n)^2
	\qquad(n\sim N),
	\label{eq:Pn_moments}
\end{align}
and such that, for every $\omega\in\Omega$, the set
\begin{align}
	I_\omega:=\{n\sim N : P_n(\omega)>0\}
	\label{eq:def_Iomega}
\end{align}
has the following property: if $n,u\in I_\omega$, then $(n,u)$ is a unitary divisor of $n$ and of $u$.
\end{lemma}
\begin{proof}
For each prime $p\le 2N$ let $V_p$ be a random variable with values in $\N$ and
\begin{align}
	\Pb{V_p=v}=\frac{1}{v(v+1)}
	\qquad(v\in\N),
	\label{eq:def_Vp}
\end{align}
and let the $V_p$ be independent. This is a probability distribution, since $\sum_{v\ge1}1/(v(v+1))=1$ by telescoping. Since there are finitely many primes $p\le 2N$, the random variables $V_p$ can be realized on the countable set $\Omega=\N^{\pi(2N)}$ with the product measure $\mathbb{P}$.

For $\omega\in\Omega$ we call an integer $n\sim N$ permitted if $v_p(n)=V_p(\omega)$ for every prime $p\mid n$; equivalently, $v_p(n)\in\{0,V_p(\omega)\}$ for every prime $p\le 2N$. By independence,
\begin{align}
	\Pb{n \text{ is permitted}}
	=
	\prod_{p\mid n}\Pb{V_p=v_p(n)}
	=
	\prod_{p^v\parallel n}\frac{1}{v(v+1)}.
	\label{eq:prob_permitted}
\end{align}
Set
\[
	P_n(\omega)
	:=
	\frac{\ind_{n \text{ is permitted}}}{\Pb{n \text{ is permitted}}}.
\]
Then $\E{P_n}=1$, and
\[
	\E{P_n^2}
	=
	\frac{1}{\Pb{n \text{ is permitted}}}
	=
	\prod_{p^v\parallel n}v(v+1)
	\le
	\prod_{p^v\parallel n}(v+1)^2
	=
	\tau(n)^2,
\]
which proves \eqref{eq:Pn_moments}. Finally, $I_\omega$ is the set of permitted integers. If $n$ and $u$ are both permitted and $p$ is a prime dividing both, then $v_p(n)=V_p(\omega)=v_p(u)$. Hence every prime factor of $(n,u)$ occurs with the same exponent in $n$ and in $u$, so that $(n,u)$ is a unitary divisor of both.
\end{proof}
For example, $12=2^2\cdot3$ and $20=2^2\cdot5$ are permitted together when $V_2=2$, $V_3=1$ and $V_5=1$, and their greatest common divisor $4$ is coprime to $3$ and to $5$; whereas $12$ and $18=2\cdot 3^2$ are never permitted together, since they disagree at $2$ and at $3$.

For a sequence $\beta$ supported on $n\sim N$ we write $P\beta$ for the random sequence $(P_n\beta_n)_{n\sim N}$. By Lemma~\ref{lem:random_weights}, $P\beta$ has pairwise unitary support for every $\omega\in\Omega$.
\begin{lemma}
\label{lem:random_average}
Let $N\ge1$ and let $P_n$, $n\sim N$, be as in Lemma~\ref{lem:random_weights}. Then for every sequence $\beta$ supported on $n\sim N$,
\begin{align}
	\Cw(\beta)
	\le
	\E{\Cw(P\beta)}.
	\label{eq:Cw_le_expectation}
\end{align}
Consequently, if $C\ge0$ is such that $\Cw(\widetilde{\beta})\le C\nn{\widetilde{\beta}}^2$ for every sequence $\widetilde{\beta}$ supported on $n\sim N$ with pairwise unitary support, then for every sequence $\beta$ supported on $n\sim N$ and every $\eps>0$,
\begin{align}
	\Cw(\beta)
	\ll_\eps
	C\nn{\beta}^2N^{\eps}.
	\label{eq:Cw_general}
\end{align}
\end{lemma}
\begin{proof}
Since $\E{P_n}=1$, we have $\beta_n=\E{P_n\beta_n}$ for every $n$. The sums over $n$ in \eqref{eq:def_Cw} are finite, so the expectation can be taken outside, and $|\E{Z}|^2\le\E{|Z|^2}$ for every bounded random variable $Z$. Hence
\begin{align*}
	\Cw(\beta)
	=
	\sum_{m\in\N}w\Big(\frac{m}{M}\Big)\bigg|\,\mathbb{E}\Big[\sum_{\substack{n\sim N\\(m,n)=1}}P_n\beta_nV_n(m)\Big]\bigg|^2
	&\le
	\sum_{m\in\N}w\Big(\frac{m}{M}\Big)\,\mathbb{E}\Big[\,\bigg|\sum_{\substack{n\sim N\\(m,n)=1}}P_n\beta_nV_n(m)\bigg|^2\Big] \\
	&=
	\E{\Cw(P\beta)},
\end{align*}
where the sum over $m$ is finite as well. This is \eqref{eq:Cw_le_expectation}.

Now suppose that $\Cw(\widetilde{\beta})\le C\nn{\widetilde{\beta}}^2$ for every $\widetilde{\beta}$ supported on $n\sim N$ with pairwise unitary support. Since $P\beta$ is such a sequence for every $\omega\in\Omega$, \eqref{eq:Cw_le_expectation} and \eqref{eq:Pn_moments} give
\[
	\Cw(\beta)
	\le
	\E{C\nn{P\beta}^2}
	=
	C\sum_{n\sim N}\E{P_n^2}|\beta_n|^2
	\le
	C\sum_{n\sim N}\tau(n)^2|\beta_n|^2
	\ll_\eps
	C\nn{\beta}^2N^{\eps},
\]
by the divisor bound.
\end{proof}

\subsection{Theorem~\ref{thm:Cw} implies Theorem~\ref{thm:main}}
\label{sec:Cw_implies_main}
By Section~\ref{sec:reduction_N_le_M} we may assume that $N\le M$, so that $X=M$, and the hypotheses of Theorem~\ref{thm:main} are then those of Theorem~\ref{thm:Cw}. Let $\eps>0$. Theorem~\ref{thm:Cw} says that
\[
	\Cw(\widetilde{\beta})
	\le
	C\nn{\widetilde{\beta}}^2
	\quad\text{with}\quad
	C=c(\eps)\nn{\nu}^2M^{15/8+\eps}A^{3/4}
\]
for every sequence $\widetilde{\beta}$ supported on $n\sim N$ with pairwise unitary support, where $c(\eps)$ depends only on $\eps$. By Lemma~\ref{lem:random_average} and $N\le M$, we therefore have
\[
	\Cw(\beta)
	\ll_\eps
	\nn{\beta}^2\nn{\nu}^2M^{15/8+2\eps}A^{3/4}
\]
for every sequence $\beta$ supported on $n\sim N$. Inserting this into \eqref{eq:firstCS} and renaming $\eps$ gives \eqref{eq:main_bound}.

Corollary~\ref{cor:bilinear} follows by taking $A=a/2$ and $\nu$ the sequence with $\nu_a=1$ and $\nu_{a'}=0$ for $a'\ne a$. Then $a\sim A$, $1/2\le A\le X^{1/10}$, $\nn{\nu}=1$, $\Bcal(A,M,N)=\Bcal_a(M,N)$, and $A^{3/8}\le a^{3/8}$.

\subsection{Theorem~\ref{thm:Dw} implies Theorem~\ref{thm:Cw}}
\label{sec:Dw_implies_Cw}
Let $M\ge N\ge2$, $1/2\le A\le M^{1/10}$, and let $\beta$ have pairwise unitary support, as in Theorem~\ref{thm:Cw}. We have to show that
\begin{align}
	\Cw(\beta)
	\ll_\eps
	\nn{\beta}^2\nn{\nu}^2\,M^{15/8+\eps}A^{3/4}.
	\label{eq:bound_Cw_again}
\end{align}
We shall use the trivial bound
\begin{align}
	\Cw(\beta)
	\le
	\sum_{M/2<m<3M}\Big(\sum_{n\sim N}|\beta_n|\sum_{a\in\Acal}|\nu_a|\Big)^2
	\le
	3M(N+1)|\Acal|\,\nn{\beta}^2\nn{\nu}^2
	\le
	9AMN\nn{\beta}^2\nn{\nu}^2.
	\label{eq:Cw_trivial}
\end{align}
Here the first inequality follows from $|V_n(m)|\le\sum_{a\in\Acal}|\nu_a|$; the second from the Cauchy-Schwarz inequality applied to the sums over $n$ and over $a$, together with the counts $\#\{m : M/2<m<3M\}\le 5M/2+1\le 3M$ and $\#\{n : n\sim N\}\le N+1$; and the third from $N+1\le 3N/2$, which holds since $N\ge 2$, and $|\Acal|\le\lfloor 2A\rfloor\le 2A$.

\subsubsection{The case \texorpdfstring{$N\le 2A$}{N <= 2A}}
Theorem~\ref{thm:Dw} does not apply in this case, since the condition $A<N/2$ in \eqref{eq:blockrange} fails. The trivial bound suffices: by \eqref{eq:Cw_trivial} and $N\le 2A$,
\[
	\Cw(\beta)
	\le
	18A^2M\nn{\beta}^2\nn{\nu}^2
	\le
	18\,M^{15/8}A^{3/4}\nn{\beta}^2\nn{\nu}^2,
\]
because $A^{5/4}\le M^{7/8}$ is equivalent to $A\le M^{7/10}$, which holds since $A\le M^{1/10}$.

\subsubsection{The case \texorpdfstring{$N>2A$}{N > 2A}}
We now apply Theorem~\ref{thm:Dw} together with \eqref{eq:amplification}. By \eqref{eq:amplification} and \eqref{eq:bound_Dw}, if $L\ge\max(L_0,8)$ and \eqref{eq:blockrange} holds, then
\begin{align}
	\Cw(\beta)
	\ll_\eps
	\nn{\beta}^2\nn{\nu}^2
	\Big(\frac{M^2}{L}+\frac{A^{1/2}MN}{L}+ALN^{7/4}+AL^2N^{3/2}+AL^{25/16}M^{7/64}N^{3/2}\Big)M^{\eps},
	\label{eq:fivebudget}
\end{align}
where the factor $(\log L)^2$ from \eqref{eq:amplification} has been absorbed into $M^\eps$, since $L\le M$. We simplify the five terms, using $L\le M^{1/4}$, which follows from \eqref{eq:blockrange}. Since $A\ge1/2$ and $N\le M$,
\[
	\frac{M^2}{L}+\frac{A^{1/2}MN}{L}
	\le
	\frac{(\sqrt{2}+1)A^{1/2}M^2}{L}.
\]
Since $N\le M$ and $L\le M^{1/4}$, the last three terms satisfy
\[
	ALN^{7/4}\le ALM^{7/4},
	\qquad
	AL^2N^{3/2}\le ALM^{7/4},
	\qquad
	AL^{25/16}M^{7/64}N^{3/2}\le ALM^{7/4},
\]
the last one because $L^{9/16}M^{7/64}\le M^{9/64+7/64}=M^{1/4}$. Hence
\begin{align}
	\Cw(\beta)
	\ll_\eps
	\nn{\beta}^2\nn{\nu}^2
	\Big(\frac{A^{1/2}M^2}{L}+ALM^{7/4}\Big)M^{\eps}.
	\label{eq:balancepair}
\end{align}
The first term decreases and the second increases with $L$; they are equal for $L=M^{1/8}A^{-1/4}$. We take three times this value,
\begin{align}
	L=3M^{1/8}A^{-1/4},
	\label{eq:amplifierchoice}
\end{align}
the factor $3$ serving only to guarantee $L>2A$ below. With this choice,
\[
	\frac{A^{1/2}M^2}{L}=\frac{1}{3}M^{15/8}A^{3/4}
	\quadand
	ALM^{7/4}=3M^{15/8}A^{3/4},
\]
so that \eqref{eq:balancepair} gives \eqref{eq:bound_Cw_again}.

It remains to check that the choice \eqref{eq:amplifierchoice} satisfies $L\ge\max(L_0,8)$ and \eqref{eq:blockrange}. Since $1/2\le A\le M^{1/10}$, we have
\begin{align}
	L\ge 3M^{1/8-1/40}=3M^{1/10},
	\qquad
	\frac{L}{A}=3M^{1/8}A^{-5/4}\ge 3,
	\qquad
	L^4=81M^{1/2}A^{-1}\le 162M^{1/2}.
	\label{eq:Lmargins}
\end{align}
The second inequality gives $L\ge 3A>2A$, and $A<N/2$ holds in the present case, so that $A<\min(N,L)/2$. The third inequality gives $M\ge 2^{26}L^4$ as soon as $M^{1/2}\ge 162\cdot 2^{26}$, and the first gives $L\ge\max(L_0,8)$ and $L\ge 2\log(2M)$ as soon as $M\ge M_0$ for a suitable absolute constant $M_0$. Thus there is an absolute constant $M_1$ such that all the required conditions hold for $M\ge M_1$, uniformly in $N$ and $A$. For $M<M_1$ we use \eqref{eq:Cw_trivial} instead: then $A\le M_1^{1/10}$ and $N\le M<M_1$, so that $9AMN\ll_{M_1}1\le M^{15/8}A^{3/4}$. This completes the deduction of Theorem~\ref{thm:Cw} from Theorem~\ref{thm:Dw}.

\subsection{An operator inequality}
\label{sec:absorption}
By Sections~\ref{sec:elementary_estimates}-\ref{sec:Dw_implies_Cw}, it remains to bound the amplified second moment $\Dw(\beta)$, that is, to prove Theorem~\ref{thm:Dw}. In \cite{DFI} and \cite{BC}, the sum $\Dw(\beta)$ is expanded by the orthogonality of characters and split into a diagonal part and an off-diagonal part, which are estimated separately, the congruence produced by the characters being exploited in the latter; the terms in which an amplifier prime divides one of the indices $n_1,n_2$ require a separate treatment in both papers. We organize this bookkeeping differently. The following lemma, which is purely algebraic, bounds the norm of an amplified operator by its off-diagonal part restricted to indices coprime to the amplifier primes, at the cost of an error which is controlled by the norm of the unamplified operator.

In the application, $\mathcal E$ is the space of sequences $\beta$ supported on a given set of indices, $\mathsf T\beta$ is the inner sum in \eqref{eq:def_Dw} before amplification, $U_{p}$ is multiplication by $\chi(p)$, and $E_{p}$ keeps the indices divisible by $p$. The quantity $H(v)$ below then collects the terms with $p_1\ne p_2$ in which both indices are coprime to $p_1p_2$.
\begin{lemma}
\label{lem:absorb}
Let $\mathcal E$ and $\mathscr H$ be complex Hilbert spaces, and let $\mathsf T:\mathcal E\to\mathscr H$ be a bounded linear map with operator norm satisfying $\|\mathsf T\|^2\le\mathcal T$. Fix a finite set
$\Pcal$.  For each $p\in\Pcal$, let $E_p:\mathcal E\to\mathcal E$ be an orthogonal projection and $U_p:\mathscr H\to\mathscr H$ a bounded linear map with $\|U_p\|\le1$, and assume that the projections $E_p$, $p\in\Pcal$, commute with one another and satisfy
\begin{align}
	\sum_{p\in\Pcal}\|E_pv\|^2
	\le
	J\|v\|^2
	\qquad(v\in\mathcal E)
	\label{eq:absorb_J}
\end{align}
for some $J\ge0$. Let $I$ denote the identity on $\mathcal E$, and set
\begin{align}
	\mathsf A:=\sum_{p\in\Pcal}U_p\mathsf T,
	\qquad
	\mathsf T_{p_1,p_2}:=\mathsf T(I-E_{p_1})(I-E_{p_2}),
	\qquad
	H(v):=\sum_{\substack{p_1,p_2\in\Pcal\\p_1\ne p_2}}\langle U_{p_1}\mathsf T_{p_1,p_2}v,\,U_{p_2}\mathsf T_{p_1,p_2}v\rangle_{\mathscr H}.
	\label{eq:absorb_defs}
\end{align}
If $\mathcal G\ge0$ is such that $|H(v)|\le\mathcal G\|v\|^2$ for all $v\in\mathcal E$, then
\begin{align}
	\|\mathsf Av\|^2
	\le
	\big(4\mathcal G+4\mathcal T|\Pcal|+14\mathcal T|\Pcal|J\big)\|v\|^2
	\qquad(v\in\mathcal E).
	\label{eq:absorption}
\end{align}
\end{lemma}
\begin{proof}
The inner product is linear in the first argument. Set
\[
	\mathsf A_0:=\sum_{p\in\Pcal}U_p\mathsf T(I-E_p),
\]
so that $\mathsf A-\mathsf A_0=\sum_{p}U_p\mathsf TE_p$. Let $v\in\mathcal E$. By the inequality $\|\sum_{i\in\mathcal I}u_i\|^2\le|\mathcal I|\sum_{i\in\mathcal I}\|u_i\|^2$, together with $\|U_p\|\le1$, $\|\mathsf T\|^2\le\mathcal T$ and \eqref{eq:absorb_J},
\[
	\|(\mathsf A-\mathsf A_0)v\|^2
	\le
	|\Pcal|\sum_{p\in\Pcal}\|U_p\mathsf TE_pv\|^2
	\le
	\mathcal T|\Pcal|\sum_{p\in\Pcal}\|E_pv\|^2
	\le
	\mathcal T|\Pcal|J\|v\|^2.
\]
Hence, by $\|u+u'\|^2\le2\|u\|^2+2\|u'\|^2$,
\begin{align}
	\|\mathsf Av\|^2
	\le
	2\|\mathsf A_0v\|^2+2\mathcal T|\Pcal|J\|v\|^2.
	\label{eq:A_vs_A0}
\end{align}
It remains to bound $\|\mathsf A_0\|^2$, which is finite since $\mathsf T$ is bounded, $\|U_p\|\le1$, $\|I-E_p\|\le1$ and $\Pcal$ is finite. Set
\[
	a_{p_1}:=U_{p_1}\mathsf T(I-E_{p_1})v,
	\qquad
	x_{p_1p_2}:=U_{p_1}\mathsf T(I-E_{p_1})E_{p_2}v,
	\qquad
	y_{p_1p_2}:=U_{p_2}\mathsf T(I-E_{p_2})E_{p_1}v.
\]
Since the projections commute, $a_{p_1}-x_{p_1p_2}=U_{p_1}\mathsf T_{p_1,p_2}v$ and $a_{p_2}-y_{p_1p_2}=U_{p_2}\mathsf T_{p_1,p_2}v$, so that
\begin{align}
	H(v)
	=
	\sum_{\substack{p_1,p_2\in\Pcal\\p_1\ne p_2}}
	\langle a_{p_1}-x_{p_1p_2},\,a_{p_2}-y_{p_1p_2}\rangle_{\mathscr H}.
	\label{eq:H_as_apaq}
\end{align}
Expanding the norm of $\mathsf A_0v=\sum_{p_1}a_{p_1}$ gives
\begin{align}
	\|\mathsf A_0v\|^2
	=
	\sum_{p_1\in\Pcal}\|a_{p_1}\|^2
	+
	\sum_{\substack{p_1,p_2\in\Pcal\\p_1\ne p_2}}\langle a_{p_1},a_{p_2}\rangle_{\mathscr H},
	\label{eq:A0_expanded}
\end{align}
and expanding $\langle a_{p_1},a_{p_2}\rangle=\langle (a_{p_1}-x_{p_1p_2})+x_{p_1p_2},\,(a_{p_2}-y_{p_1p_2})+y_{p_1p_2}\rangle$ gives
\[
	\langle a_{p_1},a_{p_2}\rangle_{\mathscr H}
	=
	\langle a_{p_1}-x_{p_1p_2},\,a_{p_2}-y_{p_1p_2}\rangle_{\mathscr H}
	+\langle x_{p_1p_2},a_{p_2}\rangle_{\mathscr H}
	+\langle a_{p_1},y_{p_1p_2}\rangle_{\mathscr H}
	-\langle x_{p_1p_2},y_{p_1p_2}\rangle_{\mathscr H}.
\]
We sum this identity over $p_1\ne p_2$ and set
\[
	S:=\sum_{\substack{p_1,p_2\in\Pcal\\p_1\ne p_2}}\langle x_{p_1p_2},a_{p_2}\rangle_{\mathscr H},
	\qquad
	Z:=\sum_{\substack{p_1,p_2\in\Pcal\\p_1\ne p_2}}\langle x_{p_1p_2},y_{p_1p_2}\rangle_{\mathscr H}.
\]
Since $y_{p_1p_2}=x_{p_2p_1}$, interchanging the names of $p_1$ and $p_2$ shows that $\sum_{p_1\ne p_2}\langle a_{p_1},y_{p_1p_2}\rangle=\overline{S}$ and that $Z$ is real. Together with \eqref{eq:H_as_apaq} and \eqref{eq:A0_expanded}, we obtain
\begin{align}
	\|\mathsf A_0v\|^2
	=
	\sum_{p_1\in\Pcal}\|a_{p_1}\|^2+H(v)+2\Re S-Z.
	\label{eq:A0_identity}
\end{align}
We now bound the terms on the right. First, $\|a_{p_1}\|^2\le\mathcal T\|v\|^2$ for each $p_1$, so that
\[
	\sum_{p_1\in\Pcal}\|a_{p_1}\|^2\le\mathcal T|\Pcal|\,\|v\|^2.
\]
Next, $\|x_{p_1p_2}\|^2\le\mathcal T\|E_{p_2}v\|^2$, and therefore, by \eqref{eq:absorb_J},
\[
	\sum_{\substack{p_1,p_2\in\Pcal\\p_1\ne p_2}}\|x_{p_1p_2}\|^2
	\le
	\mathcal T|\Pcal|\sum_{p_2\in\Pcal}\|E_{p_2}v\|^2
	\le
	\mathcal T|\Pcal|J\|v\|^2,
\]
and the same bound holds for the $y_{p_1p_2}$, since $y_{p_1p_2}=x_{p_2p_1}$. By the Cauchy-Schwarz inequality, $|Z|\le\mathcal T|\Pcal|J\|v\|^2$. For $S$, we use $(I-E_{p_2})E_{p_2}=0$ to write
\[
	\sum_{\substack{p_1\in\Pcal\\p_1\ne p_2}}x_{p_1p_2}
	=
	\sum_{p_1\in\Pcal}U_{p_1}\mathsf T(I-E_{p_1})E_{p_2}v
	=
	\mathsf A_0(E_{p_2}v),
\]
so that $S=\sum_{p_2}\langle\mathsf A_0(E_{p_2}v),a_{p_2}\rangle_{\mathscr H}$, and the Cauchy-Schwarz inequality gives
\[
	|S|
	\le
	\Big(\sum_{p_2\in\Pcal}\|\mathsf A_0E_{p_2}v\|^2\Big)^{1/2}
	\Big(\sum_{p_2\in\Pcal}\|a_{p_2}\|^2\Big)^{1/2}
	\le
	\big(\|\mathsf A_0\|^2J\|v\|^2\big)^{1/2}\big(\mathcal T|\Pcal|\,\|v\|^2\big)^{1/2}.
\]
Inserting these bounds and the hypothesis $|H(v)|\le\mathcal G\|v\|^2$ into \eqref{eq:A0_identity} and taking the supremum over all unit vectors $v\in\mathcal E$, we get
\[
	\|\mathsf A_0\|^2
	\le
	\mathcal G+\mathcal T|\Pcal|(1+J)+2\sqrt{\|\mathsf A_0\|^2\,\mathcal T|\Pcal|J}
	\le
	\mathcal G+\mathcal T|\Pcal|(1+J)+\frac{1}{2}\|\mathsf A_0\|^2+2\mathcal T|\Pcal|J,
\]
where we used $2ab\le a^2/2+2b^2$ with $a=\|\mathsf A_0\|$ and $b=(\mathcal T|\Pcal|J)^{1/2}$. Hence
\[
	\|\mathsf A_0\|^2
	\le
	2\mathcal G+2\mathcal T|\Pcal|+6\mathcal T|\Pcal|J,
\]
and \eqref{eq:absorption} follows from \eqref{eq:A_vs_A0}.
\end{proof}

\subsection{Application to the amplified second moment}
\label{sec:off_diagonal}
We now cast $\Dw(\beta)$ in the setting of Lemma~\ref{lem:absorb}. Let $\beta$ be as in Theorem~\ref{thm:Dw}, and let $\mathcal I\subseteq(N,2N]$ be the support of $\beta$. Since $\beta$ has pairwise unitary support, so does every sequence supported on a subset of $\mathcal I$. Let $\mathcal E$ be the space of sequences $v=(v_n)_{n\sim N}$ with $v_n=0$ for $n\notin\mathcal I$, with the norm $\nn{v}$. Let $\mathscr H$ be the space of complex-valued functions $F(m,\chi)$ on the finite set of pairs $(m,\chi)$ with $m\in\N$, $w(m/M)>0$, and $\chi$ a Dirichlet character modulo $m$, with the norm
\begin{align}
	\|F\|_{\mathscr H}^2
	:=
	\sum_{m\in\N}\frac{w(m/M)}{\varphi(m)}\sum_{\chi\bmod m}|F(m,\chi)|^2.
	\label{eq:def_H_norm}
\end{align}
Both are finite-dimensional Hilbert spaces. For every sequence $v$ supported on $n\sim N$, for $F\in\mathscr H$ and for $p\in\Pcal$, define
\begin{align}
	(\mathsf Tv)(m,\chi):=\sum_{\substack{n\sim N\\(m,n)=1}}v_n\chi(n)V_n(m),
	\qquad
	(U_pF)(m,\chi):=\chi(p)F(m,\chi),
	\qquad
	(E_pv)_n:=\ind_{p\mid n}v_n.
	\label{eq:def_T_U_E}
\end{align}
Then $\mathsf T v\in\mathscr H$, $U_pF\in\mathscr H$, and $E_p$ maps $\mathcal E$ into itself; we regard $\mathsf T$ as a linear map $\mathcal E\to\mathscr H$, $U_p$ as a linear map $\mathscr H\to\mathscr H$ and $E_p$ as a linear map $\mathcal E\to\mathcal E$. Note that $\beta\in\mathcal E$. For $F=\sum_{p\in\Pcal}U_p\mathsf T\beta$, the definitions \eqref{eq:def_T_U_E} give
\[
	F(m,\chi)
	=
	\sum_{p\in\Pcal}\chi(p)\,(\mathsf T\beta)(m,\chi)
	=
	\bigg(\sum_{p\in\Pcal}\chi(p)\bigg)\sum_{\substack{n\sim N\\(m,n)=1}}\beta_n\chi(n)V_n(m),
\]
and hence, by \eqref{eq:def_H_norm},
\[
	\|F\|_{\mathscr H}^2
	=
	\sum_{m\in\N}\frac{w(m/M)}{\varphi(m)}\sum_{\chi\bmod m}
	\bigg|\sum_{p\in\Pcal}\chi(p)\bigg|^2
	\bigg|\sum_{\substack{n\sim N\\(m,n)=1}}\beta_n\chi(n)V_n(m)\bigg|^2,
\]
where the terms with $w(m/M)=0$ may be included since they vanish. The right-hand side is $\Dw(\beta)$ by \eqref{eq:def_Dw}. Thus
\begin{align}
	\Dw(\beta)
	=
	\Big\|\sum_{p\in\Pcal}U_p\mathsf T\beta\Big\|_{\mathscr H}^2
	=
	\|\mathsf A\beta\|_{\mathscr H}^2
	\label{eq:Dw_as_norm}
\end{align}
in the notation \eqref{eq:absorb_defs}.

We verify the hypotheses of Lemma~\ref{lem:absorb}. The set $\Pcal$ is finite with $|\Pcal|\le L$. Since $|\chi(p)|\le1$, each $U_p$ satisfies $\|U_pF\|_{\mathscr H}\le\|F\|_{\mathscr H}$. Each $E_p$ is an orthogonal projection of $\mathcal E$, namely the projection onto the sequences supported on the multiples of $p$, and these projections commute, since $E_{p_1}E_{p_2}v=(\ind_{p_1\mid n}\ind_{p_2\mid n}v_n)_n=E_{p_2}E_{p_1}v$. Moreover, an integer $n\le2N$ has at most $\log(2N)/\log L$ prime factors exceeding $L$, so that
\begin{align}
	\sum_{p\in\Pcal}\|E_pv\|^2
	=
	\sum_{n\sim N}|v_n|^2\,\#\{p\in\Pcal : p\mid n\}
	\le
	\frac{\log(2N)}{\log L}\nn{v}^2,
	\label{eq:J_value}
\end{align}
which is \eqref{eq:absorb_J} with $J=\log(2N)/\log L$. The norm of $\mathsf T$ is bounded by the following lemma.
\begin{lemma}
\label{lem:unamplified}
Let $N\le M$. For every sequence $v$ supported on $n\sim N$,
\begin{align}
	\|\mathsf Tv\|_{\mathscr H}^2
	\le
	9(M+N)\Big(\frac{A}{N}+1\Big)\nn{v}^2\nn{\nu}^2.
	\label{eq:unamplified}
\end{align}
\end{lemma}
\begin{proof}
By the orthogonality relation
\begin{align}
	\frac{1}{\varphi(m)}\sum_{\chi\bmod m}\chi(n)\overline{\chi(n')}
	=
	\ind_{(nn',m)=1}\ind_{n\equiv n'\bmod m},
	\label{eq:orthogonality}
\end{align}
we have
\[
	\|\mathsf Tv\|_{\mathscr H}^2
	=
	\sum_{m\in\N}w\Big(\frac{m}{M}\Big)
	\sum_{\substack{n,n'\sim N\\(nn',m)=1\\n\equiv n'\bmod m}}
	v_n\overline{v_{n'}}V_n(m)\overline{V_{n'}(m)}
	=
	\sum_{m\in\N}w\Big(\frac{m}{M}\Big)
	\sum_{\substack{r\bmod m\\(r,m)=1}}
	\bigg|\sum_{\substack{n\sim N\\n\equiv r\bmod m}}v_nV_n(m)\bigg|^2.
\]
On the support of $w(m/M)$ we have $m>M/2\ge N/2$, so that a residue class modulo $m$ contains at most $N/m+1<3$ integers $n\sim N$. By the Cauchy-Schwarz inequality applied to the inner sum,
\begin{align}
	\|\mathsf Tv\|_{\mathscr H}^2
	\le
	3\sum_{n\sim N}|v_n|^2\sum_{\substack{M/2<m<3M\\(m,n)=1}}|V_n(m)|^2.
	\label{eq:T_norm_step}
\end{align}
We bound the inner sum. By \eqref{eq:def_Vn}, $V_n(m)$ depends only on the residue class of $\overline{m}$ modulo $n$, and as $m$ runs over the units modulo $n$ so does $\overline{m}$. Each residue class modulo $n$ contains at most $3M/n+1$ integers $m\in(M/2,3M)$. Hence, by Parseval's identity for the additive characters modulo $n$,
\[
	\sum_{\substack{M/2<m<3M\\(m,n)=1}}|V_n(m)|^2
	\le
	\Big(\frac{3M}{n}+1\Big)\sum_{r\bmod n}\bigg|\sum_{a\in\Acal}\nu_a\,\e\Big(\frac{ar}{n}\Big)\bigg|^2
	=
	\Big(\frac{3M}{n}+1\Big)\,n\sum_{s\bmod n}\bigg|\sum_{\substack{a\in\Acal\\a\equiv s\bmod n}}\nu_a\bigg|^2.
\]
A residue class modulo $n$ contains at most $A/n+1$ integers $a\in\Acal$, so by the Cauchy-Schwarz inequality the last sum over $s$ is at most $(A/n+1)\nn{\nu}^2$. Since $n>N$, we obtain
\[
	\sum_{\substack{M/2<m<3M\\(m,n)=1}}|V_n(m)|^2
	\le
	(3M+n)\Big(\frac{A}{n}+1\Big)\nn{\nu}^2
	\le
	3(M+N)\Big(\frac{A}{N}+1\Big)\nn{\nu}^2,
\]
and inserting this into \eqref{eq:T_norm_step} proves \eqref{eq:unamplified}.
\end{proof}
Under the hypotheses \eqref{eq:blockrange} we have $N\le M$ and $A<N/2$, so that Lemma~\ref{lem:unamplified} gives $\|\mathsf T\|^2\le\mathcal T$ with $\mathcal T=27M\nn{\nu}^2$. With $|\Pcal|\le L$, $J=\log(2N)/\log L\le\log(2M)/\log 2$ and $L\le M$, Lemma~\ref{lem:absorb} and \eqref{eq:Dw_as_norm} therefore yield
\begin{align}
	\Dw(\beta)
	\le
	4\mathcal G\nn{\beta}^2+\big(4\mathcal T L+14\mathcal T LJ\big)\nn{\beta}^2
	\ll
	\mathcal G\nn{\beta}^2+ML\nn{\beta}^2\nn{\nu}^2M^{\eps},
	\label{eq:Dw_after_absorption}
\end{align}
provided that $|H(v)|\le\mathcal G\nn{v}^2$ for every $v\in\mathcal E$. The second term in \eqref{eq:Dw_after_absorption} accounts for the first term in \eqref{eq:bound_Dw}. It remains to identify $H$ and to bound it. We write out $H(\beta)$; the same expression holds for every $v\in\mathcal E$ in place of $\beta$.

By \eqref{eq:def_T_U_E}, $\big((I-E_{p_1})(I-E_{p_2})\beta\big)_n=\ind_{(p_1p_2,n)=1}\beta_n$, so that
\[
	(\mathsf T_{p_1,p_2}\beta)(m,\chi)
	=
	\sum_{\substack{n\sim N\\(m,n)=(p_1p_2,n)=1}}\beta_n\chi(n)V_n(m).
\]
Inserting this and \eqref{eq:def_H_norm} into the definition of $H$ in \eqref{eq:absorb_defs}, we get
\[
	H(\beta)
	=
	\sum_{\substack{p_1,p_2\in\Pcal\\p_1\ne p_2}}
	\sum_{m\in\N}w\Big(\frac{m}{M}\Big)
	\sum_{\substack{n_1,n_2\sim N\\(m,n_1n_2)=(p_1p_2,n_1n_2)=1}}
	\beta_{n_1}\overline{\beta_{n_2}}V_{n_1}(m)\overline{V_{n_2}(m)}
	\frac{1}{\varphi(m)}\sum_{\chi\bmod m}\chi(p_1n_1)\overline{\chi(p_2n_2)}
\]
and by \eqref{eq:orthogonality} the last sum over $\chi$ equals $1$ if $(m,p_1p_2)=1$ and $p_1n_1\equiv p_2n_2\bmod m$, and vanishes otherwise. Expanding $V_{n_1}(m)\overline{V_{n_2}(m)}$ by \eqref{eq:def_Vn}, we conclude that
\begin{align}
	H(\beta)
	=
	\sum_{\substack{p_1,p_2\in\Pcal\\p_1\ne p_2}}
	\sum_{\substack{m\in\N\\(m,p_1p_2)=1}}w\Big(\frac{m}{M}\Big)
	\sum_{\substack{n_1,n_2\sim N\\(m,n_1n_2)=(p_1p_2,n_1n_2)=1\\p_1n_1\equiv p_2n_2\bmod m}}
	\beta_{n_1}\overline{\beta_{n_2}}
	\sumtwo_{a_1,a_2\in\Acal}\nu_{a_1}\overline{\nu_{a_2}}\,
	\e\Big(\frac{a_1\overline{m}}{n_1}-\frac{a_2\overline{m}}{n_2}\Big).
	\label{eq:H_explicit}
\end{align}
These are the off-diagonal terms of $\Dw$: the amplifier primes are distinct and coprime to the indices, and the congruence $p_1n_1\equiv p_2n_2 \bmod m$ relates the two indices through the modulus.

\subsection{Theorem~\ref{thm:H} implies Theorem~\ref{thm:Dw}}
\label{sec:H_implies_Dw}
Let $M$, $N$, $A$, $L$, $\beta$ and $\nu$ be as in Theorem~\ref{thm:Dw}, and let $\mathcal E$ be as in Section~\ref{sec:off_diagonal}. By \eqref{eq:Dw_after_absorption}, it suffices to find $\mathcal G$ with $|H(v)|\le\mathcal G\nn{v}^2$ for every $v\in\mathcal E$. We bound $H(\beta)$ by decomposing \eqref{eq:H_explicit} according to the value of $g=(n_1,n_2)$, which brings it to the form \eqref{eq:def_Hg}; the argument applies verbatim to every $v\in\mathcal E$, and this gives the required uniform bound.

In \eqref{eq:H_explicit}, write $n_1=gn_1'$ and $n_2=gn_2'$ with $g=(n_1,n_2)$, so that $(n_1',n_2')=1$ and $g\le 2N$. The conditions $(m,n_1n_2)=1$ and $(p_1p_2,n_1n_2)=1$ are equivalent to $(mp_1p_2,g)=1$ and $(mp_1p_2,n_1'n_2')=1$, and since $(m,g)=1$ the congruence $p_1gn_1'\equiv p_2gn_2'\bmod m$ is equivalent to $p_1n_1'\equiv p_2n_2'\bmod m$. For $g\le2N$ define the sequence $\beta^g$ supported on $n\sim N/g$ by
\begin{align}
	\beta^g_n:=\beta_{gn}\ind_{(n,g)=1}.
	\label{eq:def_vg}
\end{align}
If $\beta_{gn_1'}\ne0$ and $\beta_{gn_2'}\ne0$ with $(n_1',n_2')=1$, then $g=(gn_1',gn_2')$ is a unitary divisor of $gn_1'$ and of $gn_2'$, because $\beta$ has pairwise unitary support; that is, $(g,n_1')=(g,n_2')=1$. Hence, for every pair $n_1,n_2$ that contributes to \eqref{eq:H_explicit}, we have $\beta_{n_1}\overline{\beta_{n_2}}=\beta^g_{n_1'}\overline{\beta^g_{n_2'}}$, while every pair $n_1',n_2'$ with $(g,n_1'n_2')>1$ contributes nothing to either side. Comparing with \eqref{eq:def_Hg}, we obtain
\begin{align}
	H(\beta)
	=
	\sum_{g\le 2N}\Hcal_g(\beta^g,N/g).
	\label{eq:H_decomposed}
\end{align}
We check that $\beta^g$ satisfies the hypotheses of Theorem~\ref{thm:H}. It is supported on $n\sim N/g$, and $\beta^g_n\ne0$ implies $(n,g)=1$ by \eqref{eq:def_vg}. It has pairwise unitary support: if $\beta^g_n\ne0$ and $\beta^g_u\ne0$, then $\beta_{gn}\ne0$ and $\beta_{gu}\ne0$, so that $(gn,gu)=g(n,u)$ is a unitary divisor of $gn$ and of $gu$; thus $(g(n,u),n/(n,u))=1$, and in particular $((n,u),n/(n,u))=1$, and likewise for $u$. Since \eqref{eq:blockrange} holds, Theorem~\ref{thm:H} applies to each term of \eqref{eq:H_decomposed}. All four powers of $g$ in \eqref{eq:bound_H} are negative, and $g\ge1$, so that
\begin{align}
	|H(\beta)|
	\ll_\eps
	\nn{\nu}^2
	\Big(A^{1/2}LN+\frac{AL^3N^{7/4}}{M}+\frac{AL^4N^{3/2}}{M}+\frac{AL^{57/16}N^{3/2}}{M^{57/64}}\Big)M^{\eps}
	\sum_{g\le2N}\nn{\beta^g}^2.
	\label{eq:H_summed_g}
\end{align}
Finally, since each $n'\sim N$ has $\tau(n')$ representations $n'=gn$,
\[
	\sum_{g\le2N}\nn{\beta^g}^2
	=
	\sum_{g\le 2N}\sum_{\substack{n\sim N/g\\(n,g)=1}}|\beta_{gn}|^2
	\le
	\sum_{n'\sim N}\tau(n')|\beta_{n'}|^2
	\ll_\eps
	N^{\eps}\nn{\beta}^2
\]
by the divisor bound. Hence $|H(\beta)|\le\mathcal G\nn{\beta}^2$ with
\[
	\mathcal G
	=
	c(\eps)\nn{\nu}^2
	\Big(A^{1/2}LN+\frac{AL^3N^{7/4}}{M}+\frac{AL^4N^{3/2}}{M}+\frac{AL^{57/16}N^{3/2}}{M^{57/64}}\Big)M^{2\eps},
\]
where $c(\eps)$ depends only on $\eps$. Every $v\in\mathcal E$ is supported on $n\sim N$ and has pairwise unitary support, and nothing else about $\beta$ was used; so the same argument gives $|H(v)|\le\mathcal G\nn{v}^2$ for every $v\in\mathcal E$, with the same $\mathcal G$. Inserting this into \eqref{eq:Dw_after_absorption} and renaming $\eps$, we get
\[
	\Dw(\beta)
	\ll_\eps
	\nn{\beta}^2\nn{\nu}^2
	\Big(ML+A^{1/2}LN+\frac{AL^3N^{7/4}}{M}+\frac{AL^4N^{3/2}}{M}+\frac{AL^{57/16}N^{3/2}}{M^{57/64}}\Big)M^{\eps},
\]
which is \eqref{eq:bound_Dw} after the factor $L^2/M$ has been taken out, since $7/64-1=-57/64$. This proves Theorem~\ref{thm:Dw}, and with it Theorems~\ref{thm:Cw} and \ref{thm:main}, subject to Theorem~\ref{thm:H}.

\section{The correlation sum}
\label{sec:correlation}
In this section we prove Theorem~\ref{thm:H}. Throughout,
$M$, $N$, $A$ and $L$ satisfy \eqref{eq:blockrange}, and $g\in\N$.
The sequence $\beta$ is supported on $n\sim N/g$, has pairwise
unitary support and satisfies $\beta_n\ne0\Rightarrow(n,g)=1$;
the sequence $\nu$ is supported on $a\sim A$.
We write $\Hcal_g$ for $\Hcal_g(\beta,N/g)$.

Following \cite{DFI} and \cite{BC}, we replace the summation over $m$
by a summation over the complementary divisor
$d=(p_1n_1-p_2n_2)/m$. The main difference is that we exploit
cancellation among the resulting Kloosterman sums as the modulus
varies, using Lemma~\ref{lem:progression}. A key arithmetic feature
is that, in the configurations where this estimate is applied,
the numerator is coprime to the amplifier primes. These primes
therefore enlarge the spectral level without enlarging its
greatest common divisor with the numerator, which limits the
contribution of the exceptional spectrum.

\subsection{The complementary divisor and the residue classes}
\label{sec:complementary_divisor}
We first replace the summation over $m$ by a summation over its
complementary divisor. In each term of \eqref{eq:def_Hg}, with
$N/g$ in place of $N$, the conditions
\[
(n_1,n_2)=1,\qquad (p_1p_2,n_1n_2)=1,\qquad p_1\ne p_2
\]
give $(p_1n_1,p_2n_2)=1$. In particular, $p_1n_1\ne p_2n_2$. Write
\begin{align}
	p_1n_1-p_2n_2=dm.
	\label{eq:complementary_divisor}
\end{align}
Then $d$ is a nonzero integer. Moreover,
\begin{align}
	(d,p_1p_2n_1n_2)=1,
	\label{eq:d_coprime}
\end{align}
since, by \eqref{eq:complementary_divisor}, a prime dividing $d$ and one of $p_1n_1$, $p_2n_2$ would divide the other, which is impossible. 
The support of $w(m/M)$ gives $M/2<m<3M$, while $n_1,n_2\le 2N/g$ and $p_1,p_2\le 2L$ give $|p_1n_1-p_2n_2|<4LN/g$. Hence,
\begin{align}
	1\le|d|<D:=\frac{8LN}{gM}.
	\label{eq:def_D}
\end{align}
In particular, $\Hcal_g=0$ unless $N/g>M/(8L)$. Therefore, we assume 
\begin{align}
\frac{N}{g}>\frac{M}{8L}\label{eq:assumption on N over g}
\end{align}
from now on.

Conversely, fix $p_1,p_2,n_1,n_2$ satisfying
the original coprimality conditions. Each nonzero divisor $d$
of $p_1n_1-p_2n_2$ determines
\[
	m:=\frac{p_1n_1-p_2n_2}{d}.
\] 
The same argument as above gives $(m,p_1p_2n_1n_2)=1$. The remaining conditions on $m$ in \eqref{eq:def_Hg} are $(m,g)=1$ and the support of $w(m/M)$, which forces $m>0$.

We impose $(m,g)=1$ by writing $m\equiv c\bmod g$ with
$(c,g)=1$. The conditions $d\mid p_1n_1-p_2n_2$ and $(p_1n_1-p_2n_2)/d\equiv c\bmod g$ together are equivalent to
\begin{align}
	p_1n_1-p_2n_2\equiv cd\bmod g|d|.
	\label{eq:affine_progression}
\end{align}
Therefore
\begin{align}
	\Hcal_g
	=&
	\sum_{\substack{n_1,n_2\sim N/g\\(n_1,n_2)=(g,n_1n_2)=1}}\beta_{n_1}\overline{\beta_{n_2}}
	\sum_{\substack{c\bmod g\\(c,g)=1}}
	\sum_{0<|d|<D}
	\sum_{\substack{p_1,p_2\in\Pcal,\ p_1\ne p_2\\(p_1p_2,gdn_1n_2)=1\\p_1n_1-p_2n_2\equiv cd\bmod g|d|}}
	w\Big(\frac{p_1n_1-p_2n_2}{dM}\Big)
	\notag\\
	&\quad\times
	\sumtwo_{a_1,a_2\in\Acal}\nu_{a_1}\overline{\nu_{a_2}}\,
	\e\Big(\frac{a_1\overline{m}}{gn_1}-\frac{a_2\overline{m}}{gn_2}\Big),
	\label{eq:Hg_with_cd}
\end{align}
where $m=(p_1n_1-p_2n_2)/d$ in the exponential. The condition $(d,n_1n_2)=1$ has been dropped, as it is implied by \eqref{eq:d_coprime}.

We next split $n_1$ and $n_2$ into reduced residue classes
$\rho,\sigma\bmod g$. On the classes
$n_1\equiv\rho\bmod g$ and $n_2\equiv\sigma\bmod g$,
\eqref{eq:affine_progression} requires
\begin{align}
	p_1\rho-p_2\sigma\equiv cd\bmod g.
	\label{eq:residue_condition}
\end{align}
Moreover, on these classes the congruence \eqref{eq:affine_progression} can be written in a form which is homogeneous in $n_1$ and $n_2$. This will allow us to remove the restrictions
$n_1\equiv\rho$ and $n_2\equiv\sigma$ later while still expressing $n_2$ as a fixed unit multiple of $n_1$
modulo $g|d|$, as required for Lemma~\ref{lem:progression}. Let $\overline{\rho}$ be an inverse of $\rho$ modulo $g$. Since $n_1\equiv\rho\bmod g$,
\[
	g\mid\overline{\rho}n_1-1,
	\qquad
	cd\,\overline{\rho}n_1\equiv cd\bmod g|d|.
\]
Thus, \eqref{eq:affine_progression} is equivalent to
\begin{align}
	p_1n_1-p_2n_2\equiv cd\,\overline{\rho}\,n_1\bmod g|d|.
	\label{eq:homogeneous_progression}
\end{align}
Replacing $\overline{\rho}$ by $\overline{\rho}+\ell g$ for any $\ell\in\Z$ changes the right-hand side by a multiple of $g|d|$, so \eqref{eq:homogeneous_progression} is independent of the choice of $\overline{\rho}$. Unlike \eqref{eq:affine_progression}, \eqref{eq:homogeneous_progression} is preserved
when $(n_1,n_2)$ is replaced by $(\lambda n_1,\lambda n_2)$
for any integer $\lambda$.

Reducing \eqref{eq:homogeneous_progression} modulo $|d|$
and modulo $g$, and using \eqref{eq:residue_condition}
and $(p_2,g)=1$, gives
\begin{align}
	p_1n_1\equiv p_2n_2\bmod |d|
	\quad\text{and}\quad
	n_2\equiv\sigma\overline{\rho}\,n_1\bmod g.
	\label{eq:consequences}
\end{align}
These hold without restricting $n_1,n_2$ to the classes
$\rho,\sigma$. The first congruence implies $(d,n_1n_2)=1$ whenever
$(n_1,n_2)=(p_1p_2,d)=1$, and the second congruence implies
$(n_2,g)=1$ whenever $(n_1,g)=1$.

\subsection{The exponential factor}
\label{sec:exponential_factor}
In this subsection, we eliminate $m$ from the exponential in \eqref{eq:Hg_with_cd}
and rewrite the resulting expression in the exponential in terms of a single
Kloosterman fraction.

We first eliminate $m$ from  \eqref{eq:Hg_with_cd}. Since $(g,n_1)=1$ and $(m,gn_1)=1$, by Lemma~\ref{lem:crt_inverse}, 
\[
	\frac{\overline{m}}{gn_1}
	\equiv
	\frac{\overline{mn_1}}{g}+\frac{\overline{mg}}{n_1}\bmod 1,
\]
and similarly with $n_2$ in place of $n_1$. Since $m\equiv c\bmod g$,
\[
	\overline{mn_1}
	\equiv\overline{c}\,\overline{n_1}\bmod g.
\]
Moreover, \eqref{eq:complementary_divisor} gives
\[
	dm\equiv-p_2n_2\bmod n_1,
	\qquad
	dm\equiv p_1n_1\bmod n_2.
\] 
Since $d$ is invertible modulo $n_1n_2$ by \eqref{eq:d_coprime}, we obtain
\[
	\overline{mg}\equiv-d\,\overline{gp_2n_2}\bmod n_1
	\quadand
	\overline{mg}\equiv d\,\overline{gp_1n_1}\bmod n_2.
\]
Putting these together, we obtain
\begin{align}
	\frac{a_1\overline{m}}{gn_1}-\frac{a_2\overline{m}}{gn_2}
	\equiv
	\frac{\overline{c}\,(a_1\overline{n_1}-a_2\overline{n_2})}{g}
	-d\Big(\frac{a_1\overline{gp_2n_2}}{n_1}+\frac{a_2\overline{gp_1n_1}}{n_2}\Big)\bmod 1.
	\label{eq:phase_cd}
\end{align}

We next combine the two fractions in parentheses in \eqref{eq:phase_cd} into a
single Kloosterman fraction. For reduced residues $\rho,\sigma\bmod g$, define
\begin{align}
	\xi=\xi(c,d,p_1,a_1,a_2;\rho,\sigma)
	:=
	\overline{c}\,(a_1\overline{\rho}-a_2\overline{\sigma})+da_2\,\overline{p_1}\,\overline{\rho}\,\overline{\sigma},
	\label{eq:def_xi}
\end{align}
 with all inverses modulo $g$. For the primes $p_1,p_2$ in
\eqref{eq:Hg_with_cd} and coprime positive integers $n_1,n_2$ with $(gp_1p_2,n_1n_2)=1$, set
\begin{align}
	\Theta^{a_1,a_2}_{d,p_1,p_2}(n_1,n_2)
	:=
	\frac{d(a_2p_2-a_1p_1)\overline{gp_2n_2}}{p_1n_1}-\frac{da_2}{gp_1n_1n_2}.
	\label{eq:def_Theta}
\end{align}

\begin{lemma}
\label{lem:Theta_identity}
Let $g,n_1,n_2$ be positive integers with $(n_1,n_2)=(g,n_1n_2)=1$. Suppose $p_1\ne p_2$ are primes not dividing $gn_1n_2$. Then for any integers $d,a_1,a_2$,

\begin{align}
	-d\Big(\frac{a_1\overline{gp_2n_2}}{n_1}+\frac{a_2\overline{gp_1n_1}}{n_2}\Big)
	\equiv
	\frac{da_2\,\overline{p_1n_1n_2}}{g}+\Theta^{a_1,a_2}_{d,p_1,p_2}(n_1,n_2)\bmod 1.
	\label{eq:Theta_identity}
\end{align}
\end{lemma}
\begin{proof}
Apply Lemma~\ref{lem:reciprocity} with $x=gp_1n_1$ and
$y=n_2$, 
\[
	\frac{\overline{gp_1n_1}}{n_2}
	\equiv
	\frac{1}{gp_1n_1n_2}-\frac{\overline{n_2}}{gp_1n_1}
	\bmod 1.
\]
Then by Lemma~\ref{lem:crt_inverse} with $x=n_2$, $y=g$ and $z=p_1n_1$, we arrive at
\[
    \frac{\overline{gp_1n_1}}{n_2}
	\equiv\frac{1}{gp_1n_1n_2}-\frac{\overline{p_1n_1n_2}}{g}-\frac{\overline{gn_2}}{p_1n_1}	\bmod 1.
\]
Modulo $p_1n_1$, we have
\[
	\overline{gn_2}
	\equiv p_2\overline{gp_2n_2}\bmod p_1n_1.
\]
Also, an inverse of $gp_2n_2$ modulo $p_1n_1$ is an
inverse modulo $n_1$. We may therefore write
\[
	\frac{a_1\overline{gp_2n_2}}{n_1}
	\equiv
	\frac{a_1p_1\overline{gp_2n_2}}{p_1n_1}
	\bmod 1,
\]
using the inverse modulo $p_1n_1$ on the right. Substituting these identities into the left-hand side
of \eqref{eq:Theta_identity}, we obtain
\begin{align*}
	&-\frac{da_1p_1\overline{gp_2n_2}}{p_1n_1}
	-\frac{da_2}{gp_1n_1n_2}
	+\frac{da_2\overline{p_1n_1n_2}}{g}
	+\frac{da_2p_2\overline{gp_2n_2}}{p_1n_1}\\
	&\quad=
	\frac{da_2\overline{p_1n_1n_2}}{g}
	+\frac{d(a_2p_2-a_1p_1)\overline{gp_2n_2}}{p_1n_1}
	-\frac{da_2}{gp_1n_1n_2},
\end{align*}
where the last two terms are $\Theta^{a_1,a_2}_{d,p_1,p_2}(n_1,n_2)$ by \eqref{eq:def_Theta}.
\end{proof}
On the classes $n_1\equiv\rho$, $n_2\equiv\sigma\bmod g$, we have $\overline{c}(a_1\overline{n_1}-a_2\overline{n_2})+da_2\overline{p_1n_1n_2}\equiv\xi\bmod g$. Hence, \eqref{eq:phase_cd} and Lemma~\ref{lem:Theta_identity} give
\begin{align}
	\frac{a_1\overline{m}}{gn_1}-\frac{a_2\overline{m}}{gn_2}
	\equiv
	\frac{\overline{c}\,(a_1\overline{\rho}-a_2\overline{\sigma})+da_2\overline{p_1}\,\overline{\rho}\,\overline{\sigma}}{g}
	+\frac{d(a_2p_2-a_1p_1)\overline{gp_2n_2}}{p_1n_1}
	-\frac{da_2}{gp_1n_1n_2}\bmod 1,
	\label{eq:phase_classes}
\end{align}
that is, the exponential in \eqref{eq:Hg_with_cd} for \(\mathcal H_g\) equals $\e(\xi/g)\,\e(\Theta^{a_1,a_2}_{d,p_1,p_2}(n_1,n_2))$. The first factor \(\e(\xi/g)\) is constant on the classes. The second is a Kloosterman fraction with modulus $p_1n_1$, which is the term that carries the cancellation, and its rational term satisfies
\[
\Big|\frac{da_2}{gp_1n_1n_2}\Big|
<
\frac{8LN}{gM}\,
\frac{2A}{gL(N/g)^2}
=
\frac{16A}{MN},
\]
since $a_2\le2A$, $p_1>L$ and $n_1,n_2>N/g$. Its exponential will then be incorporated into the weight, leaving only the Kloosterman fraction with modulus $p_1n_1$.

\subsection{Removing the coefficients by Cauchy-Schwarz}
\label{sec:cauchy_schwarz}
Grouping the pairs $n_1,n_2$ in \eqref{eq:Hg_with_cd} according to their residue classes modulo $g$, and using
\eqref{eq:phase_classes}, we obtain
\begin{align}
	\Hcal_g
	=
	\sum_{\substack{\rho,\sigma\bmod g\\(\rho\sigma,g)=1}}
	\sum_{\substack{n_1,n_2\sim N/g\\(n_1,n_2)=(g,n_1n_2)=1\\n_1\equiv\rho,\ n_2\equiv\sigma\bmod g}}
	\sum_{a_2\in\Acal}
	\beta_{n_1}\overline{\beta_{n_2}}\,\overline{\nu_{a_2}}\,
	K_{g,\rho,\sigma}(n_1,n_2,a_2),
	\label{eq:Hg_with_K}
\end{align}
where, for reduced residues $\rho,\sigma\bmod g$, for coprime positive integers $n_1,n_2$ with $(g,n_1n_2)=1$ and for $a_2\in\Acal$, the quantity $K_{g,\rho,\sigma}(n_1,n_2,a_2)$ is defined by
\begin{align}
	K_{g,\rho,\sigma}(n_1,n_2,a_2)
	:=&
	\sum_{\substack{c\bmod g\\(c,g)=1}}
	\sum_{0<|d|<D}
	\sum_{\substack{p_1,p_2\in\Pcal,\ p_1\ne p_2\\(p_1p_2,gdn_1n_2)=1\\p_1\rho-p_2\sigma\equiv cd\bmod g\\p_1n_1-p_2n_2\equiv cd\overline{\rho}n_1\bmod g|d|}}
	w\Big(\frac{p_1n_1-p_2n_2}{dM}\Big)
	\notag\\
	&\quad\times
	\sum_{a_1\in\Acal}\nu_{a_1}\,
	\e\Big(\frac{\xi}{g}\Big)\e\big(\Theta^{a_1,a_2}_{d,p_1,p_2}(n_1,n_2)\big),
	\label{eq:def_K}
\end{align}
with $\xi=\xi(c,d,p_1,a_1,a_2;\rho,\sigma)$. On the classes $n_1\equiv\rho$ and $n_2\equiv\sigma\bmod g$,
the two congruences in   \eqref{eq:def_K} are together equivalent to 
\eqref{eq:affine_progression}, namely
\[
p_1n_1-p_2n_2\equiv cd\bmod g|d|.
\]
Thus, on these residue classes,
$K_{g,\rho,\sigma}(n_1,n_2,a_2)$ is exactly the inner sum
over $c,d,p_1,p_2,a_1$ in \eqref{eq:Hg_with_cd}.

We extend the definition \eqref{eq:def_K} of
$K_{g,\rho,\sigma}$ by dropping the restrictions
$n_1\equiv\rho$ and $n_2\equiv\sigma\bmod g$. We shall use the extended $K_{g,\rho,\sigma}$ only inside
a sum of squares, where enlarging the range of summation
can only increase the sum. This allows us to drop the
residue restrictions after applying the Cauchy-Schwarz inequality.
Moreover, \eqref{eq:consequences} implies $(d,n_1n_2)=1$
for every summand in \eqref{eq:def_K}. Together with
$(d,p_1p_2)=1$,  this gives
\[
	(d,p_1p_2n_1n_2)=1.
\]

We now apply the Cauchy-Schwarz inequality to the sum over $\rho,\sigma,n_1,n_2,a_2$ in \eqref{eq:Hg_with_K}.
 Each pair $n_1,n_2$ determines a unique pair of residues
$\rho,\sigma$, so that $|\Hcal_g|^2$ is at most the product of
\[\sum_{n_1,n_2\sim N/g}\sum_{a_2\in\mathcal A} |\beta_{n_1}|^2|\beta_{n_2}|^2|\nu_{a_2}|^2 = \nn{\beta}^4\nn{\nu}^2
\]
and the sum of $|K_{g,\rho,\sigma}(n_1,n_2,a_2)|^2$ over the range of summation of \eqref{eq:Hg_with_K}. To bound the latter sum, choose a fixed real function
\[
	\psi\in C_c^\infty((1/2,3)),\qquad
	0\le\psi\le1,\qquad
	\psi(t)=1\ \text{on }[1,2].
\]
Since the weight $\psi(gn_1/N)\psi(gn_2/N)$ equals $1$ on
the original range $n_1,n_2\sim N/g$, we may insert this weight into the sum of
$|K_{g,\rho,\sigma}(n_1,n_2,a_2)|^2$. Moreover, we may enlarge the
range to all positive $n_1,n_2$ satisfying
$(n_1,n_2)=(g,n_1n_2)=1$, without imposing
$n_1\equiv\rho$ or $n_2\equiv\sigma\bmod g$, since every added term is
nonnegative. We therefore obtain
\begin{align}
	|\Hcal_g|^2
	\le
	\nn{\beta}^4\nn{\nu}^2\,\Qcal_g,
	\label{eq:Hg_CS}
\end{align}
where
\begin{align}
	\Qcal_g
	:=
	\sum_{\substack{\rho,\sigma\bmod g\\(\rho\sigma,g)=1}}
	\sum_{a_2\in\Acal}
	\sum_{\substack{n_1,n_2\ge1\\(n_1,n_2)=(g,n_1n_2)=1}}
	\psi\Big(\frac{gn_1}{N}\Big)\psi\Big(\frac{gn_2}{N}\Big)
	|K_{g,\rho,\sigma}(n_1,n_2,a_2)|^2.
	\label{eq:def_Q}
\end{align}
The sum $\Qcal_g$ no longer involves $\beta$. It remains to bound it in terms of $g$, $M$, $N$, $A$ and $L$.
\subsection{Expanding the square}
\label{sec:expanding_square}
To estimate the sum $\Qcal_g$, we expand the square in \eqref{eq:def_Q} and sum over $n_1,n_2$, with the other variables fixed. We first restrict $d$ to a fixed sign and to a fixed dyadic interval $R\le|d|<2R$. Thus, within each squared sum, the variables $d,d'$ have the same sign and satisfy $R\le|d|,|d'|<2R$. 

For each sign $\pm$ and each $R=2^j$ with
$j\ge0$ an integer and $R<D$, define $K^{\pm,R}_{g,\rho,\sigma}$ by restricting the sum defining $K_{g,\rho,\sigma}$ in \eqref{eq:def_K} to
\[
\pm d>0,
\qquad
R\le |d|<\min(2R,D).
\]
The upper bound $D$ affects only the last dyadic interval. By this decomposition,
\[
	K_{g,\rho,\sigma}
	=
	\sum_{\pm}\sum_{1\le R<D}K^{\pm,R}_{g,\rho,\sigma},
\]
where $R$ runs over powers of two. Since $D=8LN/(gM)\le8L$, the number of pairs $(\pm,R)$ is at most 
\[
2(1+\log_2D)\le8\log_2(2L),
\]
and by the Cauchy-Schwarz inequality,
\begin{align}
	|K_{g,\rho,\sigma}(n_1,n_2,a_2)|^2
	\le
	8\log_2(2L)\sum_{\pm}\sum_{1\le R<D}|K^{\pm,R}_{g,\rho,\sigma}(n_1,n_2,a_2)|^2.
	\label{eq:dyadic_CS}
\end{align}
Define $\Qcal^{\pm,R}_g$ by replacing $K_{g,\rho,\sigma}$  with $K^{\pm,R}_{g,\rho,\sigma}$ in \eqref{eq:def_Q}. Then
\begin{align}
	\Qcal_g
	\le
	8\log_2(2L)\sum_{\pm}\sum_{1\le R<D}\Qcal^{\pm,R}_g,
	\label{eq:Q_dyadic}
\end{align}
and therefore it suffices to bound $\Qcal^{\pm,R}_g$ for a fixed sign and a fixed $R$.

To do so, we first expand the square
\[
|K^{\pm,R}_{g,\rho,\sigma}|^2
=
K^{\pm,R}_{g,\rho,\sigma}\,
\overline{K^{\pm,R}_{g,\rho,\sigma}}.
\]
We write $(c,d,p_1,p_2,a_1)$ for the variables in the first factor and $(c',d',p_1',p_2',a_1')$ for the tuple of the conjugated factor. Each of the two factors is a sum over a tuple $(c,d,p_1,p_2,a_1)$ satisfying the conditions of \eqref{eq:def_K}. We use the abbreviations
\begin{align*}
\xi&=\xi(c,d,p_1,a_1,a_2;\rho,\sigma),
&
\xi'&=\xi(c',d',p_1',a_1',a_2;\rho,\sigma),\\
\Theta&=\Theta^{a_1,a_2}_{d,p_1,p_2}(n_1,n_2),
&
\Theta'&=\Theta^{a_1',a_2}_{d',p_1',p_2'}(n_1,n_2).
\end{align*}
Note that the numerator $a_2$ and the residues $\rho,\sigma$ are the same in both factors. After expanding the square in $\Qcal^{\pm,R}_g$, we  sum over $n_1,n_2$ first, with all other variables fixed. Since $\xi-\xi'$ does not depend on $n_1,n_2$, we obtain
\begin{align}
	\Qcal^{\pm,R}_g
	=
	\sum_{\substack{\rho,\sigma\bmod g\\(\rho\sigma,g)=1}}
	\sum_{a_2\in\Acal}
	\sumtwo_{\substack{c,c'\bmod g\\(cc',g)=1}}
	\sumtwo_{\substack{\pm d,\pm d'>0\\R\le|d|,|d'|<\min(2R,D)}}
	\sum_{\substack{p_1,p_2,p_1',p_2'\in\Pcal\\p_1\ne p_2,\ p_1'\ne p_2'\\(p_1p_2,gd)=(p_1'p_2',gd')=1\\p_1\rho-p_2\sigma\equiv cd\bmod g\\p_1'\rho-p_2'\sigma\equiv c'd'\bmod g}}
	\sumtwo_{a_1,a_1'\in\Acal}
	\nu_{a_1}\overline{\nu_{a_1'}}\,
	\e\Big(\frac{\xi-\xi'}{g}\Big)\,
	\Jcal,
	\label{eq:Q_expanded}
\end{align}
where
\begin{align}
	\Jcal
	:=
	\sum_{\substack{n_1,n_2\ge1,\ (n_1,n_2)=1\\(gp_1p_2p_1'p_2',n_1n_2)=1\\p_1n_1-p_2n_2\equiv cd\overline{\rho}n_1\bmod g|d|\\p_1'n_1-p_2'n_2\equiv c'd'\overline{\rho}n_1\bmod g|d'|}}
	\psi\Big(\frac{gn_1}{N}\Big)\psi\Big(\frac{gn_2}{N}\Big)
	w\Big(\frac{p_1n_1-p_2n_2}{dM}\Big)w\Big(\frac{p_1'n_1-p_2'n_2}{d'M}\Big)
	\e\big(\Theta-\Theta'\big).
	\label{eq:def_J}
\end{align}
The sum $\Jcal$ depends on $\rho$, $a_2$, and the two tuples $(c,d,p_1,p_2,a_1)$ and $(c',d',p_1',p_2',a_1')$, but does not depend on $\sigma$. We suppress these dependences from the notation. By the triangle inequality,
\begin{align}
	\Qcal^{\pm,R}_g
	\le
	\sum_{\substack{\rho,\sigma\bmod g\\(\rho\sigma,g)=1}}
	\sum_{a_2\in\Acal}
	\sumtwo_{\substack{c,c'\bmod g\\(cc',g)=1}}
	\sumtwo_{\substack{\pm d,\pm d'>0\\R\le|d|,|d'|<\min(2R,D)}}
	\sum_{\substack{p_1,p_2,p_1',p_2'\in\Pcal\\\text{as in \eqref{eq:Q_expanded}}}}
	\sumtwo_{a_1,a_1'\in\Acal}
	|\nu_{a_1}\nu_{a_1'}|\,|\Jcal|.
	\label{eq:Q_absolute}
\end{align}

To estimate $\Jcal$, we next combine $\Theta-\Theta'$
into a single Kloosterman fraction and a rational term. Define
\begin{align}
	\Delta
	:=
	d\,p_1'p_2'\,(a_2p_2-a_1p_1)-d'\,p_1p_2\,(a_2p_2'-a_1'p_1'),
	\label{eq:def_Delta}
\end{align}
which depends on the two tuples and on $a_2$, but not on $c,c'$, $\rho$, $\sigma$, $n_1$ or $n_2$. Since
\[
	\frac{\Delta}{p_1p_1'p_2p_2'}
	=
	\frac{d(a_2p_2-a_1p_1)}{p_1p_2}
	-
	\frac{d'(a_2p_2'-a_1'p_1')}{p_1'p_2'},
\]
we have $\Delta=0$ precisely when the two fractions
on the right are equal. 
\begin{lemma}
\label{lem:phase_difference}
Let $g$, $n_1$ and $n_2$ be positive integers with $(n_1,n_2)=(g,n_1n_2)=1$, and $p_1,p_2,p_1',p_2'$ be primes with $p_1\ne p_2$, $p_1'\ne p_2'$ and $(p_1p_2p_1'p_2',gn_1n_2)=1$. Let $d,d',a_1,a_1',a_2$ be integers. Define $\Theta$, $\Theta'$
and $\Delta$ as in \eqref{eq:def_Theta} and
\eqref{eq:def_Delta}, respectively.
\begin{enumerate}
\item[(i)] There is an integer $\Delta_{n_1}$, depending on $n_1,d,d',p_1,p_2,p_1',p_2',a_1,a_1',a_2$, but not on $n_2$, and
satisfying
\[
	(\Delta_{n_1},n_1)\mid\Delta,
\]
such that
\begin{align}
	\Theta-\Theta'
	\equiv
	\frac{\Delta_{n_1}\,\overline{gn_2}}{[p_1,p_1']\,n_1}
	-\frac{a_2(dp_1'-d'p_1)}{gp_1p_1'n_1n_2}\bmod 1.
	\label{eq:phase_difference_general}
\end{align}
\item[(ii)] If $(p_1p_1',p_2p_2')=1$, then
\begin{align}
	\Theta-\Theta'
	\equiv
	\frac{\Delta\,\overline{gp_2p_2'n_2}}{p_1p_1'n_1}
	-\frac{a_2(dp_1'-d'p_1)}{gp_1p_1'n_1n_2}\bmod 1.
	\label{eq:phase_difference_generic}
\end{align}
\end{enumerate}
\end{lemma}
\begin{proof}

We first establish (ii). Suppose that $(p_1p_1',p_2p_2')=1$. Then both $gp_2n_2$ and $gp_2'n_2$ are coprime to $p_1p_1'n_1$. We may therefore take their inverses modulo $p_1p_1'n_1$ when computing the Kloosterman fractions of $\Theta$ and $\Theta'$. We then obtain
\begin{align}
	\Theta-\Theta'
	\equiv
	\frac{dp_1'(a_2p_2-a_1p_1)\overline{gp_2n_2}-d'p_1(a_2p_2'-a_1'p_1')\overline{gp_2'n_2}}{p_1p_1'n_1}
	-\frac{da_2}{gp_1n_1n_2}+\frac{d'a_2}{gp_1'n_1n_2}\bmod1.\label{eq:theta-theta'}
\end{align}
Modulo $p_1p_1'n_1$, we have
\[
\overline{gp_2n_2}
\equiv p_2'\overline{gp_2p_2'n_2},
\qquad
\overline{gp_2'n_2}
\equiv p_2\overline{gp_2p_2'n_2}.
\]
Substituting these congruences, by \eqref{eq:def_Delta}, the numerator of the first fraction in \eqref{eq:theta-theta'} is 
\[
\bigl(
dp_1'p_2'(a_2p_2-a_1p_1)
-d'p_1p_2(a_2p_2'-a_1'p_1')
\bigr)\overline{gp_2p_2'n_2}=\Delta\,\overline{gp_2p_2'n_2}.
\]
The last two terms in \eqref{eq:theta-theta'} combine to
\[
-\frac{da_2}{gp_1n_1n_2}
+\frac{d'a_2}{gp_1'n_1n_2}
=
-\frac{a_2(dp_1'-d'p_1)}{gp_1p_1'n_1n_2}.
\]
Hence
\[
\Theta-\Theta'
\equiv
\frac{\Delta\,\overline{gp_2p_2'n_2}}{p_1p_1'n_1}
-\frac{a_2(dp_1'-d'p_1)}{gp_1p_1'n_1n_2}
\bmod 1,
\]
which proves (ii).

We now prove (i). Applying Lemma~\ref{lem:Theta_identity} to $\Theta$ and $\Theta'$ and subtracting give
\[
	\Theta-\Theta'
	\equiv
	-\frac{(da_1\overline{p_2}-d'a_1'\overline{p_2'})\,\overline{gn_2}}{n_1}
	-\frac{a_2(d\overline{p_1}-d'\overline{p_1'})\,\overline{gn_1}}{n_2}
	-\frac{a_2(d\overline{p_1}-d'\overline{p_1'})\,\overline{n_1n_2}}{g}\bmod 1.
\]
Define the integers
\[
	u:=da_1\frac{[p_2,p_2']}{p_2}-d'a_1'\frac{[p_2,p_2']}{p_2'},
	\qquad
	v:=a_2\Big(d\frac{[p_1,p_1']}{p_1}-d'\frac{[p_1,p_1']}{p_1'}\Big).
\]
For any modulus coprime to $p_2p_2'$, we have 
\[
\overline{p_2}
\equiv
\frac{[p_2,p_2']}{p_2}\,\overline{[p_2,p_2']},
\qquad
\overline{p_2'}
\equiv
\frac{[p_2,p_2']}{p_2'}\,\overline{[p_2,p_2']},
\]
and similarly for $p_1,p_1'$. Thus,
\begin{align}
	\Theta-\Theta'
	\equiv
	-\frac{u\,\overline{g[p_2,p_2']n_2}}{n_1}
	-\frac{v\,\overline{g[p_1,p_1']n_1}}{n_2}
	-\frac{v\,\overline{[p_1,p_1']n_1n_2}}{g}\bmod 1.\label{eq:theta-theta' second}
\end{align}
By Lemma~\ref{lem:crt_inverse} with $x=[p_1,p_1']n_1$, $y=g$ and $z=n_2$, the last two terms in \eqref{eq:theta-theta' second} combine to 
\[
-\frac{v\overline{g[p_1,p_1']n_1}}{n_2}
-
\frac{v\overline{[p_1,p_1']n_1n_2}}{g}
\equiv
-\frac{v\overline{[p_1,p_1']n_1}}{gn_2}
\bmod 1.
\]
Then, applying Lemma~\ref{lem:reciprocity} to the coprime integers $[p_1,p_1']n_1$ and $gn_2$ gives
\begin{align}
	-\frac{v\,\overline{[p_1,p_1']n_1}}{gn_2}
	\equiv
	\frac{v\,\overline{gn_2}}{[p_1,p_1']n_1}-\frac{v}{g[p_1,p_1']n_1n_2}\bmod 1.\label{eq:lemma 4.2 lemma 2.1 applied}
\end{align}
 Since $v/[p_1,p_1']=a_2(d/p_1-d'/p_1')$, the last term in \eqref{eq:lemma 4.2 lemma 2.1 applied} is 
\[
-\frac{v}{g[p_1,p_1']n_1n_2}
=
-\frac{a_2(dp_1'-d'p_1)}{gp_1p_1'n_1n_2}.
\]
 It remains to combine the fractions with denominators $n_1$ and $[p_1,p_1']n_1$. Choose $y$ to be an inverse of $[p_2,p_2']$ modulo $n_1$, and take $\overline{gn_2}$ to be the inverse of $gn_2$ modulo $[p_1,p_1']n_1$, which is also an inverse modulo $n_1$. Then, $\overline{g[p_2,p_2']n_2}
\equiv y\,\overline{gn_2}\bmod{n_1}$.
Consequently,
\[
	-\frac{u\,\overline{g[p_2,p_2']n_2}}{n_1}
	\equiv
	-\frac{[p_1,p_1']\,u\,y\,\overline{gn_2}}{[p_1,p_1']n_1}\bmod 1.
\]
Define
\[
	\Delta_{n_1}:=v-[p_1,p_1']\,u\,y,
\]
which depends on $n_1$ through $y$ and does not depend on $n_2$. We then obtain
\[
\Theta-\Theta'
\equiv
\frac{\Delta_{n_1}\,\overline{gn_2}}{[p_1,p_1']n_1}
-\frac{a_2(dp_1'-d'p_1)}{gp_1p_1'n_1n_2}
\bmod 1.
\]
Finally, since $([p_2,p_2'],n_1)=1$ and $[p_2,p_2']y\equiv1\bmod n_1$,
\[
	(\Delta_{n_1},n_1)=([p_2,p_2']\Delta_{n_1},n_1)=([p_2,p_2']v-[p_1,p_1']u,\,n_1),
\]
and the definitions of $u$, $v$, and $\Delta$ give
\[
	[p_2,p_2']v-[p_1,p_1']u
	=
	[p_1,p_1'][p_2,p_2']\bigg(a_2\Big(\frac{d}{p_1}-\frac{d'}{p_1'}\Big)-\Big(\frac{da_1}{p_2}-\frac{d'a_1'}{p_2'}\Big)\bigg)
	=
	\frac{[p_1,p_1'][p_2,p_2']}{p_1p_1'p_2p_2'}\,\Delta.
\]
Equivalently, 
\[
\Delta
=
\frac{p_1p_1'}{[p_1,p_1']}
\frac{p_2p_2'}{[p_2,p_2']}
\bigl([p_2,p_2']v-[p_1,p_1']u\bigr).
\]
Therefore, the integer $[p_2,p_2']v-[p_1,p_1']u$ divides $\Delta$, and so does $(\Delta_{n_1},n_1)$. This completes the proof.
\end{proof}

Recall that our goal is to estimate \eqref{eq:Q_absolute}. We divide the terms of \eqref{eq:Q_absolute} into three
groups. The first group consists of the terms with $\Delta=0$, and they are bounded by counting in Section~\ref{sec:delta_zero}.
 The second group consists of the terms with $\Delta\ne0$ in which $p_1,p_2,p_1',p_2'$ are four distinct primes and $(dd',p_1p_2p_1'p_2')=1$. These are estimated in
Section~\ref{sec:generic} by Lemma~\ref{lem:progression},
using \eqref{eq:phase_difference_generic}.  The third group consists of the remaining terms with $\Delta\ne0$. These are estimated in
Section~\ref{sec:nongeneric} by Lemma~\ref{lem:completion},
using \eqref{eq:phase_difference_general}. Denote the contributions of these three groups to the
right-hand side of \eqref{eq:Q_absolute} by
$\Qcal_g^{\pm,R}(0)$, $\Qcal_g^{\pm,R}(1)$ and
$\Qcal_g^{\pm,R}(2)$, respectively. Then
\begin{align}
	\Qcal_g^{\pm,R}
	\le
	\Qcal_g^{\pm,R}(0)
	+\Qcal_g^{\pm,R}(1)
	+\Qcal_g^{\pm,R}(2).
	\label{eq:Q_three_parts}
\end{align}

Before treating these groups, we count the quadruples
$(\rho,\sigma,c,c')$ for which the summation conditions
defining $\Jcal$ admit a solution, with
$d,d',p_1,p_2,p_1',p_2'$ fixed.
The bounds for $|\Jcal|$ obtained in
Sections~\ref{sec:generic} and \ref{sec:nongeneric}
are independent of these residues, so the counts allow
us to sum those bounds over the residue classes in
\eqref{eq:Q_absolute}.

\subsection{Counting the residue classes}
\label{sec:counting}
 There are trivially at most $\varphi(g)^4$ quadruples
$(\rho,\sigma,c,c')$ of reduced residue classes modulo $g$.
The following lemma improves this bound to $g\varphi(g)$
when the congruences admit a common solution, and to
$\varphi(g)$ when $n_1,n_2$ are fixed.
In either case, we do not require $n_1\equiv\rho$
or $n_2\equiv\sigma\bmod g$.
\begin{lemma}
\label{lem:counting}
Let $g\in\N$ and $d,d'$ be nonzero integers. Suppose that $p_1,p_2,p_1',p_2'$ are primes satisfying
\[
	p_1\ne p_2,\qquad p_1'\ne p_2',\qquad
	(p_1p_2p_1'p_2',g)=1.
\]
Consider quadruples $(\rho,\sigma,c,c')$ of reduced
residue classes modulo $g$ satisfying \eqref{eq:residue_condition}, that is,
\[
	p_1\rho-p_2\sigma\equiv cd\bmod g.
\]
\begin{enumerate}
\item[(i)]The number of such quadruples $(\rho,\sigma,c,c')$ for which the congruences
\begin{align}
	p_1n_1-p_2n_2\equiv cd\,\overline{\rho}\,n_1\bmod g|d|
	\quadand
	p_1'n_1-p_2'n_2\equiv c'd'\,\overline{\rho}\,n_1\bmod g|d'|
	\label{eq:two_progressions}
\end{align}
in \eqref{eq:def_J} have a common solution $n_1,n_2$ with $(n_1,gdd')=1$ is at most $g\varphi(g)$.
\item[(ii)] For fixed integers $n_1,n_2$ with $(n_1,g)=1$, the number of such quadruples $(\rho,\sigma,c,c')$ that also satisfy \eqref{eq:two_progressions} is at most $\varphi(g)$.
\end{enumerate}
\end{lemma}
\begin{proof}
We first prove (i). Write $t=(d,d')>0$, and $d=tu$, $d'=tu'$, so that $(u,u')=1$. Consider a quadruple $(\rho,\sigma,c,c')$ counted in (i), and let $n_1,n_2$ be a common solution of \eqref{eq:two_progressions} with $(n_1,gdd')=1$. 

Both congruences in \eqref{eq:two_progressions} hold modulo $gt$. Multiplying the first by $p_2'$, the second by $p_2$, and subtracting, we obtain
\[
	(p_1p_2'-p_1'p_2)\,n_1\equiv t\,(cup_2'-c'u'p_2)\,\overline{\rho}\,n_1\bmod gt.
\]
Since $(n_1,gt)=1$, we may cancel $n_1$. It follows that $t\mid p_1p_2'-p_1'p_2$.
Dividing the congruence by $t$ and multiplying by $\rho$, we get
\begin{align}
	cup_2'-c'u'p_2\equiv\frac{p_1p_2'-p_1'p_2}{t}\,\rho\bmod g.
	\label{eq:compatibility}
\end{align}
Thus, every quadruple counted in (i) satisfies \eqref{eq:compatibility}, and moreover $\sigma\equiv\overline{p_2}(p_1\rho-cd)\bmod g$ is determined by $\rho$ and $c$. Fix one of the $\varphi(g)$ possible values of $\rho$. The map 
\[
(\Z/g\Z)^2\longrightarrow\Z/g\Z,
\qquad
(c,c')\longmapsto cup_2'-c'u'p_2
\]
 is a surjective homomorphism. Indeed, since $(p_2p_2',g)=1$, any prime dividing
$g$, $up_2'$ and $u'p_2$ would divide both $u$ and $u'$,
contradicting $(u,u')=1$. Hence, its kernel has $g^2/g =g$ elements, and each value is attained by exactly $g$ pairs $(c,c')$. Restricting $c,c'$ to reduced residue classes can only decrease this count. Thus, since each $\rho$ gives at most $g$ quadruples satisfying
\eqref{eq:residue_condition} and \eqref{eq:compatibility}, the number of total such quadruples is at most $g\varphi(g)$.

We now prove (ii). Fix one of the possible $\varphi(g)$ values of $\rho$. By \eqref{eq:consequences}, $n_2\equiv\sigma\overline{\rho}n_1\bmod g$. Since $(n_1,g)=1$, 
\[
\sigma\equiv\rho\,n_2\overline{n_1}\bmod g
\]
is determined. The first congruence in \eqref{eq:two_progressions}
also determines $c$ uniquely modulo $g$. Indeed, suppose $c$ and $c_2$ both satisfy the first congruence. Subtracting gives
\[
g|d|\mid (c-c_2)d\,\overline{\rho}n_1.
\]
Since $d=\pm|d|\ne0$, dividing by $|d|$ gives
\[
g\mid (c-c_2)\,\overline{\rho}n_1.
\]
Since $(\overline{\rho}n_1,g)=1$ we get $c\equiv c_2\bmod g$. In the same way $c'$ is determined by the second congruence in \eqref{eq:two_progressions}. Thus, at most $\varphi(g)$ quadruples satisfy the conditions of (ii).
\end{proof}

\subsection{The terms with $\Delta=0$}
\label{sec:delta_zero}
In this subsection, we bound $\Qcal_g^{\pm,R}(0)$, the contribution from $\Delta=0$. This includes, but is not restricted to, the diagonal terms in which the two tuples coincide. All these terms are bounded trivially, without using the cancellation in $\e((\xi-\xi')/g)$ or in the exponential of the rational term in \eqref{eq:phase_difference_general}.

 By the remark after \eqref{eq:def_Delta}, $\Delta=0$ is equivalent to
\begin{align}
	\frac{d(a_2p_2-a_1p_1)}{p_1p_2}=\frac{d'(a_2p_2'-a_1'p_1')}{p_1'p_2'}.
	\label{eq:Delta_zero_fractions}
\end{align}
Both sides are fractions in lowest terms. Indeed, $(d,p_1p_2)=1$ by the conditions of \eqref{eq:def_K} and
\[
a_2p_2-a_1p_1\equiv a_2p_2\bmod{p_1},
\qquad
a_2p_2-a_1p_1\equiv-a_1p_1\bmod{p_2}.
\]
Both residues are nonzero, because $p_2\ne p_1$ and $0<a_1,a_2\le2A<L<p_1,p_2$ by \eqref{eq:blockrange}. The same argument applies to the right-hand side of \eqref{eq:Delta_zero_fractions}. Therefore, their numerators are both nonzero, and the equality of the reduced fractions forces $p_1p_2=p_1'p_2'$, that is,
\begin{align}
	\{p_1',p_2'\}=\{p_1,p_2\}.
	\label{eq:same_primes}
\end{align}
For given $p_1,p_2$, there are consequently at most two choices of the ordered pair $(p_1',p_2')$. With \eqref{eq:same_primes}, the condition $\Delta=0$ becomes the linear relation
\begin{align}
	dp_1\,a_1-d'p_1'\,a_1'=a_2\,(dp_2-d'p_2')
	\label{eq:a_relation}
\end{align}
between $a_1$ and $a_1'$, where both coefficients $dp_1$ and $d'p_1'$ are nonzero. Thus, when $d,d',p_1,p_2,p_1',p_2'$ and $a_2$ are fixed, each $a_1$ determines at most one $a_1'$, and each $a_1'$ determines at most one $a_1$. The Cauchy-Schwarz inequality therefore gives
\begin{align}
	\sumtwo_{\substack{a_1,a_1'\in\Acal\\\Delta=0}}|\nu_{a_1}\nu_{a_1'}|&\le
\Big(
\sum_{\substack{a_1,a_1'\in\mathcal A\\\Delta=0}}
|\nu_{a_1}|^2
\Big)^{1/2}
\Big(
\sum_{\substack{a_1,a_1'\in\mathcal A\\\Delta=0}}
|\nu_{a_1'}|^2
\Big)^{1/2}\le\nn{\nu}^2.
	\label{eq:a_pairing}
\end{align}
Note that \eqref{eq:a_relation} involves neither $n_1,n_2$ nor the residues $\rho,\sigma,c,c'$.

We now bound $\Jcal$. Since $0\le\psi, w\le1$, $|\Jcal|$ is at most the number of pairs $n_1,n_2$ in the summation of \eqref{eq:def_J} with $\psi(gn_1/N)\psi(gn_2/N)\ne0$. Every such pair satisfies
\[
1\le n_1,n_2<\frac{3N}{g}.
\]
Insert this bound into $\mathcal Q_g^{\pm,R}(0)$, and sum over $n_1,n_2$ before $d,d'$ and the residues. By  \eqref{eq:consequences}, 
\[
d\mid p_1n_1-p_2n_2,\qquad d'\mid p_1'n_1-p_2'n_2.
\]
Among the conditions on $d,d'$ we keep only these two divisibilities and drop the others, which can only increase the sum since all terms are nonnegative. This gives
\begin{align}
	\Qcal^{\pm,R}_g(0)
	\le
	\sum_{a_2\in\Acal}\
	\sum_{\substack{p_1,p_2,p_1',p_2'\in\Pcal\\p_1\ne p_2,\ \{p_1',p_2'\}=\{p_1,p_2\}\\
(p_1p_2,g)=1}}\
	\sumtwo_{\substack{1\le n_1,n_2<3N/g,\\ (n_1,n_2)=1\\(gp_1p_2,n_1n_2)=1}}\
	\sumtwo_{\substack{\pm d,\pm d'>0\\d\mid p_1n_1-p_2n_2\\d'\mid p_1'n_1-p_2'n_2}}
	\Big(\sum_{(\rho,\sigma,c,c')}1\Big)
	\Big(\sumtwo_{\substack{a_1,a_1'\in\Acal\\\Delta=0}}|\nu_{a_1}\nu_{a_1'}|\Big),
	\label{eq:Q0_rearranged}
\end{align}
where the sum over $(\rho,\sigma,c,c')$ runs over the quadruples of reduced residues modulo $g$ satisfying the conditions of \eqref{eq:Q_expanded} and the two congruences in \eqref{eq:two_progressions} for the fixed pair $n_1,n_2$. The two factors in parentheses can be separated because \eqref{eq:a_relation} does not involve the residues and the congruences on the residues do not involve $a_1,a_1'$.

For fixed primes and $n_1,n_2$, the integers $p_1n_1-p_2n_2$ and $p_1'n_1-p_2'n_2$ are nonzero, since the two products in each difference are coprime and both exceed $1$. Moreover,
\[
|p_1n_1-p_2n_2|,
\ |p_1'n_1-p_2'n_2|<\frac{6LN}{g}\le6M^{5/4},
\]
 using $N\le M$ and $L\le M^{1/4}$. Since the signs of $d$ and $d'$ are fixed, by the divisor bound, the number of pairs $d,d'$ is at most
\[
 \tau(|p_1n_1-p_2n_2|)\,\tau(|p_1'n_1-p_2'n_2|)\ll_\eps M^\eps.
\]
For fixed primes, $n_1,n_2$ and $d,d'$, the number of quadruples $(\rho,\sigma,c,c')$ satisfying the congruences in \eqref{eq:def_J} is at most $\varphi(g)\le g$ by Lemma~\ref{lem:counting}(ii). The hypotheses of the lemma are satisfied, since $(p_1p_2p_1'p_2',g)=1$, $(n_1,g)=1$, and $n_1,n_2$ in the summation of \eqref{eq:def_J} satisfy \eqref{eq:two_progressions}.

 The number of $a_2$ is at most $2A$, the number of quadruples of primes is at most $2|\Pcal|^2\le8L^2$, and the number of pairs $n_1,n_2$ is at most $9N^2/g^2$. Finally, by \eqref{eq:a_pairing}, the last factor in \eqref{eq:Q0_rearranged} is at most $\nn{\nu}^2$. Altogether, we have
\begin{align}
	\Qcal^{\pm,R}_g(0)
	\ll_\eps
(2A)(8L^2)\frac{9N^2}{g^2}
\,M^\eps\,g\,\nn{\nu}^2\ll_\eps
	\nn{\nu}^2\,\frac{AL^2N^2}{g}\,M^{\eps}.
	\label{eq:Q0_bound}
\end{align}
Note that this bound does not depend on the sign or on $R$.

\subsection{Four distinct primes}
\label{sec:generic}
In this subsection, we bound $\Qcal^{\pm,R}_g(1)$, the contribution from the terms in \eqref{eq:Q_absolute} in which $p_1,p_2,p_1',p_2'$ are four distinct primes and $(dd',p_1p_2p_1'p_2')=1$. Fix the sign $\pm$, the number $R$, the residues $\rho,\sigma$, the numerator $a_2$, and the two tuples $(c,d,p_1,p_2,a_1)$ and $(c',d',p_1',p_2',a_1')$, with $\pm d,\pm d'>0$ and $R\le|d|,|d'|<\min(2R,D)$. Set
\begin{align}
	\eta_R:=\min\Big(1,\frac{gMR}{LN}\Big).
	\label{eq:def_eta}
\end{align}
The support conditions of the factors $w$ and $\psi$
in \eqref{eq:def_J} restrict $n_2/n_1$ to an interval
of length $O(\eta_R)$, as established in
Lemma~\ref{lem:weight}. We will apply Lemma~\ref{lem:progression} with
$\delta=\eta_R$ to estimate $\Jcal$, and then sum the resulting
bound over the residues, numerators, primes and divisors. We begin by bounding $\Delta$ and showing that
$\Delta$ is coprime to the four amplifier primes.
\begin{lemma}
\label{lem:Delta_generic}
Suppose that the primes $p_1,p_2,p_1',p_2'$ are distinct and that $(dd',p_1p_2p_1'p_2')=1$. Then 
\[
(\Delta,p_1p_2p_1'p_2')=1,\qquad|\Delta|<128ARL^3.
\]
In particular $\Delta\ne0$.
\end{lemma}
\begin{proof}
By \eqref{eq:def_Delta}, we have
\begin{align*}
	\Delta&\equiv dp_1'p_2'a_2p_2 \ \bmod{p_1},
	&\qquad
	\Delta&\equiv-d'p_1p_2a_2p_2' \ \bmod{p_1'},\\
	\Delta&\equiv-dp_1'p_2'a_1p_1 \ \bmod{p_2},
	&\qquad
	\ \Delta&\equiv d'p_1p_2a_1'p_1' \bmod{p_2'}.
\end{align*}
Each factor on each right side is a unit modulo the corresponding modulus. This follows from the distinctness of the primes, the condition $(dd',p_1p_2p_1'p_2')=1$, and $0<a_2\le2A<L<p_1$ by \eqref{eq:blockrange}. Hence, none of the four primes divides $\Delta$.

To bound the size of $\Delta$, we have $|d|,|d'|<2R$, each prime is at most $2L$, and each numerator is at most $2A$. Therefore, each of the two terms in \eqref{eq:def_Delta} has absolute value less than 
\[
2R\cdot4L^2\cdot(2A\cdot2L+2A\cdot2L)=64ARL^3.
\]
\end{proof}

We next combine the two congruences in \eqref{eq:def_J} into a single one. Note that the following lemma does not require the four primes to be distinct, and will be used again in Section~\ref{sec:nongeneric}.
\begin{lemma}
\label{lem:merged_progression}
Let $d,d'$ be nonzero integers and $p_1\ne p_2$, $p_1'\ne p_2'$ primes with $(p_1p_2,gd)=(p_1'p_2',gd')=1$. Suppose that $\rho,\sigma,c,c'$ are reduced residues modulo $g$ with 
\begin{align}
p_1\rho-p_2\sigma\equiv cd\ \bmod g,
\qquad
p_1'\rho-p_2'\sigma\equiv c'd'\ \bmod g.\label{eq:two_congruences_lem4.5}
\end{align}
Set $k:=g[d,d']$. Then either the two congruences in \eqref{eq:two_progressions}
\begin{align*}
p_1n_1-p_2n_2
&\equiv cd\,\overline\rho n_1\ \bmod{g|d|},\\
p_1'n_1-p_2'n_2
&\equiv c'd'\,\overline\rho n_1\ \bmod{g|d'|}
\end{align*}
have no common solution $n_1,n_2$ with $(n_1,k)=1$, or there is an integer $\kappa$ with $(\kappa,k)=1$ such that, for all integers $n_1,n_2$ with $(n_1,k)=1$, the two congruences in \eqref{eq:two_progressions} hold if and only if
\begin{align}
	n_2\equiv\kappa n_1\bmod k.
	\label{eq:merged_progression}
\end{align}
\end{lemma}
\begin{proof}
Since $(p_2,g|d|)=1$, there exists an integer $\kappa_1$ with
\begin{align}
p_2\kappa_1\equiv p_1-cd\,\overline{\rho}\bmod g|d|,\label{eq:kappa1 defn}
\end{align}
and the first congruence in \eqref{eq:two_progressions} is then equivalent to
\[
n_2\equiv\kappa_1n_1\bmod g|d|.
\]
The integer $\kappa_1$ is a unit modulo $g|d|$. Indeed, reducing \eqref{eq:kappa1 defn} modulo $|d|$, we have 
\[
p_2\kappa_1\equiv p_1\bmod |d|.
\]
Since $(p_1,d)=1$, this implies $(\kappa_1,d)=1$. Modulo $g$, the first congruence in \eqref{eq:two_congruences_lem4.5} implies 
\[
p_2\kappa_1\equiv p_1-(p_1\rho-p_2\sigma)\overline{\rho}\equiv p_2\sigma\overline{\rho}\bmod g.
\]
Canceling $p_2$, we obtain $\kappa_1\equiv\sigma\overline{\rho}\bmod g$. Since $\rho,\sigma$ are reduced residues, $(\kappa_1,g)=1$. Together with $(\kappa_1,d)=1$, this proves $(\kappa_1,g|d|)=1$.

Similarly, there exists an integer $\kappa_2$ such that
the second congruence in \eqref{eq:two_progressions} is then equivalent to
\[
n_2\equiv\kappa_2n_1\bmod g|d'|,
\]
and the integer $\kappa_2$ is a unit modulo $g|d'|$. 

We now determine when the two congruences in \eqref{eq:two_progressions} have common solutions. Their moduli satisfy
\[
(g|d|,g|d'|)=g(d,d'),\qquad
[g|d|,g|d'|]=k.
\]
Let $n_1$ be coprime to $k$. A common solution must satisfy
\[
\kappa_1\equiv\kappa_2\bmod g(d,d'),
\]
since otherwise $\kappa_1n_1\not\equiv\kappa_2n_1\bmod g(d,d')$. Suppose now that this condition holds. There exists a unique integer $\kappa$ modulo $k$, with 
\begin{align}
\kappa\equiv\kappa_1\bmod g|d|,\qquad \kappa\equiv\kappa_2\bmod g|d'|.\label{eq:two_equi_relation_on_kappa}
\end{align}
Indeed, write $\kappa=\kappa_1+g|d|\,x$ for some integer $x$. Then the second congruence in \eqref{eq:two_equi_relation_on_kappa} gives 
\[
g|d|\,x\equiv\kappa_2-\kappa_1\bmod g|d'|,
\]
 which is solvable because $(g|d|,g|d'|)=g(d,d')$ divides $\kappa_2-\kappa_1$. Any two choices of $\kappa$ differ by a multiple of both moduli, hence by a multiple of their least common multiple $k$.
 
 With this choice of $\kappa$, the two congruences $n_2\equiv\kappa_1n_1\bmod g|d|$ and $n_2\equiv\kappa_2n_1\bmod g|d'|$ together are equivalent to \eqref{eq:merged_progression}. Finally, $\kappa$ is a unit modulo $g|d|$ and modulo $g|d'|$, because it agrees with $\kappa_1$ and $\kappa_2$, respectively. Every prime dividing $k$ divides at least one of these moduli, so $(\kappa,k)=1$.
\end{proof}
We apply this lemma to $\Jcal$. By \eqref{eq:consequences}, every pair in the sum defining $\Jcal$ in \eqref{eq:def_J} satisfies $(n_1,dd')=1$ and $(n_1,g)=1$. Thus, $(n_1,k)=1$, where $k=g[d,d']$. If the two congruences in \eqref{eq:def_J}, which are the congruences in \eqref{eq:two_progressions}, have no common solution with $(n_1,k)=1$, then $\Jcal=0$. Otherwise, Lemma~\ref{lem:merged_progression} allows us to replace them by
\[
n_2\equiv\kappa n_1\bmod k,\qquad(\kappa,k)=1,
\]
where $\kappa$ is as defined in Lemma~\ref{lem:merged_progression}. This also implies $(n_2,k)=1$, since $\kappa$ and $n_1$ are units modulo $k$. Hence,
\begin{align}
	\Jcal
	=&
	\sum_{\substack{n_1\ge1\\(n_1,\,gp_1p_1'p_2p_2'[d,d'])=1}}\ 
	\sum_{\substack{n_2\ge1,\ n_2\equiv\kappa n_1\bmod k\\(n_2,\,gp_1p_1'p_2p_2'n_1)=1}}
	\psi\Big(\frac{gn_1}{N}\Big)\psi\Big(\frac{gn_2}{N}\Big)
	\notag\\
	&\quad\times
	w\Big(\frac{p_1n_1-p_2n_2}{dM}\Big)w\Big(\frac{p_1'n_1-p_2'n_2}{d'M}\Big)
	\e\big(\Theta-\Theta'\big).
	\label{eq:J_merged}
\end{align}

In the present case, since the four primes are distinct, $(p_1p_1',p_2p_2')=1$. Lemma~\ref{lem:phase_difference}(ii) then gives
\begin{align}
	\e\big(\Theta-\Theta'\big)
	=
	\e\Big(\frac{\Delta\,\overline{gp_2p_2'n_2}}{p_1p_1'n_1}\Big)
	\e\Big(-\frac{a_2(dp_1'-d'p_1)}{gp_1p_1'n_1n_2}\Big).
	\label{eq:phase_generic_product}
\end{align}
The second exponential is small and smooth, and we incorporate it into the weight. For $t,z\in\R$, define
\begin{align}
	W(t,z)
	:=
	\psi(t)\psi(z)\,
	w\Big(\frac{N(p_1t-p_2z)}{gMd}\Big)
	w\Big(\frac{N(p_1't-p_2'z)}{gMd'}\Big)
	\e\Big(-\frac{a_2(dp_1'-d'p_1)\,g}{p_1p_1'N^2\,tz}\Big),
	\label{eq:def_W}
\end{align}
with $W(t,z)=0$ whenever $\psi(t)\psi(z)=0$. Thus, the last exponential in \eqref{eq:def_W} is only evaluated for $t,z\in(1/2,3)$, and $W$ is smooth on $\R^2$ with support in $[1/2,3]^2$. For positive integers $n_1,n_2$,
\begin{align}
	W\Big(\frac{gn_1}{N},\frac{gn_2}{N}\Big)
	=
	\psi\Big(\frac{gn_1}{N}\Big)\psi\Big(\frac{gn_2}{N}\Big)
	w\Big(\frac{p_1n_1-p_2n_2}{dM}\Big)w\Big(\frac{p_1'n_1-p_2'n_2}{d'M}\Big)
	\e\Big(-\frac{a_2(dp_1'-d'p_1)}{gp_1p_1'n_1n_2}\Big).
	\label{eq:W_at_integers}
\end{align}
Consequently, by \eqref{eq:J_merged}, \eqref{eq:phase_generic_product} and \eqref{eq:W_at_integers}, $\Jcal$ is the sum \eqref{eq:def_J_lemma} in Lemma~\ref{lem:progression} with
\begin{align}
	P=p_1p_1',\qquad Q=gp_2p_2',\qquad k=g[d,d'],\qquad Z=\frac{N}{g},
	\label{eq:dictionary}
\end{align}
with the present $\Delta$, $\kappa$, and  the weight $W$ defined in \eqref{eq:def_W}. We now verify the hypotheses of Lemma~\ref{lem:progression}, beginning with the support and derivative bounds for the weight $W$.
\begin{lemma}
\label{lem:weight}
Let $W$ be given by \eqref{eq:def_W} and set $F(t,y):=W(t,ty)$. There is an interval $I$ of length at most $24\eta_R$ such that $F(t,y)=0$ unless $y\in I$. Moreover, for all $i,j\ge0$ and all $t,y\in\R$,
\begin{align}
	|\partial_t^i\partial_y^jF(t,y)|\le C_{ij}\,\eta_R^{-j},
	\label{eq:weight_derivatives}
\end{align}
where the constants $C_{ij}$ depend only on $i$, $j$, $w$ and $\psi$.
\end{lemma}
\begin{proof}
Set
\[
	b(y):=\frac{N(p_1-p_2y)}{gMd},
	\qquad
	\tilde b(y):=\frac{N(p_1'-p_2'y)}{gMd'},
	\qquad
	\lambda:=\frac{a_2(dp_1'-d'p_1)\,g}{p_1p_1'N^2}.
\]
Then
\begin{align}
	F(t,y)=\psi(t)\,\psi(ty)\,w\big(t\,b(y)\big)\,w\big(t\,\tilde b(y)\big)\,\e\Big(-\frac{\lambda}{t^2y}\Big).\label{eq:new_def_F}
\end{align}
If $F(t,y)\ne0$, then $t\in(1/2,3)$ and $ty\in(1/2,3)$, so that $y\in(1/6,6)$. Moreover, $w(tb(y))\ne0$ implies that $|tb(y)|<3$, and hence $|b(y)|<6$. It follows that
\[
	\Big|y-\frac{p_1}{p_2}\Big|<\frac{6gM|d|}{Np_2}<\frac{12gMR}{LN},
\]
where the last inequality follows from $|d|<2R$ and $p_2>L$. If $\eta_R<1$, we take $I=(p_1/p_2-12\eta_R,\,p_1/p_2+12\eta_R)$, which has length $24\eta_R$; if $\eta_R=1$, we take $I=(1/6,6)$, which has length less than $6\le24\eta_R$. In either case, $|I|\le24\eta_R$. This proves the assertion on the support.

For the derivatives, it suffices to consider $S=[1/2,3]\times[1/6,6]$. Indeed, since the support of $\psi$ is a compact subset of $(1/2,3)$, $\psi(t)\psi(ty)$ vanishes in a neighborhood of every point outside $S$. Therefore, all derivatives of $F$ vanish outside $S$, and it suffices to bound them on $S$. On $S$, the derivatives of $\psi(t)\psi(ty)$ of all orders are bounded by constants depending only on $\psi$ and the orders. For the factor $w(tb(y))$,  note that $b$ is an affine function of $y$ with constant derivative $\beta_1:=-Np_2/(gMd)$. By Leibniz's rule, for $i,j\ge0$,
\[
	\partial_t^i\partial_y^j\,w\big(tb(y)\big)
	=
	\beta_1^{\,j}\sum_{v=0}^{\min(i,j)}\binom{i}{v}\frac{j!}{(j-v)!}\,t^{\,j-v}\,b(y)^{\,i-v}\,w^{(i+j-v)}\big(tb(y)\big).
\]
For any nonzero term in this sum, we have $|tb(y)|\le3$. On $S$, this gives $|b(y)|\le6$ for $t\ge1/2$. Together with $t\le3$, this bounds the sum by a constant depending only on $i$, $j$ and $w$. Moreover,
\[
	|\beta_1|=\frac{Np_2}{gM|d|}\le\frac{2LN}{gMR}\le\frac{2}{\eta_R},
\]
where the last inequality is an equality if $\eta_R<1$. Hence,
\[
|\partial_t^i\partial_y^jw(tb(y))|\ll_{i,j}\eta_R^{-j}
\]
on $S$, and the same holds for $w(t\tilde b(y))$. 

Finally, 
\[
	|\lambda|\le\frac{2A\,(2R\cdot2L+2R\cdot2L)\,g}{L^2N^2}=\frac{16ARg}{LN^2}<\frac{128A}{MN}<\frac{64}{M}\le32,
\]
using $R<D=8LN/(gM)$ and $A<N/2$. Since  the function $1/(t^2y)$ has bounded derivatives of all orders on $S$, all derivatives of $\e(-\lambda/(t^2y))$ on $S$ are bounded by constants depending only on their order. Applying Leibniz's rule to the product defining $F$ in \eqref{eq:new_def_F} and using $\eta_R\le1$, we obtain \eqref{eq:weight_derivatives}.
\end{proof}
We now verify the remaining hypotheses of Lemma~\ref{lem:progression} for the parameters \eqref{eq:dictionary}, with $\delta=\eta_R$ and $C=128$. We have $N/g>M/(8L)$ by \eqref{eq:assumption on N over g}, since otherwise $\Hcal_g=0$, as shown in Section~\ref{sec:complementary_divisor}. Together with the assumptions in \eqref{eq:blockrange}, we obtain
\[
Z=N/g>M/(8L)\ge2^{23}L^3, \qquad g<8LN/M\le8L.
\]
In particular, $Z\ge1$, while $0<\eta_R\le1$ by definition.

 We have $(p_1p_1',g)=1$ by the conditions of \eqref{eq:def_K}, $(p_1p_1',p_2p_2')=1$ since the four primes are distinct, and $(p_1p_1',dd')=1$ by the definition of $\Qcal^{\pm,R}_g(1)$. Therefore, $(P,Qk)=1$. Lemma~\ref{lem:Delta_generic} gives $(\Delta,P)=1$ and $\Delta\ne0$, and 
Lemma~\ref{lem:merged_progression} gives $(\kappa,k)=1$. 

We then verify the hypotheses in \eqref{eq:progression_hypotheses}. Since $P=p_1p_1'$ and $p_1p_1'\le4L^2$, together with $Z>2^{23}L^3$, we establish $P\le Z$. To bound $\Delta$, Lemma~\ref{lem:Delta_generic} and $A<L/2$ give
\[
|\Delta|\le128ARL^3,\qquad PQk\ge L^2\cdot gL^2\cdot gR\ge L^4R>ARL^3.
\]
Hence, $|\Delta|\le128PQk$. If $\eta_R=1$, then $\eta_RQ\ge1$. If $\eta_R<1$, then since $M\ge N$ and $g,R,L\ge1$,
\[
\eta_RQ>gMRL\,g/N\ge1.
\]
For the bound on $k$, we have $k=g[d,d']<4gR^2$. If $\eta_R=1$, then the inequalities $R<D=8LN/(gM)$, $g<8L$ and $2^{11}L^3\le M^2$ imply
\[
k<4gR^2<\frac{256L^2N^2}{gM^2}=\Big(\frac{N}{g}\Big)^2\frac{256gL^2}{M^2}\le Z^2.
\]
If $\eta_R<1$, then $128\eta_RZ^2=128MRN/(gL)$. The inequality $4gR^2<128MRN/(gL)$ is equivalent to $g^2LR<32MN$, which follows from 
\[
g^2LR<8gL^2N/M\le8gL^2<64L^3\le32MN.
\]
Thus, $k\le128\eta_RZ^2$ in both cases. Lastly, $g<8L$, $R<8L/g$, and $A<L/2$ give
\[
P\le4L^2,\qquad Q<32L^3,\qquad k<256L^2,\qquad |\Delta|<2^{10}L^5.
\]
Since $Z>2^{23}L^3$, these quantities are all less than $Z^2$. In particular,
 $\max(P,Q,k,|\Delta|)\le Z^{128}$.

Together with the weight bounds in Lemma~\ref{lem:weight},
this verifies all the hypotheses of
Lemma~\ref{lem:progression}. It therefore gives, for every $\eps>0$,
\[
	|\Jcal|
	\ll_\eps
	Z^{\eps}\Big\{\frac{\eta_RZ}{Pk}+Z\sqrt{PQ\eta_R}\Big[1+\Big(\frac{Z\sqrt{(Qk,\Delta)}}{\sqrt{P}\,Qk}\Big)^{7/32}\Big]\Big\}.
\]
We may simplify this estimate using $Z\le M$, $L^2<P\le4L^2$, $gL^2<Q\le4gL^2$, and $k=g[d,d']$. Moreover, since $(\Delta,p_2p_2')=1$ by Lemma~\ref{lem:Delta_generic}, 
\[
(Qk,\Delta)=(g^2p_2p_2'[d,d'],\Delta)\le g^2[d,d'].
\]
Consequently,
\begin{align}
	|\Jcal|
	\ll_\eps
	M^{\eps}\Big\{\frac{\eta_RN}{g^2L^2[d,d']}+\frac{NL^2\sqrt{\eta_R}}{\sqrt{g}}\Big[1+\Big(\frac{N}{g^2L^3\sqrt{[d,d']}}\Big)^{7/32}\Big]\Big\}.
	\label{eq:J_generic_bound}
\end{align}
For fixed $g,M,N,A,L,R$, the right-hand side depends only on $d$ and $d'$.

It remains to sum \eqref{eq:J_generic_bound} over the variables in $\mathcal Q_g^{\pm,R}(1)$. For the sum over $d,d'$, we use the following estimate.
\begin{lemma}
\label{lem:lcm_sums}
Let $R\ge1$. For $0\le s<1$,
\begin{align}
	\sum_{R\le d,d'<2R}\frac{1}{[d,d']^{s}}\le4\zeta(2-s)\,R^{2-2s},
	\qquad
	\sum_{R\le d,d'<2R}\frac{1}{[d,d']}\le4\,(1+\log 2R).
	\label{eq:lcm_sums}
\end{align}
\end{lemma}
\begin{proof}
Group the pairs according to $r=(d,d')$ and write 
\[
d=ru,\qquad d'=ru',\qquad (u,u')=1.
\]
Then $r<2R$, and $u,u'\in [R/r,2R/r)$. For $r\le R$, this interval contains at most $R/r+1\le2R/r$ integers. For $R<r<2R$, it only contains the integer $1$. Thus, there are at most $(2R/r)^2$ pairs for each $r$. Moreover, $[d,d']=ruu'\ge r(R/r)^2=R^2/r$. Hence, 
\[
	\sum_{R\le d,d'<2R}\frac{1}{[d,d']^{s}}\le\sum_{r<2R}\Big(\frac{2R}{r}\Big)^2\Big(\frac{R^2}{r}\Big)^{-s}
	=4R^{2-2s}\sum_{r<2R}r^{s-2}.
\]
For $s<1$, the last sum is at most $\zeta(2-s)$, and for $s=1$, it is at most $1+\log2R$.
\end{proof}

We insert the bound for $|\Jcal|$ in \eqref{eq:J_generic_bound} into the sum \eqref{eq:Q_absolute} restricted to the terms of $\Qcal^{\pm,R}_g(1)$. By \eqref{eq:consequences}, every pair $n_1,n_2$ in the summation of \eqref{eq:def_J} satisfies $(n_1,gdd')=1$. Thus, for fixed primes and $d,d'$, by Lemma~\ref{lem:counting}(i), the number of quadruples $(\rho,\sigma,c,c')$ for which $\Jcal\ne0$ is at most $g\varphi(g)\le g^2$. The sum over $a_2,a_1,a_1'$ is
\[
	\sum_{a_2\in\Acal}\ \sumtwo_{a_1,a_1'\in\Acal}|\nu_{a_1}\nu_{a_1'}|
	\le|\Acal|\Big(\sum_{a\in\Acal}|\nu_a|\Big)^2\le4A^2\nn{\nu}^2,
\]
by the Cauchy-Schwarz inequality and $|\Acal|\le2A$. The number of quadruples of primes is at most $|\Pcal|^4\le16L^4$. Hence, by \eqref{eq:J_generic_bound},
\[
	\Qcal^{\pm,R}_g(1)
	\ll_\eps
	A^2\nn{\nu}^2g^2L^4M^{\eps}
	\sum_{R\le|d|,|d'|<2R}
	\Big\{\frac{\eta_RN}{g^2L^2[d,d']}+\frac{NL^2\sqrt{\eta_R}}{\sqrt{g}}+\frac{N^{39/32}L^2\sqrt{\eta_R}}{g^{15/16}L^{21/32}\,[d,d']^{7/64}}\Big\},
\]
where $d$ and $d'$ have the fixed sign. Applying Lemma~\ref{lem:lcm_sums} with $s=1,0,7/64$, respectively, and absorbing $\log2R\le\log2D\ll\log2L$ into $M^{\eps}$, we obtain
\begin{align}
	\Qcal^{\pm,R}_g(1)
	\ll_\eps
	A^2\nn{\nu}^2M^{\eps}\big\{\eta_RNL^2+g^{3/2}NL^6R^2\sqrt{\eta_R}+g^{17/16}N^{39/32}L^{171/32}R^{57/32}\sqrt{\eta_R}\big\}.
	\label{eq:Q1_bound}
\end{align}

\subsection{The remaining terms}
\label{sec:nongeneric}
In this subsection we bound $\Qcal^{\pm,R}_g(2)$, the contribution from the terms in \eqref{eq:Q_absolute} with $\Delta\ne0$ for which the primes $p_1,p_2,p_1',p_2'$ are not distinct, or are distinct but $(dd',p_1p_2p_1'p_2')>1$. Fix the sign, $R$, the residues $\rho,\sigma$, the numerator $a_2$, and the two tuples $(c,d,p_1,p_2,a_1)$ and $(c',d',p_1',p_2',a_1')$, as in Section~\ref{sec:generic}. We retain the definition of $\eta_R$ in \eqref{eq:def_eta}. 
We estimate the sum over $n_2$ for each $n_1$ separately, using Lemma~\ref{lem:completion}. This is sufficient because the number of choices of parameters $p_1,p_2,p_1',p_2',d,d'$ here is $O(L^3R^2)$, as shown below.

Recall that Lemma~\ref{lem:merged_progression} does not require the four primes to be distinct. Therefore, either $\Jcal=0$, or the two congruences in \eqref{eq:def_J} can be replaced by
\begin{align}
n_2\equiv\kappa n_1\bmod k,\label{eq:n2_congruence_in_sec4.8}
\end{align}
where $k=g[d,d']$ and $\kappa$ is a unit modulo $k$. In the latter case, by Lemma~\ref{lem:phase_difference}(i), whose hypotheses hold by the conditions of \eqref{eq:def_J},
\begin{align}
	\Jcal
	=
	\sum_{\substack{n_1\ge1\\(n_1,\,gp_1p_2p_1'p_2'[d,d'])=1}}
	\psi\Big(\frac{gn_1}{N}\Big)\,
	\Jcal(n_1),
	\label{eq:J_sum_n1}
\end{align}
where
\begin{align}
	\Jcal(n_1):=
	\sum_{\substack{n_2\equiv\kappa n_1\bmod k\\(n_2,\,p_1p_2p_1'p_2'n_1)=1}}
	G(n_2)\,
	\e\Big(\frac{\Delta_{n_1}\,\overline{gn_2}}{[p_1,p_1']\,n_1}\Big).
	\label{eq:J_by_n1}
\end{align}
For $x>0$, the weight is 
\begin{align}
	G(x):=\psi\Big(\frac{gx}{N}\Big)\,w\Big(\frac{p_1n_1-p_2x}{dM}\Big)\,w\Big(\frac{p_1'n_1-p_2'x}{d'M}\Big)\,\e\Big(-\frac{a_2(dp_1'-d'p_1)}{gp_1p_1'n_1x}\Big),
	\label{eq:def_G}
\end{align}
and $G(x)=0$ for $x\le0$. Since $\psi(gx/N)$ vanishes for $x\le N/(2g)$, the function $G$ is smooth on $\R$. The integer $\Delta_{n_1}$ is provided by Lemma~\ref{lem:phase_difference}(i) and satisfies 
\[
(\Delta_{n_1},n_1)\mid\Delta.
\]
 The condition $(n_2,g)=1$ follows from \eqref{eq:n2_congruence_in_sec4.8}.
 
 We fix $n_1$ with $(n_1,gp_1p_2p_1'p_2'[d,d'])=1$ and $\psi(gn_1/N)\ne0$, so that $N/(2g)<n_1<3N/g$. To estimate $\Jcal(n_1)$ by Lemma~\ref{lem:completion}, we first verify
the support and derivative bounds of $G$ required by the lemma. Set
\begin{align}
	\ell:=\frac{12\,\eta_RN}{g}.
	\label{eq:def_ell}
\end{align}
 If $G(x)\ne0$, then
\[
\frac{N}{2g}<x<\frac{3N}{g},
\qquad
|p_1n_1-p_2x|<3|d|M.
\]
The latter inequality restricts $x$ to an interval of length $6|d|M/p_2<12RM/L$. If $\eta_R<1$, then $12RM/L=\ell$. If $\eta_R=1$, then the support of $\psi(gx/N)$ already lies in an interval of length less than $3N/g\le\ell$. Thus, $G$ is supported in an interval of length at most $\ell$.
 
 For the derivatives, we have
\[
\Big|
\frac{d^i}{dx^i}\psi\Big(\frac{gx}{N}\Big)
\Big|
\ll_i \Big(\frac gN\Big)^i
\]
 and 
\[
\Big|
\frac{d^i}{dx^i}
w\Big(\frac{p_1n_1-p_2x}{dM}\Big)
\Big|
\ll_i
\Big(\frac{p_2}{|d|M}\Big)^i
\le
\Big(\frac{2L}{RM}\Big)^i.
\]
The same bound holds for the $i$-th derivative of $w((p_1'n_1-p_2'x)/(d'M))$. By the definition of $\eta_R$,
\[
\frac{2L}{RM}\le\frac{2g}{\eta_RN}.
\]
Indeed, if $\eta_R<1$, the right-hand side equals $2L/(MR)$. If $\eta_R=1$, then $gMR\ge LN$ gives $2L/(RM)\le2g/N$. 

For the exponential factor in $G$, write
\[
\mu_0:=a_2(dp_1'-d'p_1)/(gp_1p_1'n_1).
\]
For $x>N/(2g)$ and $i\ge1$,
\[
	\Big|\frac{d^i}{dx^i}\,\frac{\mu_0}{x}\Big|=\frac{i!\,|\mu_0|}{x^{i+1}}
	\le\frac{i!\,|\mu_0|}{x^2}\Big(\frac{2g}{N}\Big)^{i-1}.
\]
Moreover,
\[
	\frac{|\mu_0|}{x^2}<\frac{2A\cdot8RL}{gL^2\,(N/(2g))}\cdot\frac{4g^2}{N^2}=\frac{128ARg^2}{LN^3}<\frac{1024A\,g}{MN^2}<\frac{512\,g}{MN},
\]
using $|dp_1'-d'p_1|<8RL$, $R<8LN/(gM)$, and $A<N/2$. Consequently, 
\[
\Big|
\frac{d^i}{dx^i}\frac{\mu_0}{x}
\Big|
\ll_i\Big(\frac gN\Big)^i,
\]
and by the chain rule the same bound holds for the $i$-th derivative of $\e(-\mu_0/x)$. Thus, the $i$-th derivative of each of the four factors of $G$ is $\ll_i(2g/(\eta_RN))^i=(24/\ell)^i$.  Leibniz's rule therefore gives 
\[
|G^{(i)}(x)|\ll_i\ell^{-i},
\]
for $i=0,1,2$ and all $x$, with constants depending only on $w$ and $\psi$.

We next remove the condition $(n_2,p_2p_2')=1$ in \eqref{eq:J_by_n1} by M\"obius inversion, 
\[
\ind_{(n_2,p_2p_2')=1}=\sum_{f\mid(n_2,p_2p_2')}\mu(f).
\]
We keep the condition $(n_2,[p_1,p_1']n_1)=1$. For a divisor $f\mid p_2p_2'$ with
$(f,[p_1,p_1']n_1k)>1$, the corresponding sum over $n_2$ is empty, since $(n_2,[p_1,p_1']n_1k)=1$. For the remaining $f$,   write $n_2=fj$. This gives
\begin{align}
	\Jcal(n_1)
	=
	\sum_{\substack{f\mid p_2p_2'\\(f,[p_1,p_1']n_1k)=1}}\mu(f)
	\sum_{\substack{j\equiv\kappa\overline{f}n_1\bmod k\\(j,[p_1,p_1']n_1)=1}}
	G(fj)\,
	\e\Big(\frac{\Delta_{n_1}\,\overline{gf}\,\overline{j}}{[p_1,p_1']\,n_1}\Big),
	\label{eq:J_n1_mobius}
\end{align}
with at most four nonzero terms in the outer sum.

The modulus $c:=[p_1,p_1']n_1$ of the Kloosterman fraction in \eqref{eq:J_n1_mobius} need not be coprime to $k$, because a prime among $p_1,p_1'$ may divide $d$ or $d'$. To obtain coprime moduli for Lemma~\ref{lem:completion}, set
\begin{align}
	c_0:=\big([p_1,p_1'],\,k\big),\qquad c_1:=\frac{[p_1,p_1']}{c_0}\,n_1,
	\label{eq:def_c0c1}
\end{align}
so that $c=c_0c_1$. Since $[p_1,p_1']$ is squarefree and 
\[
(n_1,k)=(n_1,[p_1,p_1'])=1,
\]
we have
\[
(c_0,c_1)=1,\qquad(c_1,k)=1.
\]
 By Lemma~\ref{lem:crt_inverse},
\[
	\frac{\overline{gf}\,\overline{j}}{c}=\frac{\overline{gfj}}{c_0c_1}
	\equiv
	\frac{\overline{gfjc_1}}{c_0}+\frac{\overline{gfjc_0}}{c_1}\bmod 1.
\]
On the progression $j\equiv\kappa\overline{f}n_1\bmod k$, we have
\[
gfj\equiv g\kappa n_1\bmod c_0
\]
since $c_0\mid k$. Therefore, the factor $\e(\Delta_{n_1}\overline{gfjc_1}/c_0)$ is independent of $j$ and is a constant of modulus one. Set 
\[
u:=\Delta_{n_1}\overline{gfc_0},
\]
where the inverse is taken modulo $c_1$. Then $\Delta_{n_1}\overline{gfjc_0}\equiv u\,\overline{j}\bmod c_1$, and the inner sum in \eqref{eq:J_n1_mobius} thus has absolute value
\[
	\Big|\sum_{\substack{j\equiv\kappa\overline{f}n_1\bmod k\\(j,c_1)=1}}G(fj)\,\e\Big(\frac{u\,\overline{j}}{c_1}\Big)\Big|,
\]
where the condition $(j,c_0)=1$ has been dropped because it is implied by $j\equiv\kappa\overline f\,n_1\bmod k$, since $(\kappa\overline f\,n_1,k)=1$ and $c_0\mid k$.

We may now apply Lemma~\ref{lem:completion} with the coprime moduli $c_1$ and $k$, residue $r=\kappa\overline{f}n_1$, and weight $F(x)=G(fx)$. This weight is smooth, supported in an interval of length at most $\ell/f$, and satisfies
\[
|F^{(i)}(x)|=f^i|G^{(i)}(fx)|\ll_i(\ell/f)^{-i}
\]
 for $i=0,1,2$. The lemma gives
\[
	\Big|\sum_{\substack{j\equiv\kappa\overline{f}n_1\bmod k\\(j,c_1)=1}}G(fj)\,\e\Big(\frac{u\,\overline{j}}{c_1}\Big)\Big|
	\ll
	\frac{\ell}{f\,c_1k}\,\tau(c_1)\,(u,c_1)+\tau(c_1)^2\sqrt{c_1}.
\]
Since $\overline{gfc_0}$ is a unit modulo $c_1$, we have $(u,c_1)=(\Delta_{n_1},c_1)$. Moreover, since $c_1$ is the product of the coprime integers $[p_1,p_1']/c_0$ and $n_1$,
\[
	\frac{(u,c_1)}{c_1}
	\le\frac{([p_1,p_1']/c_0)\,(\Delta_{n_1},n_1)}{([p_1,p_1']/c_0)\,n_1}
	=\frac{(\Delta_{n_1},n_1)}{n_1}
	\le\frac{(\Delta,n_1)}{n_1},
\]
where the last inequality follows from $(\Delta_{n_1},n_1)\mid \Delta$. Furthermore, 
\[
c_1\le[p_1,p_1']n_1<\frac{12L^2N}{g},
\]
so $\tau(c_1)\ll_\eps M^{\eps}$ and $\sqrt{c_1}<2L\sqrt{3N/g}$. Summing over  $f$ in \eqref{eq:J_n1_mobius}, which has  at most four values, and using $f\ge1$, $\ell=12\eta_RN/g$, and $k=g[d,d']$, we obtain
\begin{align}
	|\Jcal(n_1)|
	\ll_\eps
	M^{\eps}\bigg\{\frac{\eta_RN}{g^2[d,d']}\cdot\frac{(\Delta,n_1)}{n_1}+L\sqrt{\frac{N}{g}}\bigg\}.
	\label{eq:J_n1_bound}
\end{align}

To sum over $n_1$, we use $(\Delta,n_1)=\sum_{s\mid(\Delta,n_1)}\varphi(s)$. Since $\Delta\ne0$,
\[
	\sum_{n_1<3N/g}\frac{(\Delta,n_1)}{n_1}
	=\sum_{s\mid|\Delta|}\varphi(s)\sum_{\substack{n_1<3N/g\\s\mid n_1}}\frac{1}{n_1}
	\le\sum_{s\mid|\Delta|}\frac{\varphi(s)}{s}\Big(1+\log\frac{3N}{g}\Big)
	\le\tau(|\Delta|)\Big(1+\log\frac{3N}{g}\Big)
	\ll_\eps M^{\eps},
\]
where the last inequality follows from $0<|\Delta|<128ARL^3<M^3$ by the size estimate in the proof of Lemma~\ref{lem:Delta_generic}, which does not require the primes to be distinct.  The second term of \eqref{eq:J_n1_bound} is summed trivially over the fewer than $3N/g$ values of $n_1$. By \eqref{eq:J_sum_n1}, \eqref{eq:J_n1_bound}, and $0\le\psi\le1$,
\begin{align}
	|\Jcal|
	\ll_\eps
	M^{\eps}\Big\{\frac{\eta_RN}{g^2[d,d']}+L\Big(\frac{N}{g}\Big)^{3/2}\Big\}.
	\label{eq:J_nongeneric_bound}
\end{align}

We now count the pairs of tuples that contribute to $\Qcal^{\pm,R}_g(2)$.
\begin{lemma}
\label{lem:nongeneric_count}
Fix the sign and the dyadic number $R$. Consider the sextuples $(p_1,p_2,p_1',p_2',d,d')$ satisfying
\[
p_1,p_2,p_1',p_2'\in\mathcal P,\qquad
p_1\ne p_2,\qquad p_1'\ne p_2',\qquad \pm d,\pm d'>0,\qquad R\le|d|,|d'|<2R.
\]
The number of such sextuples for which the four primes are not distinct, or are distinct and
\[
(dd',p_1p_2p_1'p_2')>1,
\] 
is at most $416\,L^3R^2$.
\end{lemma}
\begin{proof}
If the primes are not distinct, then at least one of the equalities
\[
p_1=p_1',\qquad p_1=p_2',\qquad p_2=p_1',\qquad p_2=p_2'
\]
holds. Each equality leaves at most $|\Pcal|^3\le(2L)^3=8L^3$ quadruples of primes. For each quadruple, there are at most $R^2$ pairs of $d,d'$, yielding at most $32L^3R^2$ sextuples. 

If the primes are distinct and $(dd',p_1p_2p_1'p_2')>1$, then one of the four primes, say $p$, divides $d$ or $d'$. The number of $d$ with $\pm d>0$, $R\le|d|<2R$, and $p\mid d$ is zero if $p>2R$, and is at most $R/p+1\le3R/L$ otherwise, since $2R\ge p>L$. For each such $d$, there are at most $R$ choices of $d'$
satisfying $\pm d'>0$ and $R\le|d'|<2R$. There are four choices of $p$, two choices of the integer it divides, and at most
\[
|\mathcal P|^4\le16L^4
\]
quadruples of primes. This gives at most 
\[
4\cdot2\cdot16L^4\cdot\frac{3R^2}{L}
=
384L^3R^2
\]
sextuples. Adding the two bounds proves the lemma.
\end{proof}

Lastly, we bound $|\Jcal|$ in each term of $\Qcal^{\pm,R}_g(2)$ by \eqref{eq:J_nongeneric_bound}. By Lemma~\ref{lem:counting}, the number of quadruples $(\rho,\sigma,c,c')$ with $\Jcal\ne0$ is at most $g^2$ for fixed tuples. The sum over $a_2,a_1,a_1'$ satisfies
\[
\sum_{a_2\in\mathcal A}
\sum_{a_1,a_1'\in\mathcal A}
|\nu_{a_1}\nu_{a_1'}|=
|\mathcal A|
\Big(\sum_{a\in\mathcal A}|\nu_a|\Big)^2\le
|\mathcal A|^2\nn{\nu}^2
\le4A^2\nn{\nu}^2.
\]
For the first term of \eqref{eq:J_nongeneric_bound}, we sum over all $16L^4$ quadruples of primes and all $d,d'$ with $R\le|d|,|d'|<2R$.  Lemma~\ref{lem:lcm_sums} with $s=1$ then gives
\[
\sum_{\substack{\pm d,\pm d'>0\\R\le|d|,|d'|<2R}}
\frac1{[d,d']}
\le4(1+\log 2R).
\]
For the second term of \eqref{eq:J_nongeneric_bound}, which is independent of $d,d'$, we use Lemma~\ref{lem:nongeneric_count}. Consequently,
\begin{align}
	\Qcal^{\pm,R}_g(2)
	&\ll_\eps
	A^2\nn{\nu}^2g^2M^{\eps}\Big\{16L^4\cdot\frac{\eta_RN}{g^2}\cdot4(1+\log2R)+416L^3R^2\cdot L\Big(\frac{N}{g}\Big)^{3/2}\Big\}
	\notag\\
	&\ll_\eps
	A^2\nn{\nu}^2M^{\eps}\big\{\eta_RNL^4+g^{1/2}L^4R^2N^{3/2}\big\}.
	\label{eq:Q2_bound}
\end{align}

\subsection{Completion of the proof of Theorem~\ref{thm:H}}
\label{sec:completion_H}
We collect the bounds of Sections~\ref{sec:delta_zero}, \ref{sec:generic} and \ref{sec:nongeneric}. By \eqref{eq:Q_three_parts}, \eqref{eq:Q0_bound}, \eqref{eq:Q1_bound} and \eqref{eq:Q2_bound}, for every fixed sign and every dyadic $R$,
\begin{align}
	\Qcal^{\pm,R}_g
	\ll_\eps
	\nn{\nu}^2M^{\eps}\Big\{&
	\frac{AL^2N^2}{g}
	+A^2\eta_RN\big(L^2+L^4\big)
	+A^2g^{3/2}NL^6R^2\sqrt{\eta_R}
	\notag\\
	&+A^2g^{17/16}N^{39/32}L^{171/32}R^{57/32}\sqrt{\eta_R}
	+A^2g^{1/2}L^4N^{3/2}R^2
	\Big\}.
	\label{eq:Q_pmR_collected}
\end{align}
We insert this into \eqref{eq:Q_dyadic}. The number of pairs $(\pm,R)$ is at most $8\log_2(2L)$, and the factor $8\log_2(2L)$ in front of \eqref{eq:Q_dyadic} is of the same size. Both factors are absorbed into $M^{\eps}$. 

The first term of \eqref{eq:Q_pmR_collected} is independent of $R$. Since $\eta_R\le1$ and $L\ge1$,  the second term is at most
\[
A^2\eta_RN(L^2+L^4)\le2A^2NL^4.
\]
For the remaining three terms, we use $\eta_R\le1$ and 
\[
	\sum_{1\le R<D}R^2<\frac{4}{3}D^2,
	\qquad
	\sum_{1\le R<D}R^{57/32}<2D^{57/32},
\]
where the second inequality follows from $2^{-57/32}<1/2$. Substituting $D=8LN/(gM)$, we have
\begin{align*}
	&g^{3/2}NL^6D^2=\frac{64\,L^8N^3}{g^{1/2}M^2},
	\qquad
	g^{1/2}L^4N^{3/2}D^2=\frac{64\,L^6N^{7/2}}{g^{3/2}M^2},
	\\
	&g^{17/16}N^{39/32}L^{171/32}D^{57/32}=\frac{8^{57/32}\,L^{57/8}N^3}{g^{23/32}M^{57/32}}.
\end{align*}
Hence
\begin{align}
	\Qcal_g
	\ll_\eps
	\nn{\nu}^2M^{\eps}\Big\{
	\frac{AL^2N^2}{g}
	+A^2NL^4
	+\frac{A^2L^6N^{7/2}}{g^{3/2}M^2}
	+\frac{A^2L^8N^3}{g^{1/2}M^2}
	+\frac{A^2L^{57/8}N^3}{g^{23/32}M^{57/32}}
	\Big\}.
	\label{eq:Q_g_five_terms}
\end{align}
The term $A^2NL^4$ can be absorbed into $A^2L^6N^{7/2}/(g^{3/2}M^2)$. Indeed, the required inequality is equivalent to $g^{3/2}M^2\le L^2N^{5/2}$. Since $N/g>M/(8L)$ by \eqref{eq:assumption on N over g}, we have 
\[
L^2N^{5/2}
>
\frac{g^{5/2}M^{5/2}}{8^{5/2}L^{1/2}}
=
g^{3/2}M^2
\frac{gM^{1/2}}{8^{5/2}L^{1/2}}.
\]
By \eqref{eq:blockrange}, $M^{1/2}\ge2^{13}L^2$. Since $g\ge1$ and $L\ge1$, we have
\[
gM^{1/2}\ge2^{13}L^2
\ge8^{5/2}L^{1/2},
\]
which leads to the required inequality. Therefore,
\begin{align}
	\Qcal_g
	\ll_\eps
	\nn{\nu}^2M^{\eps}\Big\{
	\frac{AL^2N^2}{g}
	+\frac{A^2L^6N^{7/2}}{g^{3/2}M^2}
	+\frac{A^2L^8N^3}{g^{1/2}M^2}
	+\frac{A^2L^{57/8}N^3}{g^{23/32}M^{57/32}}
	\Big\}.
	\label{eq:Q_g_final}
\end{align}

Finally, by \eqref{eq:Hg_CS}, $|\Hcal_g|\le\nn{\beta}^2\nn{\nu}\sqrt{\Qcal_g}$. Thus,
\[
	|\Hcal_g|
	\ll_\eps
	\nn{\beta}^2\nn{\nu}^2
	\Big(\frac{A^{1/2}LN}{g^{1/2}}
	+\frac{AL^3N^{7/4}}{Mg^{3/4}}
	+\frac{AL^4N^{3/2}}{Mg^{1/4}}
	+\frac{AL^{57/16}N^{3/2}}{M^{57/64}g^{23/64}}\Big)M^{\eps/2},
\]
which is \eqref{eq:bound_H} since $\eps$ is arbitrary. This was proved under the assumption $N/g>M/(8L)$. If $N/g\le M/(8L)$, then $\Hcal_g=0$ by Section~\ref{sec:complementary_divisor}, so the same bound holds. The implied constant depends only on $\eps$, $w$ and $\psi$, and the bound is therefore uniform in the sequences $\beta$ and $\nu$, as required in Section~\ref{sec:H_implies_Dw}. This completes the proof of Theorem~\ref{thm:H}.

\section{Proofs of the applications}
\label{sec:moment}

\subsection{The twisted second moment of the Riemann zeta-function}
\label{sec:app_zeta}
We prove Theorem~\ref{thm:zeta}. As in its statement, $\alpha$ and $\beta$ denote complex numbers in this subsection.

\begin{proof}[Proof of Theorem~\ref{thm:zeta}]
We follow the proof of \cite[Theorem~1.1]{PR} for coefficients $c_n\ll_\eps n^\eps$, which is modeled on \cite{BCR}. With the length $Y$, every error term there is $O(T^{1/2+\eps}Y^{1/2})$, with the exception of the terms with $0<|h|<H_d$ in \cite[Section~3]{PR}, where the estimate of \cite{BC} is applied. In the notation of \cite{PR}, for $d\le Y$ and for each of the $O((\log T)^3)$ parameters $N_1,N_2\le Y$ and $M$ of the partition of unity, these terms are bounded by
\begin{align}
	\frac{T^\eps}{d}\int\!\!\int\!\!\int_{(\eps)}
	|\mathcal M(\Phi,w)|\,
	\bigg|\sum_{0<|a|<A_0}\nu_{x,y}(a)
	\sumtwo_{\substack{n_1,n_2\\(n_1,n_2)=1}}b_{n_1}b'_{n_2}\,\e\Big(-\frac{a\overline{n_1}}{n_2}\Big)\bigg|
	\,|dw|\,dy\,dx,
	\qquad
	A_0:=\frac{N_1N_2}{d^2\,T^{1-\eps}}.
	\label{eq:zeta_offdiag}
\end{align}
Here $x$ and $y$ run over a rectangle of area $O(Td^2/(N_1N_2))$: the supports of $F_M$ and $\Phi$ give $x\asymp dM/N_2$ and $y\asymp Td/(MN_1)$ with absolute implied constants, and this restriction is imposed before the Mellin inversion of $\Phi$. Further, $\mathcal M(\Phi,w)$ is the Mellin transform of $\Phi$, $|\nu_{x,y}(a)|\le\tau(|a|)$, and the coefficients $b_{n_1}$, $b'_{n_2}$, which depend on $x$ and $w$, are supported on $N_1/(2d)\le n_1\le3N_1/d$ and $N_2/(2d)\le n_2\le3N_2/d$ and are $O(T^\eps)$. By the hypothesis on $\Phi$, the integral of $|\mathcal M(\Phi,w)|$ over the line $\Re w=\eps$ is $O((\log T)^2)$.

Set $X=\max(N_1,N_2)/d$. Splitting $a$, $n_1$ and $n_2$ into dyadic ranges, we write the sum inside the absolute value in \eqref{eq:zeta_offdiag} as a combination of $O(\log T)$ sums \eqref{eq:def_B_AMN}, or their complex conjugates, with numerators of size $A\le A_0$, with lengths between $X/4$ and $2X$ for the longer variable, and with coefficient sequences whose $\ell^2$-norms have product $O(T^\eps(N_1N_2)^{1/2}d^{-1}A^{1/2})$. Since $N_1N_2/d^2\le X^2$ and $X\le Y\le T^{14/27}$,
\[
	A_0\le X^2\,T^{-1+\eps}\le X^{1/10}\,T^{-2/135+\eps},
\]
so that $A\le(X/4)^{1/10}$ if $\eps$ is small and $T$ is large, as we may assume. If both lengths are at least $2$, Theorem~\ref{thm:main} bounds each of these sums by
\[
	O\big(T^\eps\,(N_1N_2)^{1/2}d^{-1}A_0^{7/8}X^{15/16}\big),
\]
and if one of the lengths is less than $2$ the same bound follows from the Cauchy-Schwarz inequality, since $A\le X^{1/10}$. Hence \eqref{eq:zeta_offdiag} is
\[
	\ll\frac{T^\eps}{d}\cdot\frac{Td^2}{N_1N_2}\cdot\frac{(N_1N_2)^{1/2}}{d}\,A_0^{7/8}X^{15/16}
	\ll T^{1/8+\eps}\,(N_1N_2)^{3/8}\,d^{-7/4}\,X^{15/16}
	\le T^{1/8+\eps}\,Y^{27/16}\,d^{-43/16}.
\]
Summing over $d$ and over the parameters $N_1,N_2,M$ gives the first error term in \eqref{eq:zeta_result}. The corresponding terms of the second sum of the approximate functional equation are treated in \cite{PR} in the same way and satisfy the same bound.

It remains to justify the uniformity in the shifts, since the weights of the approximate functional equation in \cite{PR} contain the factor $1-4s^2/(\alpha+\beta)^2$, and the estimates for them are not uniform as $\alpha+\beta\to0$. Let $C>0$ be such that $|\alpha|,|\beta|\le C/\log T$. For shifts with $|\alpha+\beta|=(2C+1)/\log T$ and $|\alpha|,|\beta|\le(3C+1)/\log T$ these estimates lose only powers of $\log T$, and the argument above proves \eqref{eq:zeta_result}. Now fix $\beta$ with $|\beta|\le C/\log T$, and let $T$ be large. The difference between $I(\alpha,\beta)$ and the main term of \eqref{eq:zeta_result} is a holomorphic function of $\alpha$ in the disc $|\alpha+\beta|\le(2C+1)/\log T$. This disc contains every $\alpha$ with $|\alpha|\le C/\log T$, and on its boundary $|\alpha|\le(3C+1)/\log T$. By the maximum modulus principle, \eqref{eq:zeta_result} holds in the whole disc.
\end{proof}

\subsection{Determinant equations}
\label{sec:app_det}
We prove Theorem~\ref{thm:det}, by means of Corollary~\ref{cor:bilinear}.

\begin{proof}[Proof of Theorem~\ref{thm:det}]
We apply \cite[Corollary~1]{BC} with its four lengths equal to $2N$ and with the weights $f(x/N)$, $g(x/N)$; its error term is then $O_{f,g,\eps}(\nn{\alpha}\nn{\beta}N^{19/20+\eps})$ for every nonzero integer $\Delta$. It therefore suffices to prove \eqref{eq:det_result} with the error term $O(|\Delta|^{3/8}\nn{\alpha}\nn{\beta}N^{15/16+\eps})$ when $|\Delta|\le N^{1/30}$, which we assume from now on.

Let $n_1,n_2\sim N$, set $d=(n_1,n_2)$ and write $n_1=du$, $n_2=dv$, so that $(u,v)=1$. The equation $m_1n_2-m_2n_1=\Delta$ has no solution unless $d\mid\Delta$. If $d\mid\Delta$, set $\delta=\Delta/d$. The equation becomes $m_1v-m_2u=\delta$, and a pair of integers $(m_1,m_2)$ satisfies it if and only if
\[
	m_2\equiv-\delta\overline{u}\bmod v
	\qquad\text{and}\qquad
	m_1=\frac{um_2+\delta}{v},
\]
where $\overline{u}$ is the inverse of $u$ modulo $v$. Therefore
\begin{align}
	\mathcal S_\Delta(N)
	=
	\sum_{d\mid\Delta}\,
	\sumtwo_{\substack{u,v\sim N/d\\(u,v)=1}}\alpha_{du}\beta_{dv}
	\sum_{m\equiv-\delta\overline{u}\bmod v}F_{u,v}(m),
	\qquad
	F_{u,v}(x):=f\Big(\frac{ux+\delta}{vN}\Big)g\Big(\frac{x}{N}\Big),
	\label{eq:det_progression}
\end{align}
where $d$ runs over the positive divisors of $\Delta$ and $\delta=\Delta/d$.

For a smooth function $F$ with compact support and integers $c$ and $v\ge1$, Poisson summation applied to $x\mapsto F(c+vx)$ gives $\sum_{m\equiv c\bmod v}F(m)=v^{-1}\sum_{h\in\Z}\e(hc/v)\widehat F(h/v)$, where $\widehat F(\xi)=\int_\R F(x)\e(-\xi x)\,dx$. Set $\rho=\Delta/N^2$ and, for $d\mid\Delta$, $h\in\Z$ and $s,t>0$,
\begin{align}
	K_{d,h}(s,t)
	:=
	\frac dt\int_\R f\Big(\frac{sz+\rho}{t}\Big)g(z)\,\e\Big(-\frac{hdz}{t}\Big)\,dz.
	\label{eq:det_kernel}
\end{align}
The substitution $x=Nz$ shows that $v^{-1}\widehat{F_{u,v}}(h/v)=K_{d,h}(du/N,dv/N)$, and \eqref{eq:det_progression} becomes
\begin{align}
	\mathcal S_\Delta(N)
	=
	\sum_{d\mid\Delta}\sum_{h\in\Z}\,
	\sumtwo_{\substack{u,v\sim N/d\\(u,v)=1}}\alpha_{du}\beta_{dv}\,
	\e\Big(-\frac{h\delta\overline{u}}{v}\Big)
	K_{d,h}\Big(\frac{du}{N},\frac{dv}{N}\Big).
	\label{eq:det_poisson}
\end{align}
For $h=0$ we have
\[
	K_{d,0}\Big(\frac{du}{N},\frac{dv}{N}\Big)
	=\frac1v\int_\R F_{u,v}(x)\,dx
	=\frac{d}{n_1n_2}\int_\R f\Big(\frac{x+\Delta}{Nn_2}\Big)g\Big(\frac{x}{Nn_1}\Big)\,dx,
\]
where the second equality follows by replacing $x$ by $x/n_1$. The pairs $(n_1,n_2)$ with $n_1,n_2\sim N$ and $(n_1,n_2)=d$ correspond to the coprime pairs $(u,v)$ with $u,v\sim N/d$, so the terms with $h=0$ in \eqref{eq:det_poisson} add up to $\mathcal R_\Delta(N)$.

It remains to bound the terms with $h\ne0$. We first show that for all integers $0\le i,j\le2$ and $J\ge0$,
\begin{align}
	\partial_s^i\partial_t^jK_{d,h}(s,t)\ll_{f,g,J}d\,(|h|d)^{-J}
	\qquad(h\ne0,\ \tfrac12\le s,t\le\tfrac52).
	\label{eq:det_K_bound}
\end{align}
Note that $|\rho|\le1$, since $|\Delta|\le N^{1/30}$. Differentiating under the integral sign in \eqref{eq:det_kernel}, we see that $\partial_s^i\partial_t^jK_{d,h}(s,t)$ is $d$ times a sum of boundedly many integrals of the form
\[
	(hd)^k\int_\R G(s,t,z)\,\e\Big(-\frac{hdz}{t}\Big)\,dz,
	\qquad 0\le k\le j,
\]
where $G$ is smooth, vanishes unless $z$ lies in the support of $g$, and has all its derivatives with respect to $z$ bounded in terms of $f$, $g$ and the order of the derivative, uniformly for $\tfrac12\le s,t\le\tfrac52$ and $|\rho|\le1$. Integrating by parts $J+k$ times with respect to $z$ produces the factor $(t/(2\pi|h|d))^{J+k}$ and no boundary terms, which gives \eqref{eq:det_K_bound}.

Next we separate the variables $s$ and $t$. Let $\phi$ be a fixed smooth function supported in $[3/4,9/4]$ and equal to $1$ on $[1,2]$. The function $\phi(s)\phi(t)K_{d,h}(s,t)$ is smooth and supported in $[3/4,9/4]^2$, so its extension to $\R^2$ with period $2$ in each variable is smooth, and
\begin{align}
	\phi(s)\phi(t)K_{d,h}(s,t)
	=
	\sum_{k_1,k_2\in\Z}c_{k_1,k_2}\,\e\Big(\frac{k_1s}{2}\Big)\e\Big(\frac{k_2t}{2}\Big)
	\qquad(\tfrac12\le s,t\le\tfrac52),
	\label{eq:det_fourier}
\end{align}
with
\[
	c_{k_1,k_2}=c_{k_1,k_2}(d,h)
	:=\frac14\int_{1/2}^{5/2}\int_{1/2}^{5/2}\phi(s)\phi(t)K_{d,h}(s,t)\,\e\Big(-\frac{k_1s}{2}-\frac{k_2t}{2}\Big)\,ds\,dt.
\]
Integrating by parts twice with respect to $s$ if $k_1\ne0$, and twice with respect to $t$ if $k_2\ne0$, we obtain from \eqref{eq:det_K_bound} that $c_{k_1,k_2}\ll_{f,g,J}d\,(|h|d)^{-J}(1+|k_1|)^{-2}(1+|k_2|)^{-2}$. Hence the series in \eqref{eq:det_fourier} converges absolutely, and
\begin{align}
	\sum_{k_1,k_2\in\Z}|c_{k_1,k_2}(d,h)|\ll_{f,g,J}d\,(|h|d)^{-J}
	\qquad(h\ne0).
	\label{eq:det_c_bound}
\end{align}
For $u,v\sim N/d$ the points $s=du/N$ and $t=dv/N$ lie in $(1,2]$, where $\phi=1$.

Set $H=N^{1/15}$, and let $d\mid\Delta$ and $1\le|h|\le H$. We have $N/d\ge2$: this is clear if $d=1$, and if $d\ge2$ then $N^{1/30}\ge|\Delta|\ge2$ and $N/d\ge N^{29/30}\ge2$. Moreover, $|h\Delta|/d\le N^{1/15}N^{1/30}/d=N^{1/10}/d\le(N/d)^{1/10}$. For each $(k_1,k_2)$ we may therefore apply Corollary~\ref{cor:bilinear} with both lengths equal to $N/d$, with the numerator $|h\delta|=|h\Delta|/d$, and with the coefficients $\alpha_{du}\e(k_1du/(2N))$ and $\beta_{dv}\e(k_2dv/(2N))$, after taking complex conjugates if $h\delta>0$. Their $\ell^2$-norms are at most $\nn{\alpha}$ and $\nn{\beta}$. By \eqref{eq:det_fourier}, \eqref{eq:det_c_bound} with $J=2$, and \eqref{eq:bilinear_bound},
\begin{align*}
	\sumtwo_{\substack{u,v\sim N/d\\(u,v)=1}}\alpha_{du}\beta_{dv}\,
	\e\Big(-\frac{h\delta\overline{u}}{v}\Big)
	K_{d,h}\Big(\frac{du}{N},\frac{dv}{N}\Big)
	&\ll_{f,g,\eps}
	\frac{1}{d\,h^2}\,\nn{\alpha}\nn{\beta}\Big(\frac Nd\Big)^{15/16+\eps}\Big(\frac{|h\Delta|}{d}\Big)^{3/8}\\
	&\le
	|\Delta|^{3/8}\,\nn{\alpha}\nn{\beta}\,N^{15/16+\eps}\,d^{-37/16}|h|^{-13/8}.
\end{align*}
The sum of $d^{-37/16}|h|^{-13/8}$ over $d\ge1$ and $h\ne0$ converges, so the terms with $1\le|h|\le H$ in \eqref{eq:det_poisson} contribute $O(|\Delta|^{3/8}\nn{\alpha}\nn{\beta}N^{15/16+\eps})$.

For $|h|>H$ we use \eqref{eq:det_K_bound} with $i=j=0$. There are at most $2N/d$ integers $u\sim N/d$, so the Cauchy-Schwarz inequality gives $\sum_{u,v\sim N/d}|\alpha_{du}\beta_{dv}|\le(2N/d)\nn{\alpha}\nn{\beta}$. With $J=2$ the terms with $|h|>H$ in \eqref{eq:det_poisson} therefore contribute
\[
	\ll_{f,g}\nn{\alpha}\nn{\beta}\,N\sum_{d\ge1}\sum_{|h|>H}(|h|d)^{-2}
	\ll\nn{\alpha}\nn{\beta}\,N^{14/15},
\]
which is $O(\nn{\alpha}\nn{\beta}N^{15/16})$, and \eqref{eq:det_result} follows.
\end{proof}

\subsection{Convolutions in arithmetic progressions}
\label{sec:app_conv}
We prove Theorems~\ref{thm:conv} and~\ref{thm:nb} and Corollary~\ref{cor:conv}, with the notation and the conventions of Section~\ref{sec:new_results}.

\begin{proof}[Proof of Theorem~\ref{thm:conv}]
Suppose first that $|a|>x^{1/40}$. Then \eqref{eq:conv_range} gives
\[
	Q^{33}N^{36}=(Q^{27}N^{35})^{11/9}N^{-61/9}\le x^{154/9}\,|a|^{-22/3}\le x^{17-13/180},
\]
so that $N\le Q^{-11/12}x^{17/36-1/500}$, and \cite[Corollary~1.1(i)]{FR} applies with $\min(\eps,1/500)$ in place of $\eps$. We may therefore assume that $|a|\le x^{1/40}$.

We follow \cite{FR} and indicate the changes; the notation is that of \cite{FR}, with $X=MN$. In \cite[Section~8]{FR}, Corollary~1.1(i) of \cite{FR} is deduced from \cite[Theorem~1.1]{FR}, and in the proof of that theorem the estimate of \cite{BC} is used only for the sum $W^{\mathrm{Err1}}$, for which it gives \cite[(49)]{FR}
\[
	W^{\mathrm{Err1}}\ll D^{C}X^{\eps}\big(M^{3/20}N^{59/20}Q^{33/20}+N^{23/8}Q^{15/8}\big).
\]
We show that, if
\begin{align}
	8(\log2X)^4\,|a|\,(NQ)^{19/10}\le X,
	\label{eq:conv_cap}
\end{align}
then
\begin{align}
	W^{\mathrm{Err1}}\ll D^{C}X^{\eps}\,|a|^{3/8}M^{1/8}N^{53/16}Q^{27/16}.
	\label{eq:conv_W}
\end{align}
The sum $W^{\mathrm{Err1}}$ carries the coefficients $\beta_{dd_1\nu_1'}\overline{\beta_{d\nu_2}}$, which are $O(X^{\eps})$. By \cite[(48)]{FR}, it is therefore at most $D^{C}X^{\eps}MQ^{-2}$ times a sum over $\nu_1',\nu_2\le2N$ of
\begin{align}
	\bigg|\sum_{1\le h\le H}\ \sum_{k_1'\le2Q}\ \sum_{k_2'\le2Q}\eta_0(h)\eta_1(k_1')\eta_2(k_2')\,
	\e\Big(\frac{\vartheta h\,\overline{m}}{n}\Big)\bigg|,
	\label{eq:conv_trilinear}
\end{align}
where $\vartheta=a(d_1\nu_1'-\nu_2)$, $m=\gamma dd_1\nu_2k_1'$, $n=\nu_1'k_2'$ and $H=4(\log2X)^4Q^2/M$, the coefficients are at most $1$ in modulus, the integers $d,d_1\le D$ and $\gamma\le D^3$ are fixed, and the variables are subject to \cite[(26), (27)]{FR}. These conditions give $(m,n)=1$ and $\vartheta\ne0$, and they show that for fixed $\nu_1',\nu_2$ the integers $m$ and $n$ run over subsets of intervals $(M',2M']$ and $(N',2N']$ with
\[
	NQ<M'\le2D^2NQ,\qquad D^{-4}NQ\le N'\le2NQ.
\]
We split $h$ into dyadic ranges $(H_0,2H_0]$. For each of them, the sum in \eqref{eq:conv_trilinear} is of the form \eqref{eq:def_B_AMN}, or its complex conjugate, with numerators $|\vartheta|h$ in $(|\vartheta|H_0,2|\vartheta|H_0]$ and with coefficient sequences of $\ell^2$-norms $O(Q^{1/2})$, $O(Q^{1/2})$ and $O(H_0^{1/2})$. Since $|d_1\nu_1'-\nu_2|\le2N$, the hypothesis \eqref{eq:conv_cap} gives
\[
	|\vartheta|H_0\le2|a|NH=8(\log2X)^4\,|a|\,\frac{(NQ)^2}{X}\le(NQ)^{1/10}\le\max(M',N')^{1/10}.
\]
Theorem~\ref{thm:main} therefore bounds \eqref{eq:conv_trilinear} by $O(D^{2}X^\eps\,QH^{1/2}(NQ)^{15/16}(|a|NH)^{3/8})$. Summing over $\nu_1',\nu_2$ we get
\[
	W^{\mathrm{Err1}}\ll D^{C}X^\eps\,MQ^{-2}\cdot N^2\cdot QH^{1/2}\,(NQ)^{15/16}\,(|a|NH)^{3/8},
\]
which is \eqref{eq:conv_W} by the definition of $H$.

Hence, under \eqref{eq:conv_cap}, Theorem~1.1 of \cite{FR} holds with the right-hand side of \eqref{eq:conv_W} in place of the last term of its bound. Assume now \eqref{eq:conv_range}, let $X$ be large, and take $D=\exp((\log2X)^{\eps/2}/10)$, a choice of the same kind as in the proof of \cite[Proposition~8.1(i)]{FR}. Recall that $x\le4X$. Since $NQ\ll X^{14/27}$ and $N\ll X^{2/5}$, the hypotheses $M>QX^{\eps'}$, for a small $\eps'>0$, and $M>N>D^{10}$ of \cite[Theorem~1.1]{FR} are satisfied. Since $(NQ)^{27}\le x^{14}|a|^{-6}$ and $|a|\le x^{1/40}$, we have
\[
	|a|\,(NQ)^{19/10}\le|a|^{26/45}\,x^{133/135}\le x^{1-1/3000},
\]
which gives \eqref{eq:conv_cap}. Finally, $|a|^6Q^{27}N^{35}\le x^{14-\eps}$ implies $|a|^{3/8}M^{1/8}N^{53/16}Q^{27/16}\ll MN^2X^{-\eps/16}$, so that the right-hand side of \eqref{eq:conv_W}, with its $\eps$ small enough, is $O(MN^2X^{-\eps/20})$. The rest of the proof of \cite[Corollary~1.1]{FR}, including the passage to the condition $x<mn\le2x$, is unchanged.
\end{proof}

\begin{proof}[Proof of Corollary~\ref{cor:conv}]
We may assume that $\eps<14/27$, since otherwise only the moduli $q\le2$ occur, for which the differences inside the absolute values vanish. We may also assume that $x$ is large, so that $x^\eps\le x/12$. In \cite[Sections~9-11]{FR} the three corollaries are deduced from Corollary~1.1(i) of that paper, which is applied with the same residue $a$ and with a short factor of length $N$ in a range $\exp((\log x)^{\eta_1})\le N\le x^{\eta_2}$; here $\eta_1,\eta_2>0$ are parameters of the construction, which can be chosen small in terms of $\eps$ and $k$. If $Q\le x^{14/27-\eps}$, $|a|\le x^\eps$ and $\eta_2\le\eps/2$, then
\[
	|a|^6\,Q^{27}N^{35}\le x^{14-3\eps},
\]
so that Theorem~\ref{thm:conv}, with $\min(\eta_1,3\eps)$ in place of $\eps$, can be used in place of \cite[Corollary~1.1(i)]{FR}.
\end{proof}

\begin{proof}[Proof of Theorem~\ref{thm:nb}]
We may assume that $X$ is large. By \eqref{eq:nb_range} we have $X^{1/2}\le N$ and $|a|^6N^{27}\le X^{14}$, so that
\begin{align}
	X^{13/27}\le M\le X^{1/2}\le N\le X^{14/27},\qquad |a|\le X^{1/12}.
	\label{eq:nb_sizes}
\end{align}
In \cite{FR2} the sequences are real and $n\sim N$ stands for $N\le n<2N$. Since $E_0(q)$ is bilinear in $\alpha$ and $\beta$, and since the real and imaginary parts of $\alpha$ and $\beta$ satisfy the hypotheses of the theorem, we may assume that $\alpha$ and $\beta$ are real. The terms with $m=2M$ or $n=2N$ contribute $O(X^{\eps}(M+N+Q))=O(X^{3/5})$ to the left-hand side of \eqref{eq:nb_goal}, by the divisor bound applied to $|mn-a|$ when $mn\ne a$, and the sequence $\beta$ without its term $n=2N$ still satisfies the Siegel-Walfisz condition. We may therefore use the conventions of \cite{FR2}.

Set $\eps_0=\min(\sigma,1/5000)$. If $Q\ge N^{56/23}X^{-17/23+\eps_0}$, then \eqref{eq:nb_goal} is \cite[Theorem~1.1]{FR2} with $\eps_0$ in place of $\eps$. We assume from now on that
\begin{align}
	Q<N^{56/23}X^{-17/23+\eps_0}.
	\label{eq:nb_compl}
\end{align}
Together with $Q\ge X^{1/2}$ this gives $N>X^{57/112-\eps_0/2}$, and together with \eqref{eq:nb_range} it gives
\[
	|a|^6\le X^{14}Q^8N^{-35}<X^{186/23+8\eps_0}N^{-357/23}.
\]
Hence
\begin{align}
	|a|\,N^{19/10}<X^{31/23+2\eps_0}N^{-79/115}<X^{1-1/560+3\eps_0}\le X^{1-1/1000}.
	\label{eq:nb_cap}
\end{align}

We follow the proof of \cite[Theorem~1.1]{FR2}, in its notation and with its parameter $\eps$ small in terms of $\sigma$; in particular $H=X^{\eps}Q/M$ and $R=2N/Q$. By \eqref{eq:nb_sizes} and \eqref{eq:nb_range}, the assumptions \cite[(8)]{FR2} and $Q\le NX^{-\eps}$ hold. The estimate of \cite{BC} is used in \cite{FR2} only for the sum $\widetilde W^{\mathrm{Err1}}(Q)$, for which it gives, see \cite[(41), (43)]{FR2},
\[
	\widetilde W^{\mathrm{Err1}}(Q)\ll X^{68\eps}\big(M^{3/20}N^{79/20}Q^{-43/20}+N^{31/8}Q^{-2}\big).
\]
We show that, with an absolute constant $C$,
\begin{align}
	\widetilde W^{\mathrm{Err1}}(Q)\ll X^{C\eps}\,|a|^{3/8}M^{1/8}N^{53/16}Q^{-3/2}.
	\label{eq:nb_W}
\end{align}
The last two conditions of \cite[(43)]{FR2} are then replaced by
\[
	MQ\cdot X^{C\eps}\,|a|^{3/8}M^{1/8}N^{53/16}Q^{-3/2}\ll M^2N^2X^{-\eps/4},
\]
that is, by $|a|^6N^{35}\ll X^{14-(16C+4)\eps}Q^8$, which follows from \eqref{eq:nb_range} if $\eps$ is small enough in terms of $\sigma$. The other terms are treated in \cite{FR2} for all $1\le|a|\le X$. For the convenience of the reader we record three details of that treatment. In the bound \cite[(19)]{FR2} for $V^{\mathrm{Err1}}(Q)$ the condition $(q,a)=1$ is kept by applying \cite[Lemma~2.4]{FR2} with $b=|a|n_2$, and the dependence on $a$ is controlled by the inequality $\sum_{n\sim N}(ah,n)^{1/2}\le3N\tau(|ah|)$. In \cite[(32)]{FR2} the variance estimate \cite[Lemma~2.3]{FR2} is used with $N^2(\log2N)^{-A}$ on the right-hand side. The theorem follows.

It remains to prove \eqref{eq:nb_W}. By \cite[(36), (37)]{FR2} there are integers $d,d_1\le X^\eps$ with $d_1\mid d^\infty$ such that $\widetilde W^{\mathrm{Err1}}(Q)\ll X^{2\eps}MQ^{-1}|\mathcal W|$, where, with $t=dd_1$ and with $u$, $v$ in place of $\nu_1'$, $\nu_2$,
\begin{align}
	\mathcal W=\sum_{1\le|r|\le R/d}\ \sumtwo_{u,v}\beta_{tu}\beta_{dv}\sum_{0<|h|\le H}\Phi_r(h,u,v)\,\e\Big(\frac{ah\,\overline{tu}}{q}\Big),
	\qquad q:=\frac{d_1u-v}{r}.
	\label{eq:nb_calW}
\end{align}
Here $u$ and $v$ run over the integers with $tu\sim N$ and $dv\sim N$ for which $q$ is an integer and
\begin{align}
	(d_1u,v)=1,\qquad (u,d)=1,\qquad (q,tu)=1,\qquad (q,a)=1,
	\label{eq:nb_cond}
\end{align}
and $\Phi_r$ is the product of the two smooth factors in \cite[(37)]{FR2}; it vanishes unless $Q/2\le q\le5Q/2$. The last condition in \eqref{eq:nb_cond} comes from the definition of $W(Q)$ in \cite{FR2}. By \eqref{eq:nb_cond} the integers $u$ and $q$ are coprime to $t$, and $(tv,u)=1$. As in \cite{FR2}, two applications of the reciprocity law give
\begin{align}
	\frac{ah\,\overline{tu}}{q}\equiv-\frac{ah\,\overline{u}\,\overline{q}}{t}+\frac{ahr\,\overline{tv}}{u}+\frac{ahr}{tu(d_1u-v)}\bmod 1,
	\label{eq:nb_recip}
\end{align}
where $\overline{u}$ and $\overline{q}$ are inverses modulo $t$. We fix the classes of $u$ and of $q$ modulo $t$, say $u\equiv u_0$ and $q\equiv q_0\bmod t$ with $(u_0q_0,t)=1$; there are at most $t^2\le X^{4\eps}$ choices. When $(r,t)>1$ the class of $q$ is not determined by those of $r$, $u$ and $v$. The first term on the right of \eqref{eq:nb_recip} then contributes a factor $\xi_h$ of modulus $1$ that depends only on $h$. The third term is included in the weight: we set
\[
	\Psi_r(h,u,v):=\Phi_r(h,u,v)\,\e\Big(\frac{ahr}{tu(d_1u-v)}\Big).
\]
By \cite[(39)]{FR2}, for $e_0,e_1,e_2\in\{0,1\}$,
\begin{align}
	\partial_h^{e_0}\partial_u^{e_1}\partial_v^{e_2}\Psi_r(h,u,v)\ll X^{50\eps}\,|h|^{-e_0}\Big(\frac{t}{|r|Q}\Big)^{e_1}\Big(\frac{d}{|r|Q}\Big)^{e_2}.
	\label{eq:nb_deriv}
\end{align}

Next we express the conditions on $q$ without losing a power of $|a|$. The condition $(q,a)=1$ is detected by the sum of $\mu(\ell)$ over the positive integers $\ell$ dividing both $a$ and $q$, and only divisors $\ell$ coprime to $t$ occur, since $(q,t)=1$. For such $\ell$ let $q_1$ be the class modulo $\ell t$ with $q_1\equiv q_0\bmod t$ and $q_1\equiv0\bmod\ell$. The three conditions $r\mid d_1u-v$, $q\equiv q_0\bmod t$ and $\ell\mid q$ together are equivalent to $d_1u-v\equiv rq_1\bmod|r|\ell t$, and
\[
	\ind_{d_1u-v\equiv rq_1\bmod|r|\ell t}=\frac{1}{|r|\ell t}\sum_{j\bmod|r|\ell t}\e\Big(\frac{j(d_1u-v-rq_1)}{|r|\ell t}\Big).
\]
The conditions $(q,t)=1$ and $(u,d)=1$ are contained in the choice of the classes $q_0$ and $u_0$, and $(q,u)=1$ follows from $(u,v)=1$. Therefore, for some $u_0$, $q_0$ and $\ell$, and for some $j=j(r)$,
\begin{align}
	|\mathcal W|\le X^{4\eps}\,\tau(|a|)\sum_{1\le|r|\le R/d}|\mathcal W_r|,
	\qquad
	\mathcal W_r:=\sumtwo_{\substack{u,v\\(tv,u)=1}}b_ub'_v\sum_{0<|h|\le H}\xi_h\,\Psi_r(h,u,v)\,\e\Big(\frac{ahr\,\overline{tv}}{u}\Big),
	\label{eq:nb_sep}
\end{align}
where
\[
	b_u=\beta_{tu}\,\e\Big(\frac{jd_1u}{|r|\ell t}\Big)\ind_{u\equiv u_0\bmod t},
	\qquad
	b'_v=\beta_{dv}\,\e\Big(\frac{-jv}{|r|\ell t}\Big)\ind_{(v,d_1)=1}.
\]
These coefficients are supported on $tu\sim N$ and $dv\sim N$, and they are $O(X^\eps)$.

We now fix $r$, the sign of $h$ and a range $H_0<|h|\le2H_0$, where $2H_0$ is a power of $2$ and $H_0<H$, and we localize $u$ and $v$ as in \cite[Section~4]{Wr2}. Since $tu-dv=drq$ and $\Psi_r$ vanishes unless $Q/2\le q\le5Q/2$, only pairs with $|tu-dv|\le5d|r|Q/2$ occur. We split $[N,2N)$ into at most $3N/(d|r|Q)$ intervals of length $d|r|Q$; if $tu$ lies in one of them, then $dv$ lies in one of at most seven of them. For each such pair of intervals, $u$ runs over an interval of length $|r|Q/d_1$ and $v$ over an interval of length $|r|Q$, so that by \eqref{eq:nb_deriv} the weight $\Psi_r$ is removed by partial summation in $h$, $u$ and $v$ at the cost of a factor $X^{C\eps}$. In \cite{FR2} the partial summation is carried out over the full ranges of $u$ and $v$, which costs the factor $N^2/(rQ)^2$ in \cite[(40)]{FR2}.

The sums that remain are of the form \eqref{eq:def_B_AMN}, or complex conjugates of such sums, with the variables $m=tv$ and $n=u$ and with the numerators $|ar|\,|h|$. The integers $m$ lie in $[d_1N,2d_1N)$ and the integers $n$ in $[N/t,2N/t)$, and the integers of an interval $[Y,2Y)$ lie in $(Y',2Y']$ with $Y'=\lceil Y\rceil-1/2$. Thus Theorem~\ref{thm:main} can be applied with lengths $M'$ and $N'$ such that $|M'-d_1N|\le1/2$ and $|N'-N/t|\le1/2$, and with $|ar|H_0$ in place of $A$. We have $M',N'\ge2$ and $N/2\le\max(M',N')\le2X^\eps N$, and by \eqref{eq:nb_cap},
\[
	|ar|H_0\le|a|\,\frac Rd\,H\le2X^\eps\,|a|\,\frac{N^2}{X}<2X^{\eps-1/1000}N^{1/10}\le\max(M',N')^{1/10}
\]
if $\eps<1/1000$. On a pair of intervals there are at most $|r|Q+1$ values of $v$ and at most $|r|Q/d_1+1$ values of $u$, so that the $\ell^2$-norms of the two coefficient sequences have product $O(X^{2\eps}|r|Q)$, and the sequence $\xi_h$ has norm $O(H_0^{1/2})$. Theorem~\ref{thm:main} therefore bounds each of these sums by
\[
	O\big(X^{C\eps}\,|r|Q\,H_0^{1/2}\,N^{15/16}\,(|ar|H_0)^{3/8}\big).
\]
Summing over the ranges of $h$ and over the at most $21N/(d|r|Q)$ pairs of intervals, we obtain $\mathcal W_r\ll X^{C\eps}|a|^{3/8}|r|^{3/8}H^{7/8}N^{31/16}$. Since $\tau(|a|)\ll X^\eps$, the first part of \eqref{eq:nb_sep} gives
\[
	|\mathcal W|\ll X^{C\eps}\,|a|^{3/8}H^{7/8}N^{31/16}R^{11/8}\ll X^{C\eps}\,|a|^{3/8}M^{-7/8}N^{53/16}Q^{-1/2},
\]
and \eqref{eq:nb_W} follows from $\widetilde W^{\mathrm{Err1}}(Q)\ll X^{2\eps}MQ^{-1}|\mathcal W|$.
\end{proof}

\section{Conclusion and future work}
\label{sec:conclusion}

\subsection{The critical range}
\label{sec:critical_range}
Theorem~\ref{thm:main} is of interest when $M$ and $N$ are of comparable size. Let $N\le M$. The bound \eqref{eq:main_bound} is smaller than the trivial bound only if $N>M^{7/8}A^{-1/4}$, and for a single numerator $a\ll MN$ it is smaller than the bound of \cite[Theorem~1]{BC}, whose main term is then $\nn{\alpha}\nn{\beta}M^{3/5+\eps}N^{7/20}$, only if $N>M^{27/28}$. The proof keeps the dependence on $N$ up to \eqref{eq:fivebudget}, and \eqref{eq:main_bound} results from replacing $N$ by $M$ there. For $N$ smaller than $M$ the bound \eqref{eq:fivebudget} is sharper than \eqref{eq:balancepair}, and a choice of $L$ adapted to $M$ and $N$, subject to \eqref{eq:blockrange}, would give a bound that depends on both lengths. This is not done here; it will be carried out in future research.

For $N=M$ and the choice \eqref{eq:amplifierchoice} of $L$, the five terms of \eqref{eq:fivebudget} are
\[
	\tfrac13M^{15/8}A^{1/4},\qquad
	\tfrac13M^{15/8}A^{3/4},\qquad
	3M^{15/8}A^{3/4},\qquad
	9M^{7/4}A^{1/2},\qquad
	3^{25/16}M^{231/128}A^{39/64}.
\]
The first comes from the diagonal terms and the terms controlled by the second moment without amplifier, the second from the terms with $\Delta=0$ (Section~\ref{sec:delta_zero}), the third from the terms estimated by Weil's bound in Section~\ref{sec:nongeneric}, and the last two from the spectral estimate in Section~\ref{sec:generic}, for the regular and for the exceptional spectrum. The exponent $15/16$ in \eqref{eq:main_bound} is therefore the result of the balance between the terms with $\Delta=0$, which decrease with $L$, and the terms of Section~\ref{sec:nongeneric}, which increase with $L$. The terms to which the spectral estimate is applied are smaller by a power of $M$.

\subsection{The exceptional spectrum}
\label{sec:exceptional_spectrum}
The exceptional eigenvalues enter through the factor $(1+T_0)^{2\theta(q)}$ in \eqref{eq:spectral_bound}, with the bound $\theta(q)\le7/64$, and they produce the last term of \eqref{eq:bound_H} and of \eqref{eq:fivebudget}. By the computation of Section~\ref{sec:critical_range}, this term is smaller than the third one by the factor $3^{9/16}M^{-9/128}A^{-9/64}$ when $N=M$. Consequently Selberg's eigenvalue conjecture, under which $\theta(q)=0$, would not improve \eqref{eq:main_bound} in the present argument, and neither would a stronger treatment of the exceptional spectrum in Lemma~\ref{lem:spectral} alone; for large sieve inequalities for exceptional Maass forms with arbitrary coefficients we refer to Blomer and Pascadi \cite[Theorem~1.6]{BP}. Among the two spectral terms the exceptional one is the larger. It is kept at this size by Lemma~\ref{lem:Delta_generic}: the level in Lemma~\ref{lem:progression} is divisible by the four amplifier primes, but its greatest common divisor with $\Delta$ is not.

\subsection{Larger numerators and fixed denominator factors}
\label{sec:larger_numerators}
The condition $A\le X^{1/10}$ in Theorem~\ref{thm:main} comes from the choice \eqref{eq:amplifierchoice} and from the condition $2A<L$ in \eqref{eq:blockrange}, which is used in Lemma~\ref{lem:Delta_generic}. Other choices of $L$ are possible. Let $\eta>0$ and $X^{1/10}\le A\le X^{1/4-\eta}$. Let $N\le M$, and let $M$ be large in terms of $\eta$. Then $L=3A$ satisfies $L\ge\max(L_0,8)$ and, for $N>2A$, the conditions \eqref{eq:blockrange}. In this case \eqref{eq:balancepair} gives $\Cw(\beta)\ll_\eps\nn{\beta}^2\nn{\nu}^2A^2M^{7/4+\eps}$ for sequences $\beta$ with pairwise unitary support, since $A^{-1/2}M^2\le A^2M^{7/4}$ for $A\ge M^{1/10}$. For $N\le2A$, and for the remaining values of $M$, the same bound follows from \eqref{eq:Cw_trivial}, with an implied constant that depends on $\eta$. Arguing as in Sections~\ref{sec:reduction_N_le_M}-\ref{sec:Cw_implies_main}, we obtain
\[
	|\Bcal(A,M,N)|\ll_{\eps,\eta}\nn{\alpha}\nn{\beta}\nn{\nu}\,A\,X^{7/8+\eps}
	\qquad\big(X^{1/10}\le A\le X^{1/4-\eta}\big).
\]
When $M$ and $N$ are of the same order of magnitude, this is smaller than the trivial bound in the whole range, and smaller than the bound of \cite[Theorem~1]{BC} for $A\le X^{3/26-\eta}$.

In some applications the denominator of the Kloosterman fraction has a fixed factor. This is the case in the proof of Theorem~\ref{thm:conv}, where Theorem~\ref{thm:main} is applied with the fixed factor absorbed in the variable $n$. Wright \cite[Theorem~2.1]{Wr1} proved a variant of the estimate of \cite{BC} in which a fixed factor of the denominator is taken into account. The corresponding variant of Theorem~\ref{thm:main} is not considered here; it will be taken up in future research.

\let\savedaddcontentsline\addcontentsline
\renewcommand{\addcontentsline}[3]{}
\section*{Acknowledgments}
\let\addcontentsline\savedaddcontentsline
The authors thank Alexandru Zaharescu for fruitful discussions. They are grateful to Alexandru Pascadi for pointing out an error in the first version of this paper.

The corrected argument has two new ingredients. The first is the use of the estimates of Deshouillers and Iwaniec in the form given by Pascadi \cite[Corollary~5.5]{P}, see Lemma~\ref{lem:spectral}. This was the authors' starting point for the repair. The second is the reduction to sequences with pairwise unitary support, in which an arbitrary sequence is written as an average of such sequences by choosing at random the exponent with which each prime is allowed to occur, see Lemmas~\ref{lem:random_weights} and~\ref{lem:random_average}. This idea was suggested by ChatGPT Astra~6 (OpenAI) in a conversation with one of the authors. By its own account, it drew on that author's work at the time in probability and combinatorics, which was part of the context available to it. The suggestion was not in a form that could be used directly, and the statements and proofs of Section~\ref{sec:argument} were worked out by the authors, who also chose the applications.

Claude Fable 5.2 (Anthropic) was used for adversarial auditing of the manuscript and assisted with its subsequent revision and exposition. The authors take full responsibility for the contents of this paper.

\end{document}